\documentclass{article}
\usepackage{a4wide,cite,latexsym,amsfonts,amssymb,exscale,enumerate}
\usepackage[centertags,sumlimits,intlimits,namelimits,reqno]{amsmath}
\usepackage{amsthm}
\usepackage{hyperref,picins}

\usepackage{pstricks}
\usepackage{pst-grad}
\usepackage{pst-plot}
\usepackage[tiling]{pst-fill}

\usepackage{pstcol}

\usepackage{graphics}
\usepackage[all, knot]{xy}
\xyoption{arc}

\usepackage{fancyheadings}
\newcommand{\cat}[1]{\ensuremath{\mbox{\bfseries {\upshape {#1}}}}}

\newcommand{\BOX}{\hbox {$\sqcap$ \kern -1em $\sqcup$}}

\newcommand{\Hom}{{\rm Hom}}
\newcommand{\Hess}{{\rm Hess}}

\renewcommand{\to}{\rightarrow}
\newcommand{\maps}{\colon}

\newcommand{\id}{{\rm id}}
\newcommand{\del}{\partial}

\newcommand{\Thfrob}{\cat{Th}(\cat{K-Frob})}

\newcommand{\R}{{\mathbb R}}

\newcommand{\A}{{\mathbb A}}

\newcommand{\twocob}{\cat{2Cob}^{{\rm ext}}}

\newcommand{\scs}{\scriptstyle}
\newcommand{\ten}{\otimes }

\theoremstyle{definition}
\newtheorem{thm}{Theorem}[section]
\newtheorem{cor}[thm]{Corollary}

\newtheorem{lem}[thm]{Lemma}

\newtheorem{prop}[thm]{Proposition}
\newtheorem{defn}[thm]{Definition}

\newcommand{\be}{\begin{equation}}
\newcommand{\ee}{\end{equation}}
\newcommand{\ba}{\begin{eqnarray}}
\newcommand{\ea}{\end{eqnarray}}
\newcommand{\ban}{\begin{eqnarray*}}
\newcommand{\ean}{\end{eqnarray*}}
\newcommand{\barr}{\begin{array}}
\newcommand{\earr}{\end{array}}

\newcommand{\adjunction}[4]{%
  \ensuremath{\xymatrix{%
    {#1} \ar@<4pt>[r]^{{#3}} \ar@{}[r]|{\scriptscriptstyle{\bot}} &%
    {#2} \ar@<4pt>[l]^{{#4}}%
  }}%
}

\newcommand{\radjunction}[4]{%
  \ensuremath{\xymatrix{%
    {#1} \ar@<4pt>[r]^{{#3}} \ar@{}[r]|{\scriptscriptstyle{\top}} &%
    {#2} \ar@<4pt>[l]^{{#4}}%
  }}%
}

\newcommand{\ambijunction}[4]{%
  \ensuremath{\xymatrix{%
    {#1} \ar@<4pt>[r]^{{#3}} \ar@{}[r]|{\scriptscriptstyle{\top \;\bot}} &%
    {#2} \ar@<4pt>[l]^{{#4}}%
  }}%
}

\newcommand{\padjunction}[4]{%
  \ensuremath{\xymatrix{%
    {#1} \ar@<4pt>[r]^{{#3}} \ar@{}[r]|{\scriptscriptstyle{\bot_p}} &%
    {#2} \ar@<4pt>[l]^{{#4}}%
  }}%
}

\newcommand{\pradjunction}[4]{%
  \ensuremath{\xymatrix{%
    {#1} \ar@<4pt>[r]^{{#3}} \ar@{}[r]|{\scriptscriptstyle{\top_p}} &%
    {#2} \ar@<4pt>[l]^{{#4}}%
  }}%
}

\newcommand{\pseudoambijunction}[4]{%
  \ensuremath{\xymatrix{%
    {#1} \ar@<4pt>[r]^{{#3}} \ar@{}[r]|{\scriptscriptstyle{\top_p \;\bot_p}} &%
    {#2} \ar@<4pt>[l]^{{#4}}%
  }}%
}

\SelectTips{cm}{}

\psset{linewidth=0.3pt,dimen=middle}
\psset{xunit=.70cm,yunit=0.70cm}

\newcommand{\hpeprint}[2]{%
  \href{http://www.arxiv.org/abs/#1}{\texttt{arxiv:#1#2}}}%
\newcommand{\hpmathsci}[1]{%
  \href{http://www.ams.org/mathscinet-getitem?mr=#1}{\texttt{MR #1}}}%

\numberwithin{equation}{section}

\def\emph#1{{\sl #1\/}}
\def\ie{{\sl i.e.\/}}

\def\etc{{\sl etc.\/}}
\def\cf{{\sl c.f.\/}}
\let\tilde=\widetilde

\let\phi=\varphi
\let\theta=\vartheta
\let\epsilon=\varepsilon

\usepackage{bbm}

\def\N{{\mathbbm N}}
\def\R{{\mathbbm R}}

\def\cal#1{\mathcal{#1}}%
\def\1{\mathbbm{1}}%
\def\nn{\notag}

\def\pushout#1#2#3{%
  \,\mbox{\raisebox{-.8ex}{$\scriptstyle #1$}}%
  \!\coprod\!
  \mbox{\raisebox{-.8ex}{$\scriptstyle #3$}}\,}%

\def\theauthor{}
\def\empty{}
\def\theaffiliation{}
\def\preprint#1{
  \thispagestyle{plain}
  \def\theauthor{#1}
  \ifx\theauthor\empty
  \else
    \begin{flushright}{\small #1\par}\end{flushright}
  \fi
  \begin{center}}
\def\title#1{
  {\LARGE #1\par}\vskip 1em}
\def\author#1{
  \ifx\theaffiliation\empty
  \else
    \par
  \fi
  \def\theauthor{#1}\def\theaffiliation{}}
\def\email#1{
  \vskip 1em{\large\theauthor\footnote{\small email: {\tt #1}}\par}\vskip .5em}
\def\affiliation#1{
  \ifx\theaffiliation\empty
    \def\theaffiliation{second}
  \else
    \par and\par
  \fi
  {\small\sl #1}}
\def\date#1{
  \vskip 1em{(#1)\par}\end{center}\vskip 2em}

\newgray{mygray}{0.65}
\newgray{mylightgray}{0.85}

\newcommand{\bigsaddle}{
  \pscustom[fillcolor=mygray,fillstyle=solid]{
    \psbezier(0.7,1.8)(0.7,2.0)(0.4,2.1)(0.2,2.1)
    \psline(0.2,2.1)(0.2,1.0)(0.7,1.0)(0.7,1.8)
  }
  \pscustom[fillcolor=mygray,fillstyle=solid]{
    \psbezier(2.0,1.8)(2.0,2.0)(2.5,2.1)(2.7,2.1)
    \psline(2.7,2.1)(2.7,0.6)(2.5,0.55)(2.0,0.6)(2.0,1.8)
  }
  \pscustom[fillcolor=mylightgray,fillstyle=solid]{
    \psline(0.0,1.5)(0.0,0.0)
    \psbezier(0.0,0.0)(0.5,0.2)(2.0,0.2)(2.5,0.0)
    \psline(2.5,0.0)(2.5,1.5)
    \psbezier(2.5,1.5)(2.2,1.5)(2.0,1.6)(2.0,1.8)
    \psbezier(2.0,1.8)(1.8,1.0)(0.9,1.0)(0.7,1.8)
    \psbezier(0.7,1.8)(0.7,1.6)(0.3,1.5)(0.0,1.5)
  }
  \psline[linestyle=dotted,dotsep=2pt](0.2,1.52)(0.2,0.6)
  \psbezier[linestyle=dotted,dotsep=2pt](0.2,0.6)(0.7,0.4)(2.2,0.4)(2.5,0.55)
  \psline[linestyle=dotted,dotsep=2pt](2.5,0.55)(2.7,0.6)
}

\newcommand{\bigcomposed}{
  \pscustom[fillcolor=mygray,fillstyle=solid]{
    \psbezier(0.7,1.8)(0.7,2.0)(0.4,2.1)(0.2,2.1)
    \psline(0.2,2.1)(0.2,1.0)(0.7,1.0)(0.7,1.8)
  }
  \pscustom[fillcolor=mygray,fillstyle=solid]{
    \psbezier(2.0,1.8)(2.0,2.0)(2.5,2.1)(2.7,2.1)
    \psbezier(2.7,2.1)(2.9,2.1)(3.3,2.0)(3.3,1.8)
    \psline(3.3,1.8)(3.3,1.0)(2.0,1.0)(2.0,1.8)
  }
  \pscustom[fillcolor=mylightgray,fillstyle=solid]{
    \psline(0.0,1.5)(0.0,0.0)
    \psbezier(0.0,0.0)(0.5,0.2)(2.0,0.0)(2.5,0.0)
    \psbezier(2.5,0.0)(2.7,0.0)(3.3,0.1)(3.3,0.3)
    \psline(3.3,0.3)(3.3,1.8)
    \psbezier(3.3,1.8)(3.3,1.6)(2.8,1.5)(2.6,1.5)
    \psbezier(2.6,1.5)(2.3,1.5)(2.0,1.6)(2.0,1.8)
    \psbezier(2.0,1.8)(1.8,1.0)(0.9,1.0)(0.7,1.8)
    \psbezier(0.7,1.8)(0.7,1.6)(0.3,1.5)(0.0,1.5)
  }
  \psline[linestyle=dotted,dotsep=2pt](0.2,1.52)(0.2,0.6)
  \psbezier[linestyle=dotted,dotsep=2pt](0.2,0.6)(0.7,0.4)(2.2,0.6)(2.5,0.6)
  \psbezier[linestyle=dotted,dotsep=2pt](2.5,0.6)(2.7,0.6)(3.3,0.4)(3.3,0.2)
}

\newcommand{\multl}{
  \pscustom[fillcolor=lightgray, fillstyle=solid]{
        \psbezier(1.5,2.5)(1.5,1.1)(.4,1.6)(.5,0)
        \psline(-0.5,0)
        \psbezier(-0.5,0)(-.4,1.6)(-1.5,1.1)(-1.5,2.5)
        \psline(-.5,2.5)
        \psbezier(-.5,2.5)(-.6,1.5)(0.6,1.5)(.5,2.5)
        \psline(1.5,2.5)
    }
}

\newcommand{\comultl}{
  \pscustom[fillcolor=lightgray, fillstyle=solid]{
        \psbezier(1.5,0)(1.5,1.4)(.4,.9)(.5,2.5)
        \psline(-0.5,2.5)
        \psbezier(-0.5,2.5)(-.4,.9)(-1.5,1.4)(-1.5,0)
        \psline(-.5,0)
        \psbezier(-.5,0)(-.6,1)(0.6,1)(.5,0)
        \psline(1.5,0)
    }
}

\newcommand{\widecomultl}{
  \pscustom[fillcolor=lightgray, fillstyle=solid]{
        \psbezier(2.5,0)(2.5,1.4)(0.4,.9)(.5,2.5)
        \psline(-.5,2.5)
        \psbezier(-.5,2.5)(-0.4,.9)(-2.5,1.4)(-2.5,0)
        \psline(-1.5,0)
        \psbezier(-1.5,0)(-1.6,1)(1.6,1)(1.5,0)
        \psline(2.5,0)
    }
}

\newcommand{\holel}{
      \pscustom[fillstyle=solid,
    fillcolor=lightgray]{
        \psbezier(1.5,2)(1.5,.6)(.4,1.1)(.5,0)
        \psline(.5,0)(-.5,0)
        \psbezier(-0.5,0)(-.4,1.1)(-1.5,.6)(-1.5,2)
        \psbezier(-1.5,2)(-1.5,3.4)(-.4,2.9)(-0.5,4)
        \psline(0.5,4)
        \psbezier(.5,4)(.4,2.9)(1.5,3.4)(1.5,2)
    }
 \pscustom[fillstyle=solid,
    fillcolor=white]{
        \psbezier(-.5,2)(-.55,1.2)(0.55,1.2)(.5,2)
        \psbezier(.5,2)(0.55,2.8)(-.55,2.8)(-.5,2)
    }
 }

\newcommand{\ctl}{
  \begin{psclip}{
    \pscustom{
        \psline(-.58,2)(-.58,0)
        \psline(-.58,0)(.42,0)
        \psline(.42,0)(.42,2)
        \psellipse(-.08,2)(.5,.2)
    }
  }
    \pspolygon[fillcolor=lightgray,fillstyle=gradient,
    gradbegin=lightgray,gradend=gray,gradmidpoint=1,gradangle=110](-.58,0)(-.58,2.4)(.42,2.4)(.42,0)(-.58,0)
 \end{psclip}
 \pscustom[fillcolor=lightgray,fillstyle=gradient,
        gradbegin=white, gradend=gray,gradmidpoint=0,gradangle=88]{
    \psline(-.58,2)(-.58,0)
    \psbezier(-.58,0)(-.48,.5)(-.48,.7)(-.08,1)
    \psbezier(-.08,1)(.32,.7)(.32,.5)(.42,0)
    \psellipse(-.08,2)(.5,.2)
 }
 \psellipse[fillcolor=lightgray,fillstyle=gradient,
        gradbegin=lightgray, gradend=gray,gradmidpoint=1,gradangle=110](-.08,2)(.5,.2)
}

\newcommand{\twistl}{
  \pscustom[fillcolor=lightgray, fillstyle=solid]{
    \psbezier(.5,2.5)(.6,2.0)(-.4,2.1)(-.5,1.25)
    \psline(.5,1.25)
    \psbezier(.5,1.25)(.4,2.1)(-.6,2.0)(-.5,2.5)
    \psline(.5,2.5)
 }
   \pscustom[fillcolor=lightgray, fillstyle=solid]{
    \psbezier(-.5,1.25)(-.4,.4)(.6,.5)(.5,0)
    \psline(-.5,0)
    \psbezier(-.5,0)(-.6,.5)(.4,.4)(.5,1.25)
    \psline(-.5,1.25)
 }
    \begin{psclip}{
    \pscustom{
        \psbezier(.5,2.5)(.6,2.0)(-.4,2.1)(-.5,1.25)
        \psbezier(-.5,1.25)(-.4,.4)(.6,.5)(.5,0)
    }
  }
  \pscustom[fillcolor=darkgray, fillstyle=solid]{
   \psbezier(-.5,2.5)(-.6,2.0)(.4,2.1)(.5,1.25)
    \psbezier(.5,1.25)(.4,.4)(-.6,.5)(-.5,0)}
 \end{psclip}
 }

\newcommand{\holew}{
    \pspolygon[fillcolor=lightgray,fillstyle=gradient,
    gradbegin=lightgray,gradend=gray,gradmidpoint=1,gradangle=110](.58,0)(.58,4)(-.42,4)(-.42,0)(.58,0)
   \pscustom[fillcolor=lightgray,fillstyle=gradient,
        gradbegin=white, gradend=gray,gradmidpoint=0,gradangle=88]{
    \psline(.58,4)(.58,2)
    \psbezier(.58,2)(.48,2.5)(.48,2.7)(.08,3)
    \psbezier(.08,3)(-.32,2.7)(-.32,2.5)(-.42,2)
    \psline(-.42,4)
    \psellipse(.08,4)(.5,.2)
 }
 \psellipse[fillcolor=lightgray,fillstyle=gradient,
        gradbegin=lightgray, gradend=gray,gradmidpoint=1,gradangle=110](.08,4)(.5,.2)
 \pscustom[fillcolor=lightgray,fillstyle=gradient,
        gradbegin=white, gradend=gray,gradmidpoint=0,gradangle=88]{
    \psline(.58,0)(.58,2)
    \psbezier(.58,2)(.48,1.5)(.48,1.3)(.08,1)
    \psbezier(.08,1)(-.32,1.3)(-.32,1.5)(-.42,2)
    \psline(-.42,0)
    \psbezier(-.42,0)(-.32,-.25)(.48,-.25)(.58,0)
 }
 \begin{psclip}{
 \pspolygon[linestyle=none](.58,0)(.58,.3)(-.42,.3)(-.42,0)(.58,0)
 }
 \psellipse[linestyle=dotted](.08,0)(.5,0.2)
 \end{psclip}
}

\newcommand{\ltc}{
    \pspolygon[fillcolor=lightgray,fillstyle=gradient,
    gradbegin=lightgray,gradend=gray,gradmidpoint=1,gradangle=60](.58,2)(.58,.4)(-.42,.4)(-.42,2)(.58,2)
 \pscustom[fillcolor=lightgray,fillstyle=gradient,
        gradbegin=white, gradend=gray,gradmidpoint=0,gradangle=88]{
    \psline(.58,0)(.58,2)
    \psbezier(.58,2)(.48,1.5)(.48,1.3)(.08,1)
    \psbezier(.08,1)(-.32,1.3)(-.32,1.5)(-.42,2)
    \psline(-.42,0)
    \psbezier(-.42,0)(-.32,-.25)(.48,-.25)(.58,0)
 }
 \begin{psclip}{
 \pspolygon[linestyle=none](.58,0)(.58,.3)(-.42,.3)(-.42,0)(.58,0)
 }
 \psellipse[linestyle=dotted](.08,0)(.5,0.2)
 \end{psclip}
}

\newcommand{\ltcnew}{
    \pspolygon[fillcolor=lightgray,fillstyle=gradient,
    gradbegin=lightgray,gradend=gray,gradmidpoint=1,gradangle=60](.5,2.5)(.5,.9)(-.5,.9)(-.5,2.5)(.5,2.5)
 \pscustom[fillcolor=lightgray,fillstyle=gradient,
        gradbegin=white, gradend=gray,gradmidpoint=0,gradangle=88]{
    \psline(.5,0)(.5,2.5)
    \psbezier(.5,2.5)(.4,2)(.4,1.8)(0,1.5)
    \psbezier(0,1.5)(-.4,1.8)(-.4,2)(-.5,2.5)
    \psline(-.5,0)
    \psbezier(-.5,0)(-.4,-.25)(.4,-.25)(.5,0)
 }
 \begin{psclip}{
 \pspolygon[linestyle=none](.5,0)(.5,.3)(-.5,.3)(-.5,0)(.5,0)
 }
 \psellipse[linestyle=dotted](0,0)(.5,0.2)
 \end{psclip}
}

 \newcommand{\birthl}{
 \pscustom[fillcolor=lightgray, fillstyle=solid]{
        \psbezier(-.5,0)(-.5,.9)(0.5,.9)(.5,0)
        \psline(-.5,0)
    }
 }

  \newcommand{\deathl}{
 \pscustom[fillcolor=lightgray, fillstyle=solid]{
        \psbezier(-.5,0)(-.5,-.9)(0.5,-.9)(.5,0)
        \psline(-.5,0)    }
 }

\newcommand{\zagl}{
   \pscustom[fillcolor=lightgray, fillstyle=solid]{
        \psbezier(1.5,0)(1.6,2)(-1.6,2)(-1.5,0)
        \psline(-.5,0)
        \psbezier(-.5,0)(-.6,.8)(0.6,.8)(.5,0)
        \psline(1.5,0)
    }
}

\newcommand{\zigl}{
       \pscustom[fillcolor=lightgray, fillstyle=solid]{
        \psbezier(1.5,2)(1.6,0)(-1.6,0)(-1.5,2)
        \psline(-.5,2)
        \psbezier(-.5,2)(-.6,1.1)(0.6,1.1)(.5,2)
        \psline(1.5,2)
    }
}

\newcommand{\identl}{
    \pspolygon[fillcolor=lightgray,fillstyle=solid](-.5,0)(.5,0)(.5,2.5)(-.5,2.5)(-.5,0)
}

\newcommand{\medidentl}{
    \pspolygon[fillcolor=lightgray,fillstyle=solid](-.5,0)(.5,0)(.5,2)(-.5,2)(-.5,0)
}

\newcommand{\smallidentl}{
     \pscustom[fillcolor=lightgray, fillstyle=solid]{
        \pspolygon(-.5,0)(.5,0)(.5,1)(-.5,1)(-.5,0)
    }
}

\newcommand{\Xl}[1]{%
  \pspolygon[fillcolor=lightgray,fillstyle=solid](-.5,0)(.5,0)(.5,1.5)(-.5,1.5)(-.5,0)
   \pspolygon[fillcolor=white,fillstyle=solid](-.8,0)(.8,0)(.8,1)(-.8,1)(-.8,0)
  \rput(0,.5){$#1$}
}

\newcommand{\Generalf}{
     \rput(-.8,1.5){\smallidentc}
      \rput(1.2,1.5){\smallidentc}
      \rput(-2.8,-.8){\smallidentc}
      \rput(3.2,-.8){\smallidentc}
 \pspolygon[fillcolor=lightgray,fillstyle=solid](-3.45,1)(-2.45,1)(-2.45,2.5)(-3.45,2.5)(-3.45,0)
 \pspolygon[fillcolor=lightgray,fillstyle=solid](3.45,1)(2.45,1)(2.45,2.5)(3.45,2.5)(3.45,0)
  \pspolygon[fillcolor=lightgray,fillstyle=solid](1.5,-1)(.5,-1)(.5,0)(1.5,0)(1.5,0)
  \pspolygon[fillcolor=lightgray,fillstyle=solid](-1.5,-1)(-.5,-1)(-.5,0)(-1.5,0)(-1.5,0)
  \pspolygon[fillcolor=white,fillstyle=solid](-3.8,0)(3.8,0)(3.8,1.5)(-3.8,1.5)(-3.8,0)
   \rput(0,.75){$\scs [f]$}
   }

\newcommand{\GeneralNF}{
     \rput(-.8,1.5){\smallidentc}
      \rput(1.2,1.5){\smallidentc}
      \rput(-2.8,-.8){\smallidentc}
      \rput(3.2,-.8){\smallidentc}
 \pspolygon[fillcolor=lightgray,fillstyle=solid](-3.45,1)(-2.45,1)(-2.45,2.5)(-3.45,2.5)(-3.45,0)
 \pspolygon[fillcolor=lightgray,fillstyle=solid](3.45,1)(2.45,1)(2.45,2.5)(3.45,2.5)(3.45,0)
  \pspolygon[fillcolor=lightgray,fillstyle=solid](1.5,-1)(.5,-1)(.5,0)(1.5,0)(1.5,0)
  \pspolygon[fillcolor=lightgray,fillstyle=solid](-1.5,-1)(-.5,-1)(-.5,0)(-1.5,0)(-1.5,0)
  \pspolygon[fillcolor=white,fillstyle=solid](-3.8,0)(3.8,0)(3.8,1.5)(-3.8,1.5)(-3.8,0)
   \rput(0,.75){$\scs [NF(f)]$}
   }

 \newcommand{\sigmaOne}{
         \rput(3.2,1.5){\smallidentc}
      \rput(1.2,1.5){\smallidentc}
      \rput(-.8,-.8){\smallidentc}
      \rput(1.2,-.8){\smallidentc}
  \pspolygon[fillcolor=lightgray,fillstyle=solid](-3.45,1)(-2.45,1)(-2.45,2.5)(-3.45,2.5)(-3.45,0)
  \pspolygon[fillcolor=lightgray,fillstyle=solid](-1.5,1)(-.5,1)(-.5,2.5)(-1.5,2.5)(-1.5,1)
   \pspolygon[fillcolor=lightgray,fillstyle=solid](-3.45,-1)(-2.45,-1)(-2.45,0)(-3.45,0)(-3.45,-1)
    \pspolygon[fillcolor=lightgray,fillstyle=solid](3.45,-1)(2.45,-1)(2.45,0)(3.45,0)(3.45,-1)
   \pspolygon[fillcolor=white,fillstyle=solid](-3.8,0)(3.8,0)(3.8,1.5)(-3.8,1.5)(-3.8,0)
   \rput(0,.7){$\scs \sigma_1^{-1}$}
   }

 \newcommand{\normalf}{
           \rput(-.8,-.8){\identc}
        \rput(1.2,-.8){\identc}
        \rput(-2.8,-.8){\identc}
        \rput(3.2,-.8){\identc}
    \pspolygon[fillcolor=lightgray,fillstyle=solid](-3.48,3)(-2.48,3)(-2.48,5.5)(-3.48,5.5)(-3.48,3)
  \pspolygon[fillcolor=lightgray,fillstyle=solid](-1.5,3)(-.5,3)(-.5,5.5)(-1.5,5.5)(-1.5,3)
    \pspolygon[fillcolor=lightgray,fillstyle=solid](3.45,3)(2.45,3)(2.45,5.5)(3.45,5.5)(3.45,3)
  \pspolygon[fillcolor=lightgray,fillstyle=solid](1.46,3)(.46,3)(.46,5.5)(1.46,5.5)(1.46,3)
   \pspolygon[fillcolor=white,fillstyle=solid](-3.91,1)(3.89,1)(3.89,4)(-3.91,4)(-3.91,1)
   \rput(.1,2.5){$\scs {\rm NF_\mathrm{O\rightarrow C}}(\Lambda([f]))$}
   }

   \newcommand{\sigmaOneINV}{
         \rput(-.8,1.5){\smallidentc}
      \rput(1.2,1.5){\smallidentc}
      \rput(3.28,-.8){\smallidentc}
      \rput(1.28,-.8){\smallidentc}
  \pspolygon[fillcolor=lightgray,fillstyle=solid](-3.45,1)(-2.45,1)(-2.45,2.5)(-3.45,2.5)(-3.45,0)
  \pspolygon[fillcolor=lightgray,fillstyle=solid](3.45,1)(2.45,1)(2.45,2.5)(3.45,2.5)(3.45,0)
   \pspolygon[fillcolor=lightgray,fillstyle=solid](-3.45,-1)(-2.45,-1)(-2.45,0)(-3.45,0)(-3.45,-1)
    \pspolygon[fillcolor=lightgray,fillstyle=solid](-1.45,-1)(-.45,-1)(-.45,0)(-1.45,0)(-1.45,-1)
   \pspolygon[fillcolor=white,fillstyle=solid](-3.8,0)(3.8,0)(3.8,1.5)(-3.8,1.5)(-3.8,0)
   \rput(0,.75){$\scs \sigma_1$}
   }

   \newcommand{\sigmaOnePrime}{
         \rput(-.8,1.5){\smallidentc}
      \rput(1.2,1.5){\smallidentc}
      \rput(3.28,-.8){\smallidentc}
      \rput(1.28,-.8){\smallidentc}
  \pspolygon[fillcolor=lightgray,fillstyle=solid](-3.45,1)(-2.45,1)(-2.45,2.5)(-3.45,2.5)(-3.45,0)
  \pspolygon[fillcolor=lightgray,fillstyle=solid](3.45,1)(2.45,1)(2.45,2.5)(3.45,2.5)(3.45,0)
   \pspolygon[fillcolor=lightgray,fillstyle=solid](-3.45,-1)(-2.45,-1)(-2.45,0)(-3.45,0)(-3.45,-1)
    \pspolygon[fillcolor=lightgray,fillstyle=solid](-1.45,-1)(-.45,-1)(-.45,0)(-1.45,0)(-1.45,-1)
   \pspolygon[fillcolor=white,fillstyle=solid](-3.8,0)(3.8,0)(3.8,1.5)(-3.8,1.5)(-3.8,0)
   \rput(0,.75){$\scs \sigma_1'$}
   }

   \newcommand{\sigmaTwoINV}{
         \rput(3.2,1.5){\smallidentc}
      \rput(1.2,1.5){\smallidentc}
      \rput(3.2,-.8){\smallidentc}
      \rput(-2.8,-.8){\smallidentc}
  \pspolygon[fillcolor=lightgray,fillstyle=solid](-3.45,1)(-2.45,1)(-2.45,2.5)(-3.45,2.5)(-3.45,0)
  \pspolygon[fillcolor=lightgray,fillstyle=solid](-1.45,1)(-.45,1)(-.45,2.5)(-1.45,2.5)(-1.45,0)
    \pspolygon[fillcolor=lightgray,fillstyle=solid](1.45,-1)(.45,-1)(.45,0)(1.45,0)(1.45,-1)
    \pspolygon[fillcolor=lightgray,fillstyle=solid](-1.45,-1)(-.45,-1)(-.45,0)(-1.45,0)(-1.45,-1)
   \pspolygon[fillcolor=white,fillstyle=solid](-3.8,0)(3.8,0)(3.8,1.5)(-3.8,1.5)(-3.8,0)
   \rput(0,.7){$\scs \sigma_2^{-1}$}
   }

      \newcommand{\sigmaTwoINVPrime}{
         \rput(3.2,1.5){\smallidentc}
      \rput(1.2,1.5){\smallidentc}
      \rput(3.2,-.8){\smallidentc}
      \rput(-2.8,-.8){\smallidentc}
  \pspolygon[fillcolor=lightgray,fillstyle=solid](-3.45,1)(-2.45,1)(-2.45,2.5)(-3.45,2.5)(-3.45,0)
  \pspolygon[fillcolor=lightgray,fillstyle=solid](-1.45,1)(-.45,1)(-.45,2.5)(-1.45,2.5)(-1.45,0)
    \pspolygon[fillcolor=lightgray,fillstyle=solid](1.45,-1)(.45,-1)(.45,0)(1.45,0)(1.45,-1)
    \pspolygon[fillcolor=lightgray,fillstyle=solid](-1.45,-1)(-.45,-1)(-.45,0)(-1.45,0)(-1.45,-1)
   \pspolygon[fillcolor=white,fillstyle=solid](-3.8,0)(3.8,0)(3.8,1.5)(-3.8,1.5)(-3.8,0)
   \rput(0,.7){$\scs \sigma_2'$}
   }

  \newcommand{\sigmaTwo}{
          \rput(-2.8,1.5){\smallidentc}
      \rput(3.2,1.5){\smallidentc}
      \rput(3.2,-.8){\smallidentc}
      \rput(1.2,-.8){\smallidentc}
  \pspolygon[fillcolor=lightgray,fillstyle=solid](-1.5,1)(-.5,1)(-.5,2.5)(-1.5,2.5)(-1.5,0)
   \pspolygon[fillcolor=lightgray,fillstyle=solid](1.5,1)(.5,1)(.5,2.5)(1.5,2.5)(1.5,0)
   \pspolygon[fillcolor=lightgray,fillstyle=solid](-3.45,-1)(-2.45,-1)(-2.45,0)(-3.45,0)(-3.45,-1)
  \pspolygon[fillcolor=lightgray,fillstyle=solid](-1.5,-1)(-.5,-1)(-.5,0)(-1.5,0)(-1.5,0)
   \pspolygon[fillcolor=white,fillstyle=solid](-3.8,0)(3.8,0)(3.8,1.5)(-3.8,1.5)(-3.8,0)
   \rput(0,.75){$\scs \sigma_2$}
  }

\newcommand{\BottomX}[1]{%
  \pspolygon[fillcolor=white,fillstyle=solid](-1,0)(1,0)(1,1)(-1,1)(-1,0)
  \psline[linewidth=1pt](-2,0)(2,0)
  \rput(0,.5){$#1$}
}

\newcommand{\Lambdaf}{
\rput(10.2,-3.3){\identc}
  \rput(8.2,-3.3){\identc}
   \rput(12.2,-3.3){\identc}
  \rput(14.2,-3.3){\identc}
  \rput(12.2,-1){\identc}
  \rput(14.2,-1){\identc}
  \rput(12.2,1.3){\identc}
  \rput(14.2,1.3){\identc}
  \rput(12.2,3.8){\identc}
  \rput(14.2,3.8){\identc}
   \rput(12.2,6.3){\identc}
  \rput(14.2,6.3){\identc}
   \rput(7,0){\sigmaTwo}
  \rput(7,3.2){\Generalf}
  \rput(7,6.4){\sigmaOne}
  \rput(11.2,8.8){\zagc}
  \rput(11.2,8.8){\bigzagc}
   \rput(3.05,-3){\zigl}
   \rput(3,-1){\bigzigl}
   \rput(2.05,-1){\identl}
   \rput(.05,-1){\identl}
   \rput(2.05,1.3){\identl}
   \rput(.05,1.3){\identl}
   \rput(2.05,3.8){\identl}
   \rput(.05,3.8){\identl}
   \rput(2.05,6.3){\identl}
   \rput(.05,6.3){\identl}
   \rput(2.05,8.8){\identl}
   \rput(.05,8.8){\identl}
   \rput(4.05,8.8){\identl}
   \rput(6.0,8.8){\identl}
   }

\newcommand{\LambdaINVi}{
    \rput(3.78,-.5){\sigmaTwoINVPrime}
   \rput(12.2,1.2){\zigc}
   \rput(12.2,1.2){\bigzigc}
    \rput(5,2){\identc}
    \rput(7,2){\identc}
    \rput(15.4,3.2){\medidentc}
    \rput(13.4,3.2){\medidentc}
    \rput(15.4,5.2){\identc}
    \rput(13.4,5.2){\identc}
    \rput(3.8,6.5){\zagl}
    \rput(3.8,6.5){\bigzagl}
    \rput(2.8,4.5){\medidentl}
    \rput(0.8,4.5){\medidentl}
    \rput(2.8,2){\identl}
    \rput(0.8,2){\identl}
    \rput(11.15,5.2){\identl}
    \rput(9.15,5.2){\identl}
 \rput(8,4){\GeneralG}
 \rput(12.1,8.5){\sigmaOnePrime}
 }

   \newcommand{\LambdaLambda}{
   \rput(8,-13.2){\sigmaTwoINV}
 \rput(20,-5.3){\zigc}
 \rput(20,-5.3){\bigzigc}
 \rput(17.2,-3.3){\identc}
   \rput(9.2,-10.8){\curverightc}
   \rput(11.2,-10.8){\curverightc}
   \rput(10.2,-8.3){\curverightc}
   \rput(12.2,-8.3){\curverightc}
  \rput(11.2,-5.8){\curverightc}
   \rput(13.2,-5.8){\curverightc}
  \rput(12.2,-3.3){\curverightc}
   \rput(14.2,-3.3){\curverightc}
  \rput(19.2,-3.3){\identc}
  \rput(17.2,-1){\identc}
  \rput(19.2,-1){\identc}
  \rput(17.2,1.3){\identc}
  \rput(19.2,1.3){\identc}
  \rput(17.2,3.8){\identc}
  \rput(19.2,3.8){\identc}
   \rput(17.2,6.3){\identc}
  \rput(19.2,6.3){\identc}
     \rput(23.2,-3.3){\identc}
  \rput(21.2,-3.3){\identc}
  \rput(23.2,-1){\identc}
  \rput(21.2,-1){\identc}
  \rput(23.2,1.3){\identc}
  \rput(21.2,1.3){\identc}
  \rput(23.2,3.8){\identc}
  \rput(21.2,3.8){\identc}
   \rput(23.2,6.3){\identc}
  \rput(21.2,6.3){\identc}
  \rput(21.2,8.8){\curveleftc}
  \rput(23.2,8.8){\curveleftc}
  \rput(20.2,11.3){\curveleftc}
  \rput(22.2,11.3){\curveleftc}
  \rput(19.2,13.8){\curveleftc}
  \rput(21.2,13.8){\curveleftc}
  \rput(18.2,16.3){\curveleftc}
  \rput(20.2,16.3){\curveleftc}
   \rput(12,0){\sigmaTwo}
  \rput(12,3.2){\Generalf}
  \rput(12,6.4){\sigmaOne}
  \rput(16.2,8.8){\zagc}
  \rput(16.2,8.8){\bigzagc}
   \rput(8.05,-3){\zigl}
   \rput(8,-1){\bigzigl}
   \rput(7.05,-1){\identl}
   \rput(5.05,-1){\identl}
   \rput(7.05,1.3){\identl}
   \rput(5.05,1.3){\identl}
   \rput(7.05,3.8){\identl}
   \rput(5.05,3.8){\identl}
   \rput(7.05,6.3){\identl}
   \rput(5.05,6.3){\identl}
   \rput(7.05,8.8){\identl}
   \rput(5.05,8.8){\identl}
   \rput(9.05,8.8){\curverightl}
   \rput(11,8.8){\curverightl}
   \rput(10.05,11.3){\curverightl}
   \rput(12,11.3){\curverightl}
   \rput(11.05,13.8){\curverightl}
   \rput(13,13.8){\curverightl}
   \rput(12.05,16.3){\curverightl}
   \rput(14,16.3){\curverightl}
   \rput(7.05,-11){\curveleftl}
   \rput(5.05,-11){\curveleftl}
   \rput(4.05,-8.5){\curveleftl}
   \rput(6.05,-8.5){\curveleftl}
   \rput(3.05,-6){\curveleftl}
   \rput(5.05,-6){\curveleftl}
   \rput(2.05,-3.5){\curveleftl}
   \rput(4.05,-3.5){\curveleftl}
   \rput(1.05,-1){\identl}
   \rput(3.05,-1){\identl}
   \rput(3.05,1.3){\identl}
   \rput(1.05,1.3){\identl}
   \rput(3.05,3.8){\identl}
   \rput(1.05,3.8){\identl}
   \rput(3.05,6.3){\identl}
   \rput(1.05,6.3){\identl}
   \rput(3.05,8.8){\identl}
   \rput(1.05,8.8){\identl}
   \rput(4.0,11.2){\zagl}
   \rput(4.0,11.2){\bigzagl}
   \rput(15.9,19.6){\sigmaOneINV}
   }

  \newcommand{\Xc}[1]{ %
     \pscustom[fillcolor=lightgray,fillstyle=gradient,
        gradbegin=white, gradend=gray,gradmidpoint=0,gradangle=88]{
        \psline(-.5,1.6)(-.5,0)
        \psline(.5,0)
        \psline(.5,1.6)
        \psline(-.5,1.6)
    }
\psellipse[fillcolor=lightgray,fillstyle=gradient,
        gradbegin=lightgray, gradend=gray,gradmidpoint=1,gradangle=110](0,1.6)(.5,.2)
 \pspolygon[fillcolor=white,fillstyle=solid](-.8,0)(.8,0)(.8,1)(-.8,1)(-.8,0)
  \rput(0,.5){$#1$}
}

\newcommand{\crossl}{
    \pspolygon[fillcolor=lightgray,fillstyle=solid](-1.5,0)(-.5,0)(1.5,2.5)(.5,2.5)(-1.5,0)
    \pspolygon[fillcolor=lightgray,fillstyle=solid](1.5,0)(.5,0)(-1.5,2.5)(-.5,2.5)(1.5,0)
}

\newcommand{\ucrossl}{
    \pspolygon[fillcolor=lightgray,fillstyle=solid](1.5,0)(.5,0)(-1.5,2.5)(-.5,2.5)(1.5,0)
    \pspolygon[fillcolor=lightgray,fillstyle=solid](-1.5,0)(-.5,0)(1.5,2.5)(.5,2.5)(-1.5,0)
}

\newcommand{\curverightl}{
  \pscustom[fillcolor=lightgray, fillstyle=solid]{
        \psbezier(1.5,2.5)(1.5,1.5)(.4,1.3)(.5,0)
        \psline(-0.5,0)
        \psbezier(-.5,0)(-.6,1.3)(.5,1.5)(.5,2.5)
        \psline(1.5,2.5)
    }
}

\newcommand{\curveleftl}{
  \pscustom[fillcolor=lightgray, fillstyle=solid]{
        \psbezier(-1.5,2.5)(-1.5,1.5)(-.4,1.3)(-.5,0)
        \psline(0.5,0)
        \psbezier(.5,0)(.6,1.3)(-.5,1.5)(-.5,2.5)
        \psline(-1.5,2.5)
    }
}

\newcommand{\topbits}{
  \rput(6,0){\ltc}
  \rput(5.8,2){\multl}
  \rput(4.8,4.5){\multl}
  \rput(2.5,9){$\ddots$}
  \rput(1.5,10){\multl}
}

\newcommand{\medbits}{
  \rput(5.8,0){\multc}
  \rput(4.8,2.5){\multc}
  \rput(2,7){$\ddots$}
  \rput(0,8){\hugemultc}
}

\newcommand{\multc}{
      \pscustom[fillstyle=gradient,
    gradbegin=white, gradend=gray,gradmidpoint=0,gradangle=70]{
        \psbezier(1.5,2.5)(1.5,1.1)(.4,1.6)(.5,0)
        \psbezier(.5,0)(.4,-.25)(-.4,-.25)(-.5,0)
        \psbezier(-0.5,0)(-.4,1.6)(-1.5,1.1)(-1.5,2.5)
        \psline(-.5,2.5)
        \psbezier(-.5,2.5)(-.6,1.5)(0.6,1.5)(.5,2.5)
        \psline(1.5,2.5)
    }
    \psellipse[fillcolor=lightgray,fillstyle=gradient,
        gradbegin=lightgray, gradend=gray,gradmidpoint=1,gradangle=110](-1,2.5)(.5,.2)
    \psellipse[fillcolor=lightgray,fillstyle=gradient,
        gradbegin=lightgray, gradend=gray,gradmidpoint=1,gradangle=110](1,2.5)(.5,.2)
     \begin{psclip}{
 \pspolygon[linestyle=none](.5,0)(.5,.3)(-.5,.3)(-.5,0)(.5,0)
 }
 \psellipse[linestyle=dotted](0,0)(.5,0.2)
 \end{psclip}
 }

\newcommand{\widemultc}{
  \pscustom[fillstyle=gradient,gradbegin=white,gradend=gray,gradmidpoint=0,gradangle=70]{
     \psbezier(2,2.5)(2,1.1)(.4,1.6)(.5,0)
     \psbezier(.5,0)(.4,-.25)(-.4,-.25)(-.5,0)
     \psbezier(-0.5,0)(-.4,1.6)(-2,1.1)(-2,2.5)
     \psline(-1,2.5)
     \psbezier(-1,2.5)(-.6,1.5)(0.6,1.5)(1,2.5)
     \psline(2,2.5)
  }
  \psellipse[fillcolor=lightgray,fillstyle=gradient,gradbegin=lightgray,gradend=gray,gradmidpoint=1,gradangle=110](-1.5,2.5)(.5,.2)
  \psellipse[fillcolor=lightgray,fillstyle=gradient,gradbegin=lightgray,gradend=gray,gradmidpoint=1,gradangle=110](1.5,2.5)(.5,.2)
  \begin{psclip}{
    \pspolygon[linestyle=none](.5,0)(.5,.3)(-.5,.3)(-.5,0)(.5,0)
   }
   \psellipse[linestyle=dotted](0,0)(.5,0.2)
   \end{psclip}
}

\newcommand{\hugemultc}{
    \pscustom[fillstyle=gradient,gradbegin=white,gradend=gray,gradmidpoint=0,gradangle=70]{
     \psbezier(3,2.5)(2.2,.9)(.4,1.2)(.5,0)
     \psbezier(.5,0)(.4,-.25)(-.4,-.25)(-.5,0)
     \psbezier(-0.5,0)(-.4,1.2)(-2.2,.9)(-3,2.5)
     \psline(-2,2.5)
     \psbezier(-2,2.5)(-1.6,1.2)(1.6,1.2)(2,2.5)
     \psline(3,2.5)
  }
  \psellipse[fillcolor=lightgray,fillstyle=gradient,gradbegin=lightgray,
  gradend=gray,gradmidpoint=1,gradangle=110](-2.5,2.5)(.5,.2)
  \psellipse[fillcolor=lightgray,fillstyle=gradient,gradbegin=lightgray,
  gradend=gray,gradmidpoint=1,gradangle=110](2.5,2.5)(.5,.2)
  \begin{psclip}{
    \pspolygon[linestyle=none](.5,0)(.5,.3)(-.5,.3)(-.5,0)(.5,0)
   }
   \psellipse[linestyle=dotted](0,0)(.5,0.2)
   \end{psclip}
   }

\newcommand{\comultc}{
  \pscustom[fillstyle=gradient,
    gradbegin=white, gradend=gray,gradmidpoint=0,gradangle=110]{
        \psbezier(1.5,0)(1.5,1.4)(.4,.9)(.5,2.5)
        \psline(-0.5,2.5)
        \psbezier(-0.5,2.5)(-.4,.9)(-1.5,1.4)(-1.5,0)
        \psbezier(-1.5,0)(-1.4,-.25)(-.6,-.25)(-.5,0)
        \psbezier(-.5,0)(-.6,1)(0.6,1)(.5,0)
        \psbezier(.5,0)(.6,-.25)(1.4,-.25)(1.5,0)
    }
  \psellipse[fillcolor=lightgray,fillstyle=gradient,
        gradbegin=lightgray, gradend=gray,gradmidpoint=1,gradangle=110](0,2.5)(.5,.2)
\begin{psclip}{
 \pspolygon[linestyle=none](1.5,0)(1.5,.3)(-1.5,.3)(-1.5,0)(1.5,0)
 }
 \psellipse[linestyle=dotted](1,0)(.5,0.2)
 \psellipse[linestyle=dotted](-1,0)(.5,0.2)
 \end{psclip}
 }

 \newcommand{\holec}{
      \pscustom[fillstyle=gradient,
    gradbegin=white, gradend=gray,gradmidpoint=0,gradangle=70]{
        \psbezier(1.5,2)(1.5,.6)(.4,1.1)(.5,0)
        \psbezier(.5,0)(.4,-.25)(-.4,-.25)(-.5,0)
        \psbezier(-0.5,0)(-.4,1.1)(-1.5,.6)(-1.5,2)
        \psline(-.5,2)
        \psbezier(-.5,2)(-.55,1.2)(0.55,1.2)(.5,2)
        \psline(1.5,2)
    }
     \begin{psclip}{
 \pspolygon[linestyle=none](.5,0)(.5,.3)(-.5,.3)(-.5,0)(.5,0)
 }
 \psellipse[linestyle=dotted](0,0)(.5,0.2)
 \end{psclip}
 \pscustom[fillstyle=gradient,
    gradbegin=white, gradend=gray,gradmidpoint=0,gradangle=110,linestyle=none]{
        \psbezier(1.5,2)(1.5,3.4)(.4,2.9)(.5,4)
        \psline(-0.5,4)
        \psbezier(-0.5,4)(-.4,2.9)(-1.5,3.4)(-1.5,2)
        \psbezier(-1.5,2)(-1.4,1.75)(-.6,1.75)(-.5,2)
        \psbezier(-.5,2)(-.55,2.8)(0.55,2.8)(.5,2)
        \psbezier(.5,2)(.6,1.75)(1.4,1.75)(1.5,2)
    }
 \psbezier(1.5,2)(1.5,3.4)(.4,2.9)(.5,4)
 \psbezier(-0.5,4)(-.4,2.9)(-1.5,3.4)(-1.5,2)
        \psbezier(-.5,2)(-.55,2.8)(0.55,2.8)(.5,2)
  \psellipse[fillcolor=lightgray,fillstyle=gradient,
        gradbegin=lightgray, gradend=gray,gradmidpoint=1,gradangle=110](0,4)(.5,.2)
  }

\newcommand{\birthc}{
 \pscustom[fillstyle=gradient,
    gradbegin=white, gradend=gray,gradmidpoint=0,gradangle=110]{
        \psbezier(-.5,0)(-.5,.9)(0.5,.9)(.5,0)
        \psbezier(.5,0)(.4,-.25)(-.4,-.25)(-.5,0)
    }
 \begin{psclip}{
 \pspolygon[linestyle=none](.5,0)(.5,.3)(-.5,.3)(-.5,0)(.5,0)
 }
 \psellipse[linestyle=dotted](0,0)(.5,0.2)
 \end{psclip}
 }

\newcommand{\deathc}{
 \pscustom[fillstyle=gradient,
    gradbegin=white, gradend=gray,gradmidpoint=0,gradangle=70]{
        \psbezier(-.5,1)(-.5,.1)(0.5,.1)(.5,1)
        \psline(-.5,1)
 }
  \psellipse[fillcolor=lightgray,fillstyle=gradient,
        gradbegin=lightgray, gradend=gray,gradmidpoint=1,gradangle=110](0,1)(.5,.2)
 }

\newcommand{\zagc}{
   \pscustom[fillstyle=gradient,
    gradbegin=white, gradend=gray,gradmidpoint=0,gradangle=110]{
        \psbezier(1.5,0)(1.6,2)(-1.6,2)(-1.5,0)
        \psbezier(-1.5,0)(-1.4,-.25)(-.6,-.25)(-.5,0)
        \psbezier(-.5,0)(-.6,.8)(0.6,.8)(.5,0)
        \psbezier(.5,0)(.6,-.25)(1.4,-.25)(1.5,0)
    }
  \begin{psclip}{
 \pspolygon[linestyle=none](1.5,0)(1.5,.3)(-1.5,.3)(-1.5,0)(1.5,0)
 }
 \psellipse[linestyle=dotted](1,0)(.5,0.2)
 \psellipse[linestyle=dotted](-1,0)(.5,0.2)
 \end{psclip}
 }

 \newcommand{\bigzagc}{
  \pscustom[fillstyle=gradient,
    gradbegin=white, gradend=gray,gradmidpoint=0,gradangle=70]{
        \psbezier(3.5,0)(3.6,4)(-3.6,4)(-3.5,0)
        \psbezier(-3.5,0)(-3.4,-.25)(-2.6,-.25)(-2.5,0)
        \psbezier(-2.5,0)(-2.6,2.9)(2.6,2.9)(2.5,0)
        \psbezier(2.5,0)(2.6,-.25)(3.4,-.25)(3.5,0)
    }
 \begin{psclip}{
 \pspolygon[linestyle=none](3.5,0)(3.5,.3)(-3.5,.3)(-3.5,0)(3.5,0)
 }
 \psellipse[linestyle=dotted](3,0)(.5,0.2)
 \psellipse[linestyle=dotted](-3,0)(.5,0.2)
 \end{psclip}
 }

 \newcommand{\bigzagl}{
  \pscustom[fillstyle=solid,fillcolor=lightgray]{
        \psbezier(3.5,0)(3.6,4)(-3.6,4)(-3.5,0)
        \psline(-2.5,0)
        \psbezier(-2.5,0)(-2.6,2.9)(2.6,2.9)(2.5,0)
        \psline(3.5,0)
    }
 }

  \newcommand{\bigzigl}{
  \pscustom[fillstyle=solid,fillcolor=lightgray]{
        \psbezier(3.5,0)(3.6,-4)(-3.6,-4)(-3.5,0)
        \psline(-2.5,0)
        \psbezier(-2.5,0)(-2.6,-2.9)(2.6,-2.9)(2.5,0)
        \psline(3.5,0)
    }
 }

 \newcommand{\GeneralG}{
         \rput(-1,-.8){\smallidentc}
      \rput(1.4,-.8){\smallidentc}
      \rput(3.4,-.8){\smallidentc}
      \rput(-3,-.8){\smallidentc}
  \pspolygon[fillcolor=lightgray,fillstyle=solid](-3.65,1)(-2.65,1)(-2.65,2.5)(-3.65,2.5)(-3.65,0)
  \pspolygon[fillcolor=lightgray,fillstyle=solid](3.65,1)(2.65,1)(2.65,2.5)(3.65,2.5)(3.65,0)
    \pspolygon[fillcolor=lightgray,fillstyle=solid](-1.65,1)(-.65,1)(-.65,2.5)(-1.65,2.5)(-1.65,0)
  \pspolygon[fillcolor=lightgray,fillstyle=solid](1.65,1)(.65,1)(.65,2.5)(1.65,2.5)(1.65,0)
   \pspolygon[fillcolor=white,fillstyle=solid](-4,0)(4,0)(4,1.5)(-4,1.5)(-4,0)
   \rput(0,.75){$\scs [g]$}
   \rput(0,2.05){$;$}
   \rput(0,-.7){$;$}
   }

\newcommand{\zigc}{
       \pscustom[fillstyle=gradient,
    gradbegin=white, gradend=gray,gradmidpoint=0,gradangle=70]{
        \psbezier(1.5,2)(1.6,0)(-1.6,0)(-1.5,2)
        \psline(-.5,2)
        \psbezier(-.5,2)(-.6,1.2)(0.6,1.2)(.5,2)
        \psline(1.5,2)
    }
 \psellipse[fillcolor=lightgray,fillstyle=gradient,
        gradbegin=lightgray, gradend=gray,gradmidpoint=1,gradangle=110](1,2)(.5,.2)
        \psellipse[fillcolor=lightgray,fillstyle=gradient,
        gradbegin=lightgray, gradend=gray,gradmidpoint=1,gradangle=110](-1,2)(.5,.2)
}

\newcommand{\bigzigc}{
       \pscustom[fillstyle=gradient,
    gradbegin=white, gradend=gray,gradmidpoint=0,gradangle=70]{
        \psbezier(3.5,2)(3.6,-2)(-3.6,-2)(-3.5,2)
        \psline(-2.5,2)
        \psbezier(-2.5,2)(-2.6,-.9)(2.6,-.9)(2.5,2)
        \psline(3.5,2)
    }
 \psellipse[fillcolor=lightgray,fillstyle=gradient,
        gradbegin=lightgray, gradend=gray,gradmidpoint=1,gradangle=110](3,2)(.5,.2)
        \psellipse[fillcolor=lightgray,fillstyle=gradient,
        gradbegin=lightgray, gradend=gray,gradmidpoint=1,gradangle=110](-3,2)(.5,.2)
 }

\newcommand{\identc}{
 \pscustom[fillcolor=lightgray,fillstyle=gradient,
        gradbegin=white, gradend=gray,gradmidpoint=0,gradangle=88]{
 \psline(.5,0)(.5,2.5)
 \psline(-.5,2.5)
 \psline(-.5,0)
 \psbezier(-.5,0)(-.4,-.25)(.4,-.25)(.5,0)
 }
\psellipse[fillcolor=lightgray,fillstyle=gradient,
        gradbegin=lightgray, gradend=gray,gradmidpoint=1,gradangle=110](0,2.5)(.5,.2)
 \begin{psclip}{
 \pspolygon[linestyle=none](.5,0)(.5,.3)(-.5,.3)(-.5,0)(.5,0)
 }
 \psellipse[linestyle=dotted](0,0)(.5,0.2)
 \end{psclip}
 }

  \newcommand{\medidentc}{
     \pscustom[fillcolor=lightgray,fillstyle=gradient,
        gradbegin=white, gradend=gray,gradmidpoint=0,gradangle=88]{
        \psline(-.5,2)(-.5,0)
        \psbezier(-.5,0)(-.4,-.25)(.4,-.25)(.5,0)
        \psline(.5,2)
        \psline(-.5,2)
    }
\psellipse[fillcolor=lightgray,fillstyle=gradient,
        gradbegin=lightgray, gradend=gray,gradmidpoint=1,gradangle=110](0,2)(.5,.2)
 \begin{psclip}{
 \pspolygon[linestyle=none](.5,0)(.5,.3)(-.5,.3)(-.5,0)(.5,0)
 }
 \psellipse[linestyle=dotted](0,0)(.5,0.2)
 \end{psclip}
}

 \newcommand{\smallidentc}{
     \pscustom[fillcolor=lightgray,fillstyle=gradient,
        gradbegin=white, gradend=gray,gradmidpoint=0,gradangle=88]{
        \psline(-.5,1)(-.5,0)
        \psbezier(-.5,0)(-.4,-.25)(.4,-.25)(.5,0)
        \psline(.5,1)
        \psline(-.5,1)
    }
\psellipse[fillcolor=lightgray,fillstyle=gradient,
        gradbegin=lightgray, gradend=gray,gradmidpoint=1,gradangle=110](0,1)(.5,.2)
 \begin{psclip}{
 \pspolygon[linestyle=none](.5,0)(.5,.3)(-.5,.3)(-.5,0)(.5,0)
 }
 \psellipse[linestyle=dotted](0,0)(.5,0.2)
 \end{psclip}
}

\newcommand{\crossc}{
   \pscustom[fillcolor=lightgray,fillstyle=gradient,
        gradbegin=white, gradend=gray,gradmidpoint=0,gradangle=125]{
 \psline(-.5,0)(1.5,2.5)
 \psline(.5,2.5)
 \psline(-1.5,0)
 \psbezier(-1.5,0)(-1.4,-.25)(-.6,-.25)(-.5,0)
 }
 \pscustom[fillcolor=lightgray,fillstyle=gradient,
        gradbegin=white, gradend=gray,gradmidpoint=0,gradangle=125]{
 \psline(.5,0)(-1.5,2.5)
 \psline(-.5,2.5)
 \psline(1.5,0)
 \psbezier(1.5,0)(1.4,-.25)(.6,-.25)(.5,0)
 }
\psellipse[fillcolor=lightgray,fillstyle=gradient,
        gradbegin=lightgray, gradend=gray,gradmidpoint=1,gradangle=110](-1,2.5)(.5,.2)
\psellipse[fillcolor=lightgray,fillstyle=gradient,
        gradbegin=lightgray, gradend=gray,gradmidpoint=1,gradangle=70](1,2.5)(.5,.2)
 \psline[linestyle=dotted](1.5,2.5)(-.5,0)
 \psline[linestyle=dotted](.5,2.5)(-1.5,0)
 \begin{psclip}{
 \pspolygon[linestyle=none](1.5,0)(1.5,.3)(-1.5,.3)(-1.5,0)(1.5,0)
 }
 \psellipse[linestyle=dotted](1,0)(.5,0.2)
 \psellipse[linestyle=dotted](-1,0)(.5,0.2)
 \end{psclip}
 }

\newcommand{\crossmixlc}{
   \pscustom[fillcolor=lightgray,fillstyle=gradient,
        gradbegin=white, gradend=gray,gradmidpoint=0,gradangle=125]{
 \psline(-.5,0)(1.5,2.5)
 \psline(.5,2.5)
 \psline(-1.5,0)
 \psbezier(-1.5,0)(-1.4,-.25)(-.6,-.25)(-.5,0)
 }
 \pscustom[fillcolor=lightgray,fillstyle=solid]{
 \psline(.5,0)(-1.5,2.5)
 \psline(-.5,2.5)
 \psline(1.5,0)
 \psline(.5,0)
 }
\psellipse[fillcolor=lightgray,fillstyle=gradient,
        gradbegin=lightgray, gradend=gray,gradmidpoint=1,gradangle=110](1,2.5)(.5,.2)
 \begin{psclip}{
 \pspolygon[linestyle=none](1.5,0)(1.5,.3)(-1.5,.3)(-1.5,0)(1.5,0)
 }
 \psellipse[linestyle=dotted](-1,0)(.5,0.2)
 \end{psclip}
 }

\newcommand{\crossmixcl}{
   \pscustom[fillcolor=lightgray,fillstyle=solid]{
 \psline(-.5,0)(1.5,2.5)
 \psline(.5,2.5)
 \psline(-1.5,0)
 \psline(-.5,0)
 }
 \pscustom[fillcolor=lightgray,fillstyle=gradient,
        gradbegin=white, gradend=gray,gradmidpoint=0,gradangle=125]{
 \psline(.5,0)(-1.5,2.5)
 \psline(-.5,2.5)
 \psline(1.5,0)
 \psbezier(1.5,0)(1.4,-.25)(.6,-.25)(.5,0)
 }
\psellipse[fillcolor=lightgray,fillstyle=gradient,
        gradbegin=lightgray, gradend=gray,gradmidpoint=1,gradangle=70](-1,2.5)(.5,.2)
 \begin{psclip}{
 \pspolygon[linestyle=none](1.5,0)(1.5,.3)(0,.3)(0,0)(1.5,0)
 }
 \psellipse[linestyle=dotted](1,0)(.5,0.2)
 \end{psclip}
 }

\newcommand{\curverightc}{
  \pscustom[fillstyle=gradient,
    gradbegin=white, gradend=gray,gradmidpoint=0,gradangle=65]{
        \psbezier(1.5,2.5)(1.5,1.5)(.4,1.3)(.5,0)
        \psbezier(.5,0)(.4,-.25)(-.4,-.25)(-.5,0)
        \psbezier(-.5,0)(-.6,1.3)(.5,1.5)(.5,2.5)
        \psline(1.5,2.5)
    }
     \psellipse[fillcolor=lightgray,fillstyle=gradient,
        gradbegin=lightgray, gradend=gray,gradmidpoint=1,gradangle=110](1,2.5)(.5,.2) \begin{psclip}{
 \pspolygon[linestyle=none](.5,0)(.5,.3)(-.5,.3)(-.5,0)(.5,0)
 }
 \psellipse[linestyle=dotted](0,0)(.5,0.2)
 \end{psclip}
}

\newcommand{\curveleftc}{
  \pscustom[fillstyle=gradient,
    gradbegin=white, gradend=gray,gradmidpoint=0,gradangle=115]{
        \psbezier(-1.5,2.5)(-1.5,1.5)(-.4,1.3)(-.5,0)
        \psbezier(-.5,0)(-.4,-.25)(.4,-.25)(.5,0)
        \psbezier(.5,0)(.6,1.3)(-.5,1.5)(-.5,2.5)
        \psline(-1.5,2.5)
    }
    \psellipse[fillcolor=lightgray,fillstyle=gradient,
        gradbegin=lightgray, gradend=gray,gradmidpoint=1,gradangle=110](-1,2.5)(.5,.2)
 \begin{psclip}{
 \pspolygon[linestyle=none](.5,0)(.5,.3)(-.5,.3)(-.5,0)(.5,0)
 }
 \psellipse[linestyle=dotted](0,0)(.5,0.2)
 \end{psclip}
}

\newcommand{\ctlab}[1]{%
  \rput(0,0){\ctl}
  \rput(-0.7,-0.25){$\scriptstyle #1$}
  \rput(0.7,-0.25){$\scriptstyle #1$}
}
\newcommand{\ltcab}[1]{%
  \rput(0,0){\ltc}
  \rput(-0.7,2.25){$\scriptstyle #1$}
  \rput(0.7,2.25){$\scriptstyle #1$}
}

\begin{document}
%

\preprint{AEI-2005-153\\ DAMTP-2005-80}

\title{Open-closed strings: Two-dimensional extended TQFTs and Frobenius algebras}

\author{Aaron D.\ Lauda}
\email{A.Lauda@dpmms.cam.ac.uk}
\affiliation{Department of Pure Mathematics and Mathematical Statistics,\\
  University of Cambridge, Cambridge CB3 0WB, United Kingdom}

\author{Hendryk Pfeiffer}
\email{pfeiffer@math.ubc.ca}
\affiliation{Department of Mathematics, University of British Columbia,\\
  1984 Mathematics Road, Vancouver, BC, V2T 1Z2, Canada}

\date{January 6, 2008}

%
\begin{abstract}
%

We study a special sort of $2$-dimensional extended Topological
Quantum Field Theories (TQFTs). These are defined on open-closed
cobordisms by which we mean smooth compact oriented $2$-manifolds with
corners that have a particular global structure in order to model the
smooth topology of open and closed string worldsheets. We show that
the category of open-closed TQFTs is equivalent to the category of
knowledgeable Frobenius algebras. A knowledgeable Frobenius algebra
$(A,C,\imath,\imath^\ast)$ consists of a symmetric Frobenius algebra
$A$, a commutative Frobenius algebra $C$, and an algebra homomorphism
$\imath\colon C\to A$ with dual $\imath^\ast\colon A\to C$, subject to
some conditions. This result is achieved by providing a description of
the category of open-closed cobordisms in terms of generators and the
well-known Moore--Segal relations. In order to prove the sufficiency
of our relations, we provide a normal form for such cobordisms which
is characterized by topological invariants. Starting from an arbitrary
such cobordism, we construct a sequence of moves (generalized handle
slides and handle cancellations) which transforms the given cobordism
into the normal form. Using the generators and relations description
of the category of open-closed cobordisms, we show that it is
equivalent to the symmetric monoidal category freely generated by a
knowledgeable Frobenius algebra. Our formalism is then generalized to
the context of open-closed cobordisms with labeled free boundary
components, \ie\ to open-closed string worldsheets with D-brane labels
at their free boundaries.

\end{abstract}

\noindent
\begin{small}
\end{small}

%
\section{Introduction}
%

Motivated by open string theory, boundary conformal field theory, and
extended topological quantum field theory, open-closed cobordisms have
been a topic of considerable interest to mathematicians and
physicists. By open-closed cobordisms we mean the morphisms of a
category $\twocob$ whose objects are compact oriented smooth
$1$-manifolds, \ie\ free unions of circles $S^1$ and unit intervals
$I=[0,1]$. The morphisms are certain compact oriented smooth
2-manifolds with corners. The corners of such a manifold $f$ are
required to coincide with the boundary points $\del I$ of the
intervals. The boundary of $f$ viewed as a topological manifold, minus
the corners, consists of components that are either `black' or
'coloured'. Each corner is required to separate a black component from
a coloured one. The black part of the boundary coincides with the
union of the source and the target objects. Two such manifolds with
corners are considered equivalent if they are related by an
orientation preserving diffeomorphism which restricts to the identity
on the black part of the boundary.  An example of such an open-closed
cobordism is depicted here\footnote{In order to get a feeling for
these diagrams, the reader might wish to verify that this cobordism is
diffeomorphic to the one depicted in Figure~1 of~\cite{BCR}.},
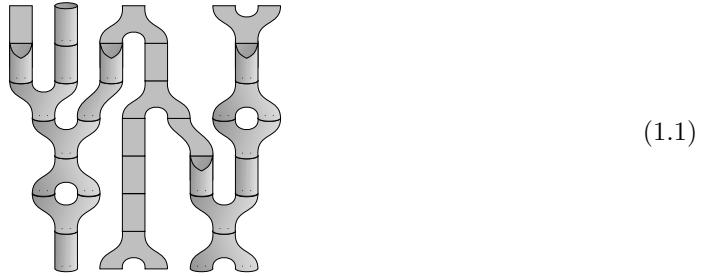
\begin{equation}
\label{eq_figure}
\begin{aligned}
\psset{xunit=.3cm,yunit=.2cm}
\begin{pspicture}[0.5](10,18)
  \rput(  2, 0){\identc}
  \rput(4.8, 0){\comultl}
  \rput(  9, 0){\comultc}
  \rput(  2, 2.5){\multc}
  \rput(4.8, 2.5){\identl}
  \rput(  9, 2.5){\multc}
  \rput(  2, 5){\comultc}
  \rput(4.8, 5){\identl}
  \rput(  8, 5){\ltcnew}
  \rput( 10, 5){\identc}
  \rput(  2, 7.5){\multc}
  \rput(4.8, 7.5){\identl}
  \rput(7.8, 7.5){\curveleftl}
  \rput( 10, 7.5){\multc}
  \rput(  1,10){\multc}
  \rput(  3,10){\curverightc}
  \rput(5.8,10){\comultl}
  \rput( 10,10){\comultc}
  \rput(  0,12.5){\ltcnew}
  \rput(  2,12.5){\identc}
  \rput(  4,12.5){\ltcnew}
  \rput(5.8,12.5){\identl}
  \rput( 10,12.5){\ltcnew}
  \rput(-0.2,15){\identl}
  \rput(  2,15){\identc}
  \rput(4.8,15){\comultl}
  \rput(9.8,15){\multl}
\end{pspicture}
\end{aligned}
\end{equation}
where the boundaries at the top and at the bottom of the diagram
are the black ones. In Section~\ref{sect_cob}, we present a formal
definition which includes some additional technical properties.
Gluing such cobordisms along their black boundaries, \ie\ putting
the building blocks of~\eqref{eq_figure} on top of each other, is
the composition of morphisms. The free union of manifolds, \ie\
putting the building blocks of~\eqref{eq_figure} next to each
other, provides $\twocob$ with the structure of a strict symmetric
monoidal category.

Open-closed cobordisms can be seen as a generalization of the
conventional $2$-dimensional cobordism category $\cat{2Cob}$. The
objects of this symmetric monoidal category are compact oriented
smooth $1$-manifolds \emph{without boundary}; the morphisms are
compact oriented smooth cobordisms between them, modulo
orientation-preserving diffeomorphisms that restrict to the identity
on the boundary.

The study of open-closed cobordisms plays an important role in
conformal field theory if one is interested in boundary conditions,
and open-closed cobordisms have a natural string theoretic
interpretation. The intervals in the black boundaries are interpreted
as open strings, the circles as closed strings, and the open-closed
cobordisms as string worldsheets. Here we consider only the underlying
smooth manifolds, but not any additional conformal or complex
structure. Additional labels at the coloured boundaries are
interpreted as D-branes or boundary conditions on the open strings.

An open-closed Topological Quantum Field Theory (TQFT), which we
formally define in Section~\ref{sect_tqft} below, is a symmetric
monoidal functor $\twocob\to\cal{C}$ into a symmetric monoidal
category $\cal{C}$. If $\cal{C}$ is the category of vector spaces over
a fixed field $k$, then the open-closed TQFT assigns vector spaces to
the $1$-manifolds $I$ and $S^1$, it assigns tensor products to free
unions of these manifolds, and $k$-linear maps to open-closed
cobordisms.

Such an open-closed TQFT can be seen as an extension of the notion of
a $2$-dimensional TQFT~\cite{Atiyah} which is a symmetric monoidal
functor $\cat{2Cob}\to\cal{C}$. We refer to this conventional notion
of $2$-dimensional TQFT as a \emph{closed TQFT} and to the morphisms
of $\cat{2Cob}$ as \emph{closed cobordisms}. For the classic results
on $2$-dimensional closed TQFTs, we recommend the original
works~\cite{Dij1,Abrams1,Sawin95} and the book~\cite{Kock}.

The most powerful results on closed TQFTs crucially depend on results
from Morse theory. Morse theory provides a generators and relations
description of the category $\cat{2Cob}$. First, any compact cobordism
$\Sigma$ can be obtained by gluing a finite number of elementary
cobordisms along their boundaries. In order to see this, one chooses a
Morse function $f\colon\Sigma\to\R$ such that all critical points have
distinct critical values and considers the pre-images
$f^{-1}([x_0-\epsilon,x_0+\epsilon])\subseteq\Sigma$ of intervals that
contain precisely one critical value $x_0\in\R$.  Each such pre-image
is the free union of one of the elementary cobordisms,
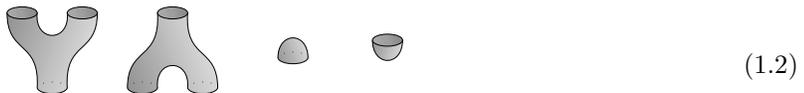
\begin{equation}
\label{eq_closedgen}
  \psset{xunit=.4cm,yunit=.4cm}
  \begin{pspicture}[.2](12,2.5)
    \rput(0,0){\multc}
    \rput(4,0){\comultc}
    \rput(8,1){\birthc}
    \rput(11,.6){\deathc}
  \end{pspicture}
\end{equation}
with zero or more cylinders over $S^1$. The different elementary
cobordisms~\eqref{eq_closedgen} are precisely the Morse data that
characterize the critical points, and the way they are glued
corresponds to the handle decomposition associated with $f$. The
Morse data of~\eqref{eq_closedgen} provide the \emph{generators}
for the morphisms of $\cat{2Cob}$. Our diagrams, for
example~\eqref{eq_figure}, are organized in such a way that the
vertical axis of the drawing plane serves as a Morse function, and
the cobordisms are composed of building blocks that contain at
most one critical point.

Second, given two Morse functions $f_1,f_2\colon\Sigma\to\R$, the
handle decompositions associated with $f_1$ and $f_2$ are related by a
finite sequence of \emph{moves}, \ie\ handle slides and handle
cancellations. This means that there are diffeomorphisms such as,
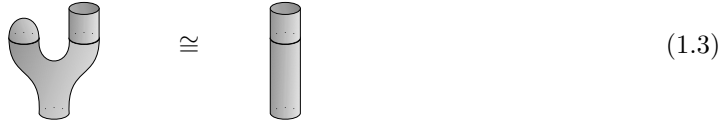
\begin{equation}
\label{eq_samplerel}
  \psset{xunit=.4cm,yunit=.4cm}
  \begin{pspicture}[0.5](4,4)
     \rput(2,0){\multc}
     \rput(3,2.5){\smallidentc}
     \rput(1,2.5){\birthc}
   \end{pspicture}
   \qquad\cong\qquad
   \begin{pspicture}[0.5](2,4)
     \rput(1,0){\identc}
     \rput(1,2.5){\smallidentc}
   \end{pspicture}
\end{equation}
which provide us with the \emph{relations} of $\cat{2Cob}$. When we
explicitly construct the diffeomorphism that relates two handle
decompositions of some manifold, we call these diffeomorphisms
\emph{moves}. The example~\eqref{eq_samplerel} corresponds to a
cancellation of a $1$-handle and a $2$-handle. Below is an example
of sliding a 1-handle past another 1-handle.
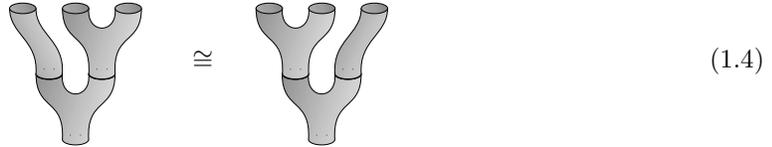
\begin{equation}
    \label{eq_samplerel2}
 \psset{xunit=.35cm,yunit=.35cm}
\begin{pspicture}[0.5](4,5.5)
  \rput(2,0){\multc}
  \rput(3,2.5){\multc}
  \rput(1,2.5){\curveleftc}
\end{pspicture}
\qquad\cong\qquad
\begin{pspicture}[0.5](4,5.5)
  \rput(2,0){\multc}
  \rput(1,2.5){\multc}
  \rput(3,2.5){\curverightc}
\end{pspicture}
\end{equation}

\parpic[r]{
$
\psset{xunit=.15cm,yunit=.15cm}
\begin{pspicture}[.2](5,25)
  \rput(1,0){\comultc}
  \rput(4,0){\curveleftc}
  \rput(6,0){\curveleftc}
  \rput(2,2.5){\comultc}
  \rput(5,2.5){\curveleftc}
  \rput(3,5){\comultc}
  \rput(2.65,7.5){\holec}
  \rput(2.65,11.5){\holec}
  \rput(2.65,15.5){\holec}
  \rput(3,19.5){\multc}
  \rput(2,22){\multc}
  \rput(4,22){\curverightc}
\end{pspicture}
$
}
Whereas it is not too difficult to construct by brute-force a set of
diffeomorphisms between manifolds such as those
in~\eqref{eq_samplerel} and~\eqref{eq_samplerel2}, \ie\ to show that
a set of relations is necessary, it is much harder to show that they
are also sufficient, \ie\ that any two handle decompositions are
related by a finite sequence of moves such as~\eqref{eq_samplerel}
and~\eqref{eq_samplerel2}. In order to establish this result, one
strategy is to prove that there exists a \emph{normal form} for the
morphisms of $\cat{2Cob}$ which is characterized by topological
invariants, and then to show that the relations suffice in order to
transform an arbitrary handle decomposition into this normal form.
The normal form for closed cobordisms is determined by the number of
incoming and outgoing boundary components together with the genus.
The example to the right shows the normal form of a closed cobordism
with three incoming boundary components, four outgoing boundary
components, and genus three. For closed cobordisms, the normal form
and proof of the sufficiency of the relations is done in detail
in~\cite{Carmody,Abrams1,Kock}.

Rather than employing the normal form, one could try to make precise,
in the context of manifolds with corners, the obvious Morse theoretic
ideas that underly the Moore--Segal relations. The advantage of the
normal form is, however, that it results in a constructive proof which
delivers all relevant diffeomorphisms in terms of sequences of
relations being applied to the relevant handlebodies (up to smooth
isotopies of the attaching sets).

In order to describe open-closed cobordisms using generators and
relations, one would need a generalization of Morse theory for
manifolds with corners. Such a generalization of Morse theory can be
used in order to find the generators of $\twocob$, and brute force can
be used to establish the necessity of certain relations.  However, we
are not aware of any abstract theorem that would guarantee the
sufficiency of these relations.

The first main result of this article is a normal form for open-closed
cobordisms with an inductive proof that the relations suffice in order
to transform any handle decomposition into the normal form. As a
consequence, for any two diffeomorphic open-closed cobordisms whose
handle decompositions are given, we explicitly construct a
diffeomorphism relating the two by constructing the corresponding
sequence of moves.

The description of $\twocob$ in terms of generators and relations
has emerged over the last couple of years from consistency
conditions in boundary conformal field theory, going back to the
work of Cardy and Lewellen~\cite{CL,Lew}, Lazaroiu~\cite{lar2}, and
Moore and Segal, see, for example~\cite{MS,MS2}, and these results
have been known to the experts for some time. More recently, the
unoriented case has also been considered by Alexeevski and
Natanzon~\cite{AN}. The first result of this article which we claim
is new, is the normal form and our inductive proof that the
relations are sufficient.

This, in turn, implies the following result in Morse theory for our
sort of compact 2-manifolds with corners which has so far not been
available by other means: The handle decompositions associated with
any two Morse functions on the same manifold are related by a finite
sequence of handle slides and handle cancellations. In particular,
these include the handle cancellation depicted on the right
of~\eqref{ziphomo} below in which a critical point on the coloured
boundary `eats up' a critical point of the interior.

Once a description of $\twocob$ in terms of generators and relations
is available, it is possible to find an algebraic characterization
for the symmetric monoidal category of open-closed TQFTs. Whereas
the category of closed TQFTs is equivalent as a symmetric monoidal
category to the category of commutative Frobenius
algebras~\cite{Abrams1}, we prove that the category of open-closed
TQFTs is equivalent as a symmetric monoidal category to the category
of \emph{knowledgeable Frobenius algebras}. We define knowledgeable
Frobenius algebras in Section~\ref{sect_knowfrob} precisely for this
purpose. A knowledgeable Frobenius algebra
$(A,C,\imath,\imath^\ast)$ consists of a symmetric Frobenius algebra
$A$, a commutative Frobenius algebra $C$, and an algebra
homomorphism $\imath\colon C\to A$ with dual $\imath^\ast\colon A\to
C$, subject to some conditions. This is the second main result of
the present article. The structure that emerges is consistent with
the results of Moore and Segal~\cite{MS,MS2}.

The algebraic structures relevant to boundary conformal field theory
have been studied by Fuchs and Schweigert~\cite{FS}. In a series of
papers, for example~\cite{FRS}, Fuchs, Runkel, and Schweigert study
Frobenius algebra objects in ribbon categories.  Topologically, this
corresponds to a situation in which the surfaces are embedded in some
$3$-manifold and studied up to ambient isotopy. In the present
article, in contrast, we consider Frobenius algebra objects in a
symmetric monoidal category, and our 2-manifolds are considered
equivalent as soon as they are diffeomorphic (as abstract manifolds)
relative to the boundary.

Various extensions of open-closed topological field theories have also
been studied. Baas, Cohen, and Ram\'{\i}rez have extended the
symmetric monoidal category of open-closed cobordisms to a symmetric
monoidal 2-category whose 2-morphisms are certain diffeomorphisms of
the open-closed cobordisms \cite{BCR}.  This work extends the work of
Tillmann who defined a symmetric monoidal 2-category extending the
closed cobordism category \cite{Till3}. She used this 2-category to
introduce an infinite loop space structure on the plus construction of
the stable mapping class group of closed cobordisms
\cite{Till2}. Using a similar construction to Tillmann's, Baas,
Cohen, and Ram\'{\i}rez have defined an infinite loop space
structure on the plus construction of the stable mapping class group
of open-closed cobordisms, showing that infinite loop space
structures are a valuable tool in studying the mapping class group.

Another extension of open-closed TQFT comes from open-closed
Topological {\it Conformal} Field Theory (TCFT). It was shown by
Costello~\cite{Costello} that the category of open Topological
Conformal Field Theories is homotopy equivalent to the category of
certain $A_{\infty}$ categories with extra structure. Ignoring the
conformal structure, or equivalently taking $H_0$ of the $\Hom$ spaces
in the corresponding category, reduces this to the case of Topological
Quantum Field Theory. Costello associates to a given open TCFT an
open-closed TCFT where the homology of the closed states is the
Hochschild homology of the $A_{\infty}$ category describing the open
states. This work is also useful for providing generators and
relations for the category of open Riemann surfaces and, when
truncated, this result also agrees with the characterization of open
cobordisms and their diffeomorphisms up to isotopy given in
\cite{Lau2} where a smaller list of generators and relations is
given. In the present article, we aim directly for an explicit
description of the category of open-closed cobordisms.

The present article is structured as follows: In
Section~\ref{sect_knowfrob}, we define the notion of a
knowledgeable Frobenius algebra and introduce the symmetric
monoidal category $\cat{K-Frob}(\cal{C})$ of knowledgeable
Frobenius algebras in a symmetric monoidal category $\cal{C}$. We
provide an abstract description in terms of generators and
relations of this category by defining a category $\Thfrob$,
called the \emph{theory of knowledgeable Frobenius algebras}, and
by showing that the category of symmetric monoidal functors and
monoidal natural transformations $\Thfrob\to\cal{C}$ is equivalent
as a symmetric monoidal category to $\cat{K-Frob}(\cal{C})$. In
Section~\ref{sect_cob}, we introduce the category $\twocob$ of
open-closed cobordisms. We present a normal form for such
cobordisms and characterize the category in terms of generators
and relations. In Section~\ref{sect_tqft}, we define open-closed
TQFTs as symmetric monoidal functors $\twocob\to\cal{C}$ into some
symmetric monoidal category $\cal{C}$. We show that the category
$\twocob$ is equivalent as a symmetric monoidal category to
$\Thfrob$ which in turn implies that the category of open-closed
TQFTs in $\cal{C}$ is equivalent as a symmetric monoidal category
to the category of knowledgeable Frobenius algebras
$\cat{K-Frob}(\cal{C})$. In Section~\ref{sect_labels}, we
generalize our results to the case of labeled free boundaries.
Section~\ref{sect_conclusion} contains a summary and an outlook on
open problems.

%
\section{Knowledgeable Frobenius algebras}
%
\label{sect_knowfrob}

\subsection{Definitions}

In this section, we define the notion of a knowledgeable Frobenius
algebra. We consider these Frobenius algebras not only in the
symmetric monoidal category $\cat{Vect}_k$ of vector spaces over
some fixed field $k$, but in any generic symmetric monoidal
category. Other examples include the symmetric monoidal categories
of Abelian groups, graded-vector spaces, and chain complexes.

\begin{defn}
Let $(\cal{C},\otimes,\1,\alpha,\lambda,\rho,\tau)$ be a symmetric
monoidal category with tensor product
$\otimes\colon\cal{C}\times\cal{C}\to\cal{C}$, unit object
$\1\in|\cal{C}|$, associativity constraint $\alpha_{X,Y,Z}\colon
(X\otimes Y)\otimes Z\to X\otimes(Y\otimes Z)$, left- and
right-unit constraints $\lambda_X\colon \1\otimes X\to X$ and
$\rho_X\colon X\otimes\1\to X$, and the symmetric braiding
$\tau_{X,Y}\colon X\otimes Y\to Y\otimes X$, for objects $X,Y,Z$
of $\cal{C}$ (in symbols $X,Y,Z\in|\cal{C}|$).

\begin{enumerate}
\item
  An \emph{algebra object} $(A,\mu,\eta)$ in $\cal{C}$ consists of an
  object $A$ and morphisms $\mu\colon A\otimes A\to A$ and
  $\eta\colon\1\to A$ of $\cal{C}$ such that:
\begin{equation}
\begin{aligned}
\xymatrix{
  (A\otimes A)\otimes A\ar[rr]^{\alpha_{A,A,A}}\ar[d]_{\mu\otimes\id_A}
    &&A\otimes(A\otimes A)\ar[d]^{\id_A\otimes\mu}\\
  A\otimes A\ar[dr]_{\mu}&&A\otimes A\ar[dl]^{\mu}\\
  & A
}
\end{aligned}
\end{equation}
and
\begin{equation}
\begin{aligned}
\xymatrix{
  \1\otimes A\ar[rr]^{\eta\otimes\id_A}\ar[ddrr]_{\lambda_A}
    &&A\otimes A\ar[dd]_{\mu}
    &&A\otimes\1\ar[ll]_{\id_A\otimes\eta}\ar[ddll]^{\rho_A}\\
  \\
  &&A
}
\end{aligned}
\end{equation}
  commute.
\item
  A \emph{coalgebra object} $(A,\Delta,\epsilon)$ in $\cal{C}$
  consists of an object $A$ and morphisms $\Delta\colon A\to A\otimes
  A$ and $\epsilon\colon A\to\1$ of $\cal{C}$ such that:
\begin{equation}
\begin{aligned}
\xymatrix{
  &A\ar[dl]_{\Delta}\ar[dr]^{\Delta}\\
  A\otimes A\ar[d]_{\Delta\otimes\id_A}
    &&A\otimes A\ar[d]^{\id_A\otimes\Delta}\\
  (A\otimes A)\otimes A\ar[rr]_{\alpha_{A,A,A}}&&A\otimes(A\otimes A)
}
\end{aligned}
\end{equation}
and
\begin{equation}
\begin{aligned}
\xymatrix{
  &&A\ar[ddll]_{\lambda_A^{-1}}\ar[ddrr]^{\rho_A^{-1}}\ar[dd]_{\Delta}\\
  \\
  \1\otimes A
    &&A\otimes A\ar[ll]^{\epsilon\otimes\id_A}\ar[rr]_{\id_A\otimes\epsilon}
    &&A\otimes 1
}
\end{aligned}
\end{equation}
  commute.
\item
  A \emph{homomorphism of algebras} $f\colon A\to A^\prime$ between
  two algebra objects $(A,\mu,\eta)$ and
  $(A^\prime,\mu^\prime,\eta^\prime)$ in $\cal{C}$ is a morphism $f$
  of $\cal{C}$ such that:
\begin{equation}
\begin{aligned}
\xymatrix{
  A\otimes A\ar[rr]^{\mu}\ar[dd]_{f\otimes f}
    &&A\ar[dd]^{f}\\
  \\
  A^\prime\otimes A^\prime\ar[rr]_{\mu^\prime}&&A^\prime
}
\end{aligned}
\qquad\mbox{and}\qquad
\begin{aligned}
\xymatrix{
  \1\ar[rr]^{\eta}\ar[rrdd]_{\eta^\prime}
    &&A\ar[dd]^{f}\\
  \\
  &&A^\prime
}
\end{aligned}
\end{equation}
  commute.
\item
  A \emph{homomorphism of coalgebras} $f\colon A\to A^\prime$ between
  two coalgebra objects $(A,\Delta,\epsilon)$ and
  $(A^\prime,\Delta^\prime,\epsilon^\prime)$ in $\cal{C}$ is a
  morphism $f$ of $\cal{C}$ such that:
\begin{equation}
\begin{aligned}
\xymatrix{
  A\ar[rr]^{\Delta}\ar[dd]_{f}
    &&A\otimes A\ar[dd]^{f\otimes f}\\
  \\
  A^\prime\ar[rr]_{\Delta^\prime}
  &&A^\prime\otimes A^\prime
}
\end{aligned}
\qquad\mbox{and}\qquad
\begin{aligned}
\xymatrix{
  A\ar[rrdd]^{\epsilon}\ar[dd]_{f}\\
  \\
  A^\prime\ar[rr]_{\epsilon^\prime}
    &&\1
}
\end{aligned}
\end{equation}
  commute.
\end{enumerate}
\end{defn}

\begin{defn}
Let $(\cal{C},\otimes,\1,\alpha,\lambda,\rho,\tau)$ be a symmetric
monoidal category.
\begin{enumerate}
\item
  A \emph{Frobenius algebra object} $(A,\mu,\eta,\Delta,\epsilon)$ in
  $\cal{C}$ consists of an object $A$ and of morphisms $\mu$, $\eta$,
  $\Delta$, $\epsilon$ of $\cal{C}$ such that:
\begin{enumerate}
\item
  $(A,\mu,\eta)$ is an algebra object in $\cal{C}$,
\item
  $(A,\Delta,\epsilon)$ is a coalgebra object in $\cal{C}$, and
\item
  the following compatibility condition, called the \emph{Frobenius
  relation}, holds,
\begin{equation}
\begin{aligned}
\xymatrix{
  &&A\otimes A\ar[dd]_{\mu}\ar[dll]_{\Delta\otimes\id_A}\ar[drr]^{\id_A\otimes\Delta}\\
  (A\otimes A)\otimes A\ar[dd]_{\alpha_{A,A,A}}
    &&&& A\otimes(A\otimes A)\ar[dd]^{\alpha^{-1}_{A,A,A}}\\
  &&A\ar[dd]_{\Delta}\\
  A\otimes (A\otimes A)\ar[drr]_{\id_A\otimes\mu}
    &&&& (A\otimes A)\otimes A\ar[dll]^{\mu\otimes\id_A}\\
  &&A\otimes A
}
\end{aligned}
\end{equation}
\end{enumerate}
\item
  A Frobenius algebra object $(A,\mu,\eta,\Delta,\epsilon)$ in
  $\cal{C}$ is called \emph{symmetric} if:
\begin{equation}
  \epsilon\circ\mu = \epsilon\circ\mu\circ\tau.
\end{equation}
  It is called \emph{commutative} if:
\begin{equation}
  \mu = \mu\circ\tau.
\end{equation}
\item
  Let $(A,\mu,\eta,\Delta,\epsilon)$ and
  $(A^\prime,\mu^\prime,\eta^\prime,\Delta^\prime,\epsilon^\prime)$ be
  Frobenius algebra objects in $\cal{C}$. A \emph{homomorphism of
  Frobenius algebras} $f\colon A\to A^\prime$ is a morphism $f$ of
  $\cal{C}$ which is both a homomorphism of algebra objects and a
  homomorphism of coalgebra objects.
\end{enumerate}
\end{defn}

Notice that for any Frobenius algebra object
$(A,\mu,\eta,\Delta,\epsilon)$ in $\cal{C}$, the object $A$ is
always a rigid object of $\cal{C}$. In $\cat{Vect}_k$ this means
that $A$ is finite-dimensional.

The unit object $\1\in|\cal{C}|$ forms an algebra object
$(\1,\lambda_\1,\id_\1)$ in $\cal{C}$ with multiplication
$\lambda_\1\colon\1\otimes\1\to\1$ and unit $\id_\1\colon\1\to\1$
as well as a coalgebra object $(\1,\lambda_\1^{-1},\id_\1)$
defining a commutative Frobenius algebra object in $\cal{C}$.
Given two algebra objects $(A,\mu_A,\eta_A)$ and
$(B,\mu_B,\eta_B)$ in $\cal{C}$, the tensor product $(A\otimes
B,\mu_{A\otimes B},\eta_{A\otimes B})$ forms an algebra object in
$\cal{C}$ with,
\begin{eqnarray}
  \mu_{A\otimes B} &=& (\mu_A\otimes\mu_B)
    \circ\alpha_{A,A,B\otimes B}^{-1}\circ (\id_A\otimes\alpha_{A,B,B})
    \circ(\id_A\otimes (\tau_{B,A}\otimes\id_B))\nn\\
   &&\quad\circ(\id_A\otimes\alpha_{B,A,B}^{-1})\circ\alpha_{A,B,A\otimes B},\\
  \eta_{A\otimes B} &=& (\eta_A\otimes\eta_B)\circ\lambda_\1^{-1}.
\end{eqnarray}
A similar result holds for coalgebra objects and for Frobenius algebra
objects in $\cal{C}$. Given two homomorphisms of algebra objects
$f\colon(A,\mu_A,\eta_A)\to(A^\prime,\mu_{A^\prime},\eta_{A^\prime})$
and
$g\colon(B,\mu_B,\eta_B)\to(B^\prime,\mu_{B^\prime},\eta_{B^\prime})$,
their tensor product $f\otimes g\colon(A\otimes B,\mu_{A\otimes
B},\eta_{A\otimes B})\to(A^\prime\otimes B^\prime,\mu_{A^\prime\otimes
B^\prime},\eta_{A^\prime\otimes B^\prime})$ forms a homomorphism of
algebra objects. A similar result holds for homomorphisms of coalgebra
and homomorphisms of Frobenius algebra objects.

The following definition plays a central role in the structure of
open-closed TQFTs.

\begin{defn}
Let $(\cal{C},\otimes,\1,\alpha,\lambda,\rho,\tau)$ be a symmetric
monoidal category. A \emph{knowledgeable Frobenius algebra}
$\A=(A,C,\imath,\imath^\ast)$ in $\cal{C}$ consists of,
\begin{itemize}
\item
  a symmetric Frobenius algebra $A=(A,\mu_A,\eta_A,\Delta_A,\epsilon_A)$,
\item
  a commutative Frobenius algebra $C=(C,\mu_C,\eta_C,\Delta_C,\epsilon_C)$,
\item
  morphisms $\imath\colon C\to A$ and $\imath^\ast\colon A\to C$ of $\cal{C}$,
\end{itemize}
such that $\imath\colon C\to A$ is a homomorphism of algebra objects
in $\cal{C}$ and,
\begin{alignat}{2}
\label{eq_kfrob1}
  \mu_A\circ(\imath\otimes\id_A)
    &= \mu_A\circ\tau_{A,A}\circ(\imath\otimes\id_A)
    &\qquad& \mbox{(knowledge),}\\
\label{eq_kfrob2}
  \epsilon_C\circ\mu_C\circ(\id_C\otimes\imath^\ast)
    &= \epsilon_A\circ\mu_A\circ(\imath\otimes\id_A)
    &\qquad& \mbox{(duality),}\\
\label{eq_kfrob3}
  \mu_A\circ\tau_{A,A}\circ\Delta_A
    &= \imath\circ\imath^{\ast}
    &\qquad& \mbox{(Cardy condition).}
\end{alignat}
\end{defn}
Condition~\eqref{eq_kfrob2} says that $\imath^\ast$ is the morphism
dual to $\imath$. Together with the fact that $\imath$ is an algebra
homomorphism, this implies that $\imath^\ast\colon A\to C$ is a
homomorphism of coalgebras in $\cal{C}$.

If $\cal{C}=\cat{Vect}_k$, the condition~\eqref{eq_kfrob1} states
that the image of $C$ under $\imath$ is contained in the centre of
$A$, $\imath(C)\subseteq Z(A)$. The name \emph{knowledgeable
Frobenius algebra} is meant to indicate that the symmetric
Frobenius algebra $A$ knows something about its centre. This is
specified precisely by $C$, $\imath$ and $\imath^\ast$. Notice
that the centre $Z(A)$ itself cannot be characterized\footnote{We
thank James Dolan and John Baez for pointing this out.} by
requiring the commutativity of diagrams labeled by objects and
morphisms of $\cal{C}$.

It is not difficult to see that every strongly separable algebra $A$
can be equipped with the structure of a Frobenius algebra such that
$(A,Z(A),\imath,\imath^\ast)$ forms a knowledgeable Frobenius algebra
with the inclusion $\imath\colon Z(A)\to A$ and an appropriately
chosen Frobenius algebra structure on $Z(A)$. There are also examples
$(A,C,\imath,\imath^\ast)$ of knowledgeable Frobenius algebras in
$\cat{Vect}_k$ in which $C$ is not the centre of $A$. For more
details, we refer to~\cite{LP}.

\begin{defn}
A \emph{homomorphism of knowledgeable Frobenius algebras}
\begin{equation}
  \phi\maps (A,C,\imath,\imath^\ast) \to (A',C',\imath',\imath'^{\ast})
\end{equation}
in the symmetric monoidal category $\cal{C}$ is a pair
$\phi=(\phi_1,\phi_2)$ of Frobenius algebra homomorphisms $\phi_1
\maps A \to A'$ and $\phi_2\maps C
\to C'$ such that
\begin{equation}
\begin{aligned}
  \xymatrix{
    C \ar[r]^{\phi_2} \ar[d]_{\imath} & C' \ar[d]^{\imath'} \\
    A \ar[r]_{\phi_1} & A'
  }
\end{aligned}
\qquad\mbox{and}\qquad
\begin{aligned}
  \xymatrix{
    A \ar[r]^{\phi_1} \ar[d]_{\imath^{\ast}} & A'\ar[d]^{\imath'^{\ast}} \\
    C \ar[r]_{\phi_2} & C'
  }
\end{aligned}
\end{equation}
commute.
\end{defn}

\begin{defn}
Let $\cal{C}$ be a symmetric monoidal category. By
\cat{K-Frob}$(\cal{C})$ we denote the category of knowledgeable
Frobenius algebras in $\cal{C}$ and their homomorphisms.
\end{defn}

\begin{prop}
Let $\cal{C}$ be a symmetric monoidal category. The category
$\cat{K-Frob}(\cal{C})$ forms a symmetric monoidal category as
follows. The tensor product of two knowledgeable Frobenius algebra
objects $\A=(A,C,\imath,\imath^{\ast})$ and
$\A'=(A',C',\imath',\imath'^{\ast})$ is defined as
$\A\otimes\A^\prime:=(A\otimes A^\prime,C\otimes
C^\prime,\imath\otimes\imath^\prime,\imath^\ast\otimes{\imath^\prime}^\ast)$. The
unit object is given by $\underline{\1}:=(\1,\1,\id_\1,\id_\1)$, and
the associativity and unit constraints and the symmetric braiding are
induced by those of $\cal{C}$. Given two homomorphisms
$\phi=(\phi_1,\phi_2)$ and $\psi=(\psi_1,\psi_2)$ of knowledgeable
Frobenius algebras, their tensor product is defined as
$\phi\otimes\psi:=(\phi_1\otimes\psi_1,\phi_2\otimes\psi_2)$.
\end{prop}

\subsection{The category $\Thfrob$}

In this section, we define the category $\Thfrob$, called the
\emph{theory of knowledgeable Frobenius algebras}. The description
that follows is designed to make $\Thfrob$ the symmetric monoidal
category freely generated by a knowledgeable Frobenius algebra, and
the terminology `theory of~$\ldots$' indicates that knowledgeable
Frobenius algebras in any symmetric monoidal category $\cal{C}$
arise precisely as the symmetric monoidal functors
$\Thfrob\to\cal{C}$. This is in analogy to the theory of algebraic
theories in which one uses `with finite products' rather than
'symmetric monoidal'. Readers who are interested in the topology of
open-closed cobordisms rather than in the abstract description of
knowledgeable Frobenius algebras may wish to look briefly at
Proposition~\ref{PROPthfrob} and then directly proceed to
Section~\ref{sect_cob}.

The subsequent definition follows the construction of the `free
category with group structure' given by Laplaza~\cite{Laplaza}.  It
forms an example of a symmetric monoidal sketch, a structure slightly
more general than an operad or a PROP, see for example~\cite{Hyland}
for the definition of symmetric monoidal sketches and a discussion of
their freeness properties.

\begin{defn}
The category $\Thfrob$ is defined as follows. Its objects are the
elements of the free $\{\1,\otimes\}$-algebra over the two element set
$\{A,C \}$. These are words of a formal language that are defined by
the following requirements,
\begin{itemize}
\item
  The symbols $\1$, $A$ and $C$ are objects of $\Thfrob$.
\item
  If $X$ and $Y$ are objects of $\Thfrob$, then $(X\otimes Y)$ is an
  object of $\Thfrob$.
\end{itemize}
We now describe the edges of a graph $\mathcal{G}$ whose vertices are the objects of
$\Thfrob$. There are edges,
\begin{gather}
 \mu_A \maps A \ten A \to A, \quad \eta_A \maps \1 \to A, \quad
 \Delta_A \maps A \to A \ten A, \quad \epsilon_A \maps \1 \to A,\nn\\
 \mu_C \maps C \ten C \to C, \quad \eta_C \maps \1 \to C, \quad
 \Delta_C \maps C \to C \ten C, \quad \epsilon_C \maps \1 \to C,\\
 \imath \maps C \to A, \qquad \imath^{\ast} \maps A\to C,\nn
\end{gather}
and for all objects $X$,$Y$,$Z$ there are to be edges
\begin{gather}
 \alpha_{X,Y,Z} \maps (X \ten Y) \ten Z \to X \ten (Y \ten Z),
 \quad
 \tau_{X,Y} \maps X \ten Y \to Y \ten X,\nn\\
 \lambda_X \maps \1 \ten X \to X ,
 \quad
 \rho_X \maps X \ten \1 \to X,\\
 \bar{\alpha}_{X,Y,Z} \maps X \ten (Y \ten Z) \to (X \ten Y) \ten Z,
 \quad
 \bar{\tau}_{X,Y} \maps Y \ten X \to X \ten Y,\nn\\
 \bar{\lambda}_X \maps  X  \to\1 \ten X,
 \quad
 \bar{\rho}_X \maps  X \to X \ten \1.
\end{gather}
For every edge $f\maps X\to Y$ and for every object $Z$, there are to
be edges $Z \ten f \maps Z \ten X \to Z \ten Y$, $f \ten Z \maps X
\ten Z \to Y \ten Z$. These edges are to be interpreted as words in a
formal language and are considered distinct if they have distinct
names.

Let $\mathcal{H}$ be the category freely generated by the graph
$\mathcal{G}$. We now describe a congruence on the category
$\mathcal{H}$. We define a relation $\sim$ as follows. We require the
relations making $(A,\mu_A,\eta_A,\Delta_A,\epsilon_A)$ a symmetric
Frobenius algebra object, those making
$(C,\mu_C,\eta_C,\Delta_C,\epsilon_C)$ a commutative Frobenius algebra
object, those making $\imath\maps C \to A$ an algebra homomorphism as
well as~\eqref{eq_kfrob1}, \eqref{eq_kfrob2},
and~\eqref{eq_kfrob3}. The relations making $\alpha_{X,Y,Z}$,
$\lambda_{X}$, and $\rho_{X}$ satisfy the pentagon and triangle axioms
of a monoidal category as well as those making $\tau_{X,Y}$ a
symmetric braiding, are required for all objects $X$,$Y$,$Z$. We also
require the following relations for all objects $X,Y$ and morphisms
$p,q,t,s$ of $\mathcal{H}$,
\begin{equation}
\begin{aligned}
  (X \ten p)(X \ten q) \sim X \ten (pq),\qquad
  (p \ten X)(q\ten X) \sim (pq) \ten X,\\
  (t \ten Y)(X \ten s) \sim (X \ten s)(t \ten Y),\quad
  \id_{X \ten Y}\sim X \ten \id_Y \sim \id_X \ten Y,
\end{aligned}
\end{equation}
that make $\ten$ a functor. Then we require the relations that
assert the naturality of
$\alpha,\lambda,\rho,\tau,\bar{\alpha},\bar{\lambda},\bar{\rho},\bar{\tau}$
and that each pair $e$ and $\bar{e}$ of edges of the graph form
the inverses of each other.  Finally, we have all expansions by
$\ten$, \ie\ for each relation $a \sim b$, we include the
relations $a \ten X \sim b \ten X$ and $X \ten a \sim X \ten b$
for all objects $X$, and all those relations obtained from these
by a finite number of applications of this process. The category
$\Thfrob$ is the category $\mathcal{H}$ modulo the category
congruence generated by $\sim$.
\end{defn}

It is clear from the description above that $\Thfrob$ contains a
knowledgeable Frobenius algebra object $(A,C,\imath,\imath^{\ast})$
which we call the knowledgeable Frobenius algebra object
\emph{generating} $\Thfrob$. Indeed, $\Thfrob$ is the symmetric
monoidal category freely generated by a knowledgeable Frobenius
algebra. Its basic property is that for any knowledgeable Frobenius
algebra $\A'=(A',C',\imath',\imath'^{\ast})$ in $\cal{C}$, there is
exactly one \emph{strict} symmetric monoidal functor $F_{\A'} \maps
\Thfrob \to \cal{C}$ which maps $(A,C,\imath,\imath^{\ast})$ to
$(A',C',\imath',\imath'^{\ast})$ and $\1 \in \Thfrob$ to
$\1\in\cal{C}$.

An interesting question to ask is whether or not homomorphisms of
knowledgeable Frobenius algebras are induced in some way by
$\Thfrob$. This question is answered by the following
proposition.

\begin{prop}
\label{PROPthfrob}
Let $\cal{C}$ be a symmetric monoidal category. The category
\begin{equation}
  \cat{Symm-Mon}(\Thfrob,\cal{C})
\end{equation}
of symmetric monoidal functors $\Thfrob\to\cal{C}$ and their monoidal
natural transformations is equivalent as a symmetric monoidal category
to the category $\cat{K-Frob($\cal{C}$)}$.
\end{prop}

This proposition is the reason for calling $\Thfrob$ the
\emph{theory of knowledgeable Frobenius algebras}. For easier reference,
we collect the definitions of symmetric monoidal functors, monoidal
natural transformations and of $\cat{Symm-Mon}$ in
Appendix~\ref{app_moncat}.

\begin{proof}
Let $(A,C,\imath,\imath^{\ast})$ be the knowledgeable Frobenius
algebra generating $\Thfrob$, and let $\psi\maps\Thfrob\to\cal{C}$ be
a symmetric monoidal functor. It is clear that the image of
$(A,C,\imath,\imath^{\ast})$ under $\psi$, together with the coherence
isomorphisms $\psi_0$ and $\psi_2$ of $\psi=(\psi,\psi_2,\psi_0)$,
defines a knowledgeable Frobenius algebra
$(\psi(A),\psi(C),\psi(\imath),\psi(\imath^{\ast}))$ in $\cal{C}$. The
symmetric Frobenius algebra structure on $\psi(A)$ is given by
\begin{equation}
  \psi(A)=\big(\psi(A),\psi(\mu_A)\circ\psi_2,\psi(\eta_A)\circ\psi_0,
    \psi_2^{-1}\circ\psi(\Delta_A),\psi_0^{-1}\circ\psi(\epsilon_A)\big).
\end{equation}
The commutative Frobenius algebra structure on $\psi(C)$ is given by
\begin{equation}
  \psi(C)=\big(\psi(C),\psi(\mu_C)\circ\psi_2,\psi(\eta_C)\circ\psi_0,
    \psi_2^{-1}\circ\psi(\Delta_C),\psi_0^{-1}\circ\psi(\epsilon_C)\big).
\end{equation}
This defines a mapping on objects
\begin{eqnarray}
  \Gamma\maps\cat{Symm-Mon}(\Thfrob,\cal{C}) &\to& \cat{K-Frob($\cal{C}$)}\\
  \psi &\mapsto& (\psi(A),\psi(C),\psi(\imath),\psi(\imath^{\ast})).\nn
\end{eqnarray}
We now extend $\Gamma$ to a functor by defining it on morphisms.

If $\phi\maps\psi\Rightarrow\psi'$ is a monoidal natural
transformation, then $\phi$ assigns to each object $X$ in $\Thfrob$ a
map $\phi_{X}\maps\psi(X)\to\psi'(X)$ in $\cal{C}$. However, since
every object in $\Thfrob$ is the tensor product of $A$'s and $C$'s and
$\1$'s, the fact that $\phi$ is a \emph{monoidal} natural
transformation means that the $\phi_{X}$ are completely determined by
two maps $\phi_1\maps\psi(A)\to\psi'(A)$ and
$\phi_2\maps\psi(C)\to\psi'(C)$. The \emph{naturality} of $\phi$ means
that the $\phi_i$ are compatible with the images of all the morphisms
in $\Thfrob$. Since all of the morphisms in $\Thfrob$ are built up
from the generators:
\begin{gather}
 \mu_A \maps A \ten A \to A, \quad \eta_A \maps \1 \to A, \quad
 \Delta_A \maps A \to A \ten A, \quad \epsilon_A \maps \1 \to A,\nn\\
 \mu_C \maps C \ten C \to C, \quad \eta_C \maps \1 \to C, \quad
 \Delta_C \maps C \to C \ten C, \quad \epsilon_C \maps \1 \to C,\\
 \imath \maps C \to A, \qquad \imath^{\ast} \maps A\to C,\nn
\end{gather}
(and the structure maps $\alpha$, $\rho$, $\lambda$, $\tau$),
naturality can be expressed by the commutativity of 10 diagrams
involving the 10 generating morphisms of $\Thfrob$. For example,
corresponding to $\mu_A\maps A\ten A\to A$ and $\eta_A\maps\1\to A$,
we have the two diagrams:
\begin{gather}
\begin{aligned}
 \xy
    (-45,8)*+{\psi(A \ten A)}="tl";
    (-18,8)*+{\psi(A) \ten \psi(A)}="tm";
    (18,8)*+{\psi'(A) \ten \psi'(A) }="tm'";
    (48,8)*+{\psi'(A \ten A)}="tr";
    (-45,-8)*+{\psi(A)}="bl";
    (48,-8)*+{\psi'(A)}="br";
        {\ar^-{\psi_2^{-1}} "tl";"tm"};
        {\ar^-{\phi_1 \ten \phi_1} "tm";"tm'"};
         {\ar^-{\psi_{2}'^{-1}} "tm'";"tr"};
        {\ar_-{\psi(\mu_A)} "tl";"bl"};
        {\ar_-{\phi_1} "bl";"br"};
        {\ar^-{\psi^\prime(\mu_A)} "tr";"br"};
 \endxy
\end{aligned}\\
\begin{aligned}
 \xy
    (-20,8)*+{\psi(\1)}="tl";
    (0,8)*+{ \1 }="tm";
    (20,8)*+{\psi'(\1)}="tr";
    (-20,-8)*+{\psi(A)}="bl";
    (20,-8)*+{\psi'(A)}="br";
        {\ar_-{\psi(\eta_A)} "tl";"bl"};
        {\ar^-{\psi_0^{-1}} "tl";"tm"};
        {\ar^-{\psi_0'} "tm";"tr"};
        {\ar_-{\phi_1} "bl";"br"};
        {\ar^-{\psi'(\eta_A)} "tr";"br"};
 \endxy
\end{aligned}
\end{gather}
which amount to saying that $\phi_1$ is an algebra homomorphism
$\psi(A)\to\psi'(A)$. Together with the conditions for the
generators $\Delta_A\maps A\to A\ten A$ and $\epsilon_A\maps
A\to\1$, we have that $\phi_1$ is a Frobenius algebra homomorphism
from $\psi(A)$ to $\psi'(A)$. Similarly, the diagrams corresponding
to the generators with a $C$ subscript imply that $\phi_2$ is a
Frobenius algebra homomorphism $\psi(C)\to\psi'(C)$. The conditions
on the images of the generators $\imath\maps C\to A$ and
$\imath^{\ast}\maps A\to C$ produce the requirement that the two
diagrams:
\begin{equation}
\begin{aligned}
 \xymatrix{\psi(C) \ar[r]^{\phi_2} \ar[d]_{\psi(\imath)} & \psi'(C)
    \ar[d]^{\psi'(\imath)} \\ \psi(A)  \ar[r]_{\phi_1} & \psi'(A)}
\end{aligned}
\qquad\mbox{and}\qquad
\begin{aligned}
 \xymatrix{\psi(A) \ar[r]^{\phi_1} \ar[d]_{\psi(\imath^{\ast})} & \psi'(A)
    \ar[d]^{\psi'(\imath^{\ast})} \\ \psi(C)  \ar[r]_{\phi_2} &
    \psi'(C) }
\end{aligned}
\end{equation}
commute. Hence, the monoidal natural transformation $\phi$ defines a
morphism of knowledgeable Frobenius algebras in $\cal{C}$. This
assignment clearly preserves the monoidal structure and symmetry up to
isomorphism. Thus, it is clear that one can define a symmetric
monoidal functor $\Gamma=(\Gamma,\Gamma_2,\Gamma_0)\colon
\cat{Symm-Mon}(\Thfrob)\to\cat{K-Frob}(\cal{C})$.

Conversely, given any knowledgeable Frobenius algebra
$\A'=(A',C',\imath',\imath'^{\ast})$ in $\cal{C}$, then by the
remarks preceding this proposition, there is an assignment
\begin{eqnarray}
  \overline{\Gamma}\maps\cat{K-Frob($\cal{C}$)}
    &\to& \cat{Symm-Mon}(\Thfrob,\cal{C}) \\
  (A',C',\imath',\imath'^{\ast}) &\mapsto& F_{\A'},
\end{eqnarray}
where $F_{\A'}$ is the strict symmetric monoidal functor sending the
knowledgeable Frobenius algebra $(A,C,\imath,\imath^{\ast})$
generating $\Thfrob$ to the knowledgeable Frobenius algebra
$(A',C',\imath',\imath'^{\ast})$ in $\cal{C}$. Furthermore, it is
clear from the discussion above that a homomorphism of knowledgeable
Frobenius algebras $\phi\maps\A_1\to\A_2$ defines a monoidal natural
transformation $\phi\maps F_{\A_1}\to F_{\A_2}$. Thus, it is clear
that $\overline{\Gamma}$ extends to a symmetric monoidal functor
$\overline{\Gamma}=(\overline{\Gamma},\overline{\Gamma}_2,\overline{\Gamma}_0)\colon
\cat{K-Frob}(\cal{C})\to\cat{Symm-Mon}(\Thfrob,\cal{C})$.

To see that $\Gamma$ and $\overline{\Gamma}$ define an equivalence of
categories, let $\A'=(A',C',\imath',\imath'^{\ast})$ be a
knowledgeable Frobenius algebra in $\cal{C}$. The composite
$\Gamma\overline{\Gamma}(\A')=\Gamma(F_{\A'})=\A'$ since $F_{\A'}$ is
a strict symmetric monoidal functor.  Hence,
$\Gamma\overline{\Gamma}=\id_{\cat{K-Frob}(\cal{C})}$. Now let
$\psi\maps\Thfrob\to\cal{C}$ be a symmetric monoidal functor and
consider the composite $\overline{\Gamma}\Gamma(\psi)$. Let
$\tilde{\A}=\big(\psi(A),\psi(\mu_A)\circ\psi_2,\psi(\eta_A)\circ\psi_0,
\psi_2^{-1}\circ\psi(\Delta_A),\psi_0^{-1}\circ\psi(\epsilon_A)\big)$
so that $\overline{\Gamma}\Gamma(\psi)=F_{\tilde{\A}}$. We define a
monoidal natural isomorphism $\theta\maps\psi\Rightarrow
F_{\tilde{\A}}$ on the generators as follows:
\begin{eqnarray}
  \theta_{\1}\maps \psi(\1)\to F_{\tilde{\A}}(\1)=\1 &:=& \psi_0^{-1},\nn\\
  \theta_{A}\maps \psi(A) \to F_{\tilde{\A}}(A)=\psi(A) &:=& 1_A, \\
  \theta_{C}\maps \psi(C) \to F_{\tilde{\A}}(C)=\psi(C) &:=& 1_C.\nn
\end{eqnarray}
The condition that $\theta$ be monoidal implies that $\theta_{A\ten
A}=(\psi_2^{-1}) _{A\ten A}$, $\theta_{A\ten C}=(\psi_2^{-1})
_{A\ten C}$, $\theta_{C\ten A}=(\psi_2^{-1}) _{C\ten A}$, and
$\theta_{C\ten C}=(\psi_2^{-1}) _{C\ten C}$. Since $\Thfrob$ is
generated by $\1$ ,$A$, and $C$, this assignment uniquely defines a
monoidal natural isomorphism. Hence,
$\overline{\Gamma}\Gamma(\psi)\cong\psi$ so that $\overline{\Gamma}$
and $\Gamma$ define a monoidal equivalence of categories.
\end{proof}

%
\section{The category of open-closed cobordisms}
%
\label{sect_cob}

In this section, we define and study the category $\twocob$ of
open-closed cobordisms. Open-closed cobordisms form a special sort
of compact smooth $2$-manifolds with corners that have a particular
global structure. If one decomposes their boundary minus the corners
into connected components, these components are either \emph{black}
or \emph{coloured} with elements of some given set $S$. Every corner
is required to separate a black boundary component from a coloured
one\footnote{In this terminology, black is not considered a
colour.}.

These $2$-manifolds with corners are viewed as cobordisms between
their black boundaries, and they can be composed by gluing them along
their black boundaries subject to a matching condition for the colours
of the other boundary components. In the conformal field theory
literature, the coloured boundary components are referred to as
\emph{free boundaries} and the colours as \emph{boundary conditions}.

$2$-manifolds with corners with this sort of global structure form a
special case of $\left<2\right>$-manifolds according to
J\"anich~\cite{Jan}. For an overview and a very convenient notation,
we refer to the introduction of the article~\cite{Lr} by Laures.

In the following two subsections, we present all definitions for a
generic set of colours $S$. Starting in
Subsection~\ref{sect_generators}, the generators and relations
description of $\twocob$ is developed only for the case of a single
colour, $S=\{\ast\}$. We finally return to the case of a generic set
of colours $S$ in Section~\ref{sect_labels}.

\subsection{$\left<2\right>$-manifolds}
\label{sect_2mfd}

\subsubsection{Manifolds with corners}

A $k$-dimensional \emph{manifold with corners} $M$ is a
topological manifold with boundary that is equipped with a smooth
structure with corners. A smooth structure with corners is defined
as follows. A \emph{smooth atlas with corners} is a family
${\{(U_\alpha,\phi_\alpha)\}}_{\alpha\in I}$ of coordinate systems
such that the $U_\alpha\subseteq M$ are open subsets which cover
$M$, and the
\begin{equation}
  \phi_\alpha\colon U_\alpha\to\phi_\alpha(U_\alpha)\subseteq \R^k_+
\end{equation}
are homeomorphisms onto open subsets of $\R_+^k:={[0,\infty)}^k$. The
transition functions
\begin{equation}
  \phi_\beta\circ\phi_\alpha^{-1}\colon\phi_\alpha(U_{\alpha\beta})\to\phi_\beta(U_{\alpha\beta})
\end{equation}
for $U_{\alpha\beta}:=U_\alpha\cap U_\beta\neq\emptyset$ are required
to be the restrictions to $\R_+^k$ of diffeomorphisms between open
subsets of $\R^k$. Two such atlases are considered equivalent if their
union is a smooth atlas with corners, and a \emph{smooth structure
with corners} is an equivalence class of such atlases.

A \emph{smooth map} $f\colon M\to N$ between manifolds with corners
$M$ and $N$ is a continuous map for which the following condition
holds. Let ${\{(U_\alpha,\phi_\alpha)\}}_{\alpha\in I}$ and
${\{(V_\beta,\psi_\beta)\}}_{\beta\in J}$ be atlases that represent
the smooth structures with corners of $M$ and $N$, respectively. For
every $p\in M$ and for every $\alpha\in I$, $\beta\in J$ with $p\in
U_\alpha$ and $f(p)\in V_\beta$, we require that the map
\begin{equation}
  \psi_\beta\circ f\circ\phi_\alpha^{-1}|_{\phi_\alpha(U_\alpha\cap f^{-1}(V_\beta))}\colon
  \phi_\alpha(U_\alpha\cap f^{-1}(V_\beta))\to\psi_\beta(f(U_\alpha)\cap V_\beta)
\end{equation}
is the restriction to $\R_+^m$ of a smooth map between open subsets of
$\R^m$ and $\R^p$ for $m$ and $p$ the dimensions of $M$ and $N$,
respectively.

\subsubsection{Manifolds with faces}

For each $p\in M$, we define $c(p)\in\N_0$ to be the number of zero
coefficients of $\phi_\alpha(p)\in\R^k$ for some $\alpha\in I$ for
which $p\in U_\alpha$. A \emph{connected face} of $M$ is the closure
of a component of $\{\,p\in M\colon\,c(p)=1\,\}$. A \emph{face} is a
free union of pairwise disjoint connected faces. This includes the
possibility that a face can be empty.

A $k$-dimensional \emph{manifold with faces} $M$ is a $k$-dimensional
manifold with corners such that each $p\in M$ is contained in $c(p)$
different connected faces. Notice that every face of $M$ is itself a
manifold with faces.

\subsubsection{$\left<n\right>$-manifolds}

A $k$-dimensional $\left<n\right>$-manifold $M$ is a $k$-dimensional
manifold with faces with a specified tuple $(\del_0
M,\ldots,\del_{n-1} M)$ of faces of $M$ such that the following two
conditions hold.
\begin{enumerate}
\item
  $\del_0 M\cup\cdots\cup\del_{n-1} M=\del M$. Here $\del M$ denotes
  the boundary of $M$ as a topological manifold.
\item
  $\del_j M\cap\del_\ell M$ is a face of both $\del_j M$ and
  $\del_\ell M$ for all $j\neq\ell$.
\end{enumerate}
Notice that a $\left<0\right>$-manifold is just a manifold without
boundary while a $\left<1\right>$-manifold is a manifold with
boundary. A diffeomorphism $f\colon M\to N$ between two
$\left<n\right>$-manifolds is a diffeomorphism of the underlying
manifolds with corners such that $f(\del_j M)=\del_j N$ for all $j$.

The following notation is taken from Laures~\cite{Lr}. Let
$\underline{2}$ denote the category associated with the partially
ordered set $\{0,1\}$, $0\leq 1$, \ie\ the category freely generated
by the graph $0\stackrel{\ast}{\longrightarrow} 1$. Denote by
$\underline{2}^n$ the $n$-fold Cartesian product of $\underline{2}$
and equip its set of objects ${\{0,1\}}^n$ with the corresponding
partial order. An $\left<n\right>$-\emph{diagram} is a functor
$\underline{2}^n\to\cat{Top}$. We use the term
$\left<n\right>$-diagram \emph{of inclusions} for an
$\left<n\right>$-diagram which sends each morphism of
$\underline{2}^n$ to an inclusion, and so on.

Every $\left<n\right>$-manifold $M$ gives rise to an
$\left<n\right>$-diagram $M\colon\underline{2}^n\to\cat{Top}$ of
inclusions as follows. For the objects
$a=(a_0,\ldots,a_{n-1})\in|\underline{2}^n|$, write
$a^\prime:=(1-a_0,\ldots,1-a_{n-1})$, and denote the standard basis of
$\R^n$ by $(e_0,\ldots,e_{n-1})$. The functor
$M\colon\underline{2}^n\to\cat{Top}$ is defined on the objects by,
\begin{equation}
  M(a):=\bigcap_{i\in\{\,i\colon\, a\leq e_i^\prime\,\}}\del_i M,
\end{equation}
if $a\neq(1,\ldots,1)$, and by $M((1,\ldots,1)):=M$. The functor sends
the morphisms of $\underline{2}^n$ to the obvious inclusions.

For all $a\in|\underline{2}^n|$, the face $M(a)$ of $M$ forms a
$\left<\ell\right>$-manifold itself for which
$\ell=\sum_{i=1}^{n-1}a_i$. An orientation of $M$ induces orientations
on the $M(a)$ as usual. The product of an $\left<n\right>$-manifold
$M$ with a $\left<p\right>$-manifold $N$ forms an
$\left<n+p\right>$-manifold, denoted by $M\times N$. The structure of
its faces can be read off from the functor
\begin{equation}
\label{eq_prodmfd}
  M\times N\colon\underline{2}^{n+p}\simeq\underline{2}^n\times\underline{2}^p
  \stackrel{M\times N}{\longrightarrow}\cat{Top}\times\cat{Top}
  \stackrel{\times}{\longrightarrow}\cat{Top}.
\end{equation}
The half-line $\R_+:=[0,\infty)$ is a $1$-dimensional manifold
with boundary, \ie\ a $1$-dimensional $\left<1\right>$-manifold.
The product of~\eqref{eq_prodmfd} then equips $\R_+^n$ with the
structure of an $n$-dimensional $\left<n\right>$-manifold.

The case that is relevant in the following is that of
$2$-dimensional $\left<2\right>$-manifolds. These are
$2$-dimensional manifolds with faces with a pair of specified
faces $(\del_0 M,\del_1 M)$ such that $\del_0 M\cup\del_1 M=\del
M$ and $\del_0 M\cap\del_1 M$ is a face of both $\del_0 M$ and
$\del_1 M$.  The following diagram shows the faces of one of the
typical $2$-dimensional $\left<2\right>$-manifolds $M$ that are
used below.
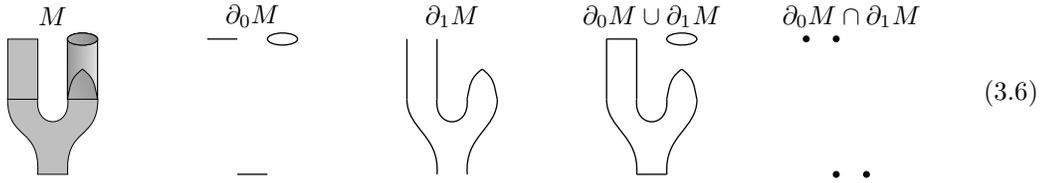
\begin{equation}
\label{eq_faces}
\begin{aligned}
\psset{xunit=.4cm,yunit=.4cm}
\begin{pspicture}[0.5](4,5.4)
  \rput(2, 0){\multl}
  \rput(1, 2.5){\medidentl}
  \rput(2.92, 2.5){\ctl}
  \rput(2, 5.25){$M$}
\end{pspicture}
\qquad\quad
\begin{pspicture}[0.5](4,5.4)
  \rput(2, 5.255){$\partial_0 M$}
  \psline[linewidth=0.5pt](1.5,4.5)(.5,4.5)
  \psline[linewidth=0.5pt](2.5,0)(1.5,0)
  \psellipse[linewidth=0.5pt](3,4.5)(.5,0.2)
\end{pspicture}
\qquad\quad
\begin{pspicture}[0.5](4,5.4)
  \rput(2, 5.25){$\partial_1 M$}
  \rput(2,0){
    \psline[linewidth=0.5pt](-1.5,4.5)(-1.5,2.5)
    \psbezier[linewidth=0.5pt](-1.5,2.5)(-1.5,1.5)(-.4,1.3)(-.5,0)
    \psline[linewidth=0.5pt](-.5,4.5)(-.5,2.5)
  }
  \rput(3,0){
    \psbezier[linewidth=0.5pt](-.5,0)(-.6,1.3)(.5,1.5)(.5,2.5)
  }
  \rput(3.08,2.5){
    \psbezier[linewidth=0.5pt](-.58,0)(-.48,.5)(-.48,.7)(-.08,1)
    \psbezier[linewidth=0.5pt](-.08,1)(.32,.7)(.32,.5)(.42,0)
  }
  \rput(2,.5){
    \psbezier[linewidth=0.5pt](-.5,2)(-.6,1)(0.6,1)(.5,2)
  }
\end{pspicture}
\qquad\quad
\begin{pspicture}[0.5](4,5.4)
  \rput(2, 5.25){$\partial_0 M \cup\partial_1 M$}
  \rput(2,0){
     \psline[linewidth=0.5pt](-1.5,4.5)(-1.5,2.5)
     \psbezier[linewidth=0.5pt](-1.5,2.5)(-1.5,1.5)(-.4,1.3)(-.5,0)
     \psline[linewidth=0.5pt](-.5,4.5)(-.5,2.5)
  }
  \rput(3,0){
    \psbezier[linewidth=0.5pt](-.5,0)(-.6,1.3)(.5,1.5)(.5,2.5)
  }
  \rput(3.08,2.5){
    \psbezier[linewidth=0.5pt](-.58,0)(-.48,.5)(-.48,.7)(-.08,1)
    \psbezier[linewidth=0.5pt](-.08,1)(.32,.7)(.32,.5)(.42,0)
  }
  \rput(2,.5){
    \psbezier[linewidth=0.5pt](-.5,2)(-.6,1)(0.6,1)(.5,2)
  }
  \psline[linewidth=0.5pt](1.5,4.5)(.5,4.5)
  \psline[linewidth=0.5pt](2.5,0)(1.5,0)
  \psellipse[linewidth=0.5pt](3,4.5)(.5,0.2)
\end{pspicture}
\qquad\quad
\begin{pspicture}[0.5](4,5.4)
  \rput(2, 5.25){$\partial_0 M \cap\partial_1 M$}
  \psdots(.5,4.5)(1.5,4.5)(1.5,0)(2.5,0)
\end{pspicture}
\end{aligned}
\end{equation}
The $\left<2\right>$-diagram $M\colon\underline{2}^2\to\cat{Top}$ of
inclusions is the following commutative square:
\begin{equation}
\begin{aligned}
\xymatrix{
  \del_0 M\cap\del_1 M\ar[rr]^{M(\id_0\times\ast)}\ar[dd]_{M(\ast\times\id_0)}&&
    \del_0 M\ar[dd]^{M(\ast\times\id_1)}\\
  \\
  \del_1 M\ar[rr]_{M(\id_1\times\ast)}&& M
}
\end{aligned}
\end{equation}
Another example of a manifold with corners $M$ which is embedded
in $\R^3$ is depicted in~\eqref{eq_figure}. It has the structure
of a $2$-dimensional $\left<2\right>$-manifold when one chooses
$\del_0 M$ to be the union of the top and bottom boundaries of the
picture, similarly to~\eqref{eq_faces}.

\subsubsection{Collars}

In order to glue $\left<2\right>$-manifolds along specified faces, we
need the following technical results.

\begin{lem}[Lemma~2.1.6 of~\cite{Lr}]
Each $\left<n\right>$-manifold $M$ admits an $\left<n\right>$-diagram $C$ of embeddings
\begin{equation}
  C(a\rightarrow b)\colon\R_+^n(a^\prime)\times M(a)\hookrightarrow\R_+^n(b^\prime)\times M(b)
\end{equation}
such that the restriction $C(a\rightarrow b)|_{\R_+^n(b^\prime)\times
M(a)}=\id_{\R_+^n(b^\prime)}\times M(a\rightarrow b)$ is the inclusion
map.
\end{lem}

In particular, for every $\left<2\right>$-manifold $M$, there is a
commutative square $C\colon\underline{2}^2\to\cat{Top}$ of
embeddings,
\begin{equation}
\begin{aligned}
\xymatrix{
  \R_+^2\times(\del_0 M\cap\del_1 M)\ar[rr]^{C(\id_0\times\ast)}\ar[dd]_{C(\ast\times\id_0)}&&
    \del_0\R_+^2\times\del_1 M\ar[dd]^{C(\ast\times\id_1)}\\
  \\
  \del_1\R_+^2\times\del_0 M\ar[rr]_{C(\id_1\times\ast)}&&\{(0,0)\}\times M
}
\end{aligned}
\end{equation}
such that the following restrictions are inclusions,
\begin{eqnarray}
  C(\id_0\times\ast)|_{\del_0\R_+^2\times(\del_0 M \cap\del_1 M)}
    &=& \id_{\del_0\R_+^2}\times M(\id_0\times\ast),\\
  C(\ast\times\id_0)|_{\del_1\R_+^2\times(\del_0 M \cap\del_1 M)}
    &=& \id_{\del_1\R_+^2}\times M(\ast\times\id_0),\\
  C(\ast\times\id_1)|_{\{(0,0)\}\times\del_1 M}
    &=& \id_{\{(0,0)\}}\times M(\ast\times\id_1),\\
  C(\id_1\times\ast)|_{\{(0,0)\}\times\del_0 M}
    &=& \id_{\{(0,0)\}}\times M(\id_1\times\ast).
\end{eqnarray}
The embedding $C(\id_1\times\ast)\colon\del_1\R_+^2\times\del_0 M\to
\{(0,0)\}\times M$ provides us with a diffeomorphism from
$([0,\epsilon]\times\{0\})\times\del_0 M\subseteq
([0,\infty)\times\{0\})\times\del_0 M=\del_1 \R_+^2\times\del_0 M$
onto a submanifold of $\{(0,0)\}\times M$ for some $\epsilon>0$. It
restricts to an inclusion on $\{(0,0)\}\times\del_0 M$ and thereby
yields a (smooth) collar neighbourhood for $\del_0 M$.

\subsection{Open-closed cobordisms}

For a topological space $M$, we denote by $\Pi_0(M)$ the set of its
connected components, and for $p\in M$, we denote by
$[p]\in\Pi_0(M)$ its component.

\subsubsection{Cobordisms}

\begin{defn}
Let $S$ be some set. An $S$-\emph{coloured open-closed cobordism}
$(M,\gamma)$ is a compact oriented $2$-dimensional
$\left<2\right>$-manifold $M$ whose distinguished faces we denote
by $(\del_0 M,\del_1 M)$, together with a map
$\gamma\colon\Pi_0(\del_1 M)\to S$. The face $\del_0 M$ is called
the \emph{black boundary}, $\del_1 M$ the \emph{coloured
boundary}, and $\gamma$ the \emph{colouring}. An \emph{open-closed
cobordism} is an $S$-coloured open-closed cobordism for which $S$
is a one-element set.

Two $S$-coloured open-closed cobordisms $(M,\gamma_M)$ and
$(N,\gamma_N)$ are considered \emph{equivalent} if there is an
orientation preserving diffeomorphism of $2$-dimensional
$\left<2\right>$-manifolds $f\colon M\to N$ that restricts to the
identity on $\del_0 M$ and that preserves the colouring, \ie\
$\gamma_N\circ f=\gamma_M$. We denote the equivalence of open-closed
cobordisms by `$\cong$' both in formulas and in diagrams.

The face $\del_0 M$ is a compact $1$-manifold with boundary and
therefore diffeomorphic to a free union of circles $S^1$ and unit
intervals $[0,1]$. For each of the unit intervals, there is thus
an orientation preserving diffeomorphism
$\phi\colon[0,1]\to\phi([0,1])\subseteq\del_0 M$ onto a component
of $\del_0 M$ such that the boundary points are mapped to the
corners, \ie\ $\phi(\{0,1\})\subseteq \del_0 M\cap\del_1 M$. We
say that the cobordism $(M,\gamma)$ \emph{equips} the unit
interval $[0,1]$ \emph{with the colours} $(\gamma_+,\gamma_-)\in
S\times S$ if $\gamma_+:=\gamma([\phi(1)])$ and
$\gamma_-:=\gamma([\phi(0)])$.
\end{defn}

\subsubsection{Gluing}

Let $(M,\gamma_M)$ and $(N,\gamma_N)$ be $S$-coloured open-closed
cobordisms and $f\colon {S^1}^\ast\to M$ and $g\colon S^1\to N$ be
orientation preserving diffeomorphisms onto components of $\del_0
M$ and $\del_0 N$, respectively. Here we have equipped the circle
$S^1$ with a fixed orientation, and ${S^1}^\ast$ denotes the one
with opposite orientation. Then we obtain an $S$-coloured
open-closed cobordism $M\pushout{f}{}{g}N$ by \emph{gluing $M$ and
$N$ along $S^1$} as follows. As a topological manifold, it is the
pushout. As mentioned in Section~\ref{sect_2mfd}, $\del_0 M$ and
thereby all its components have smooth collar neighbourhoods, and
so the standard techniques are available to equip
$M\pushout{f}{}{g}N$ with the structure of a manifold with corners
whose smooth structure is unique up to a diffeomorphism that
restricts to the identity on $\del_0 M\cup\del_0 N$. It is obvious
that $M\pushout{f}{}{g}N$ also has the structure of a
$\left<2\right>$-manifold with $\del_1 (M\pushout{f}{}{g}N)=\del_1
M\cup\del_1 N$ and, furthermore, that of an $S$-coloured
open-closed cobordism.

Similarly, let $f\colon{[0,1]}^\ast\to M$ and $g\colon[0,1]\to N$
be orientation preserving diffeomorphisms onto components of
$\del_0 M$ and $\del_0 N$, respectively, such that $(M,\gamma_M)$
equips the interval $f({[0,1]}^\ast)$ with the colours
$(\gamma_+,\gamma_-)\in S$ and $(N,\gamma_N)$ equips the interval
$g([0,1])$ precisely with the colours $(\gamma_-,\gamma_+)$. Then
we obtain the gluing of $M$ and $N$ along $[0,1]$ again as the
pushout $M\pushout{f}{}{g}N$ equipped with the smooth structure
that is unique up to a diffeomorphism which restricts to the
identity on $\del_0 M\cup\del_0 N$. It is easy to see that
$M\pushout{f}{}{g}N$ also has the structure of a
$\left<2\right>$-manifold with $\del_1 (M\pushout{f}{}{g}N)=\del_1
M\cup\del_1 N$ and, moreover, due to the matching of colours, that
of an $S$-coloured open-closed cobordism.

\subsubsection{The category $\twocob(S)$}

The following definition of the category of open-closed cobordisms is
inspired by that of Baas, Cohen and Ram{\'\i}rez~\cite{BCR}. What we
call $\twocob(S)$ in the following is in fact a skeleton of the
category of open-closed cobordisms. For this reason, we choose
particular embedded manifolds $C_{\vec{n}}$ as the objects of
$\twocob(S)$. Although these are embedded manifolds, our cobordisms
are not, and we consider two cobordisms equivalent once they are
related by an orientation preserving diffeomorphism that restricts to
the identity on their black boundaries.

\begin{defn}
Let $S$ be a set. The category $\twocob(S)$ is defined as follows. Its
objects are triples $(\vec{n},\gamma_+,\gamma_-)$ consisting of a
finite sequence $\vec{n}:=(n_1,\ldots,n_k)$, $k\in\N_0$, with
$n_j\in\{0,1\}$, $1\leq j\leq k$, and maps
$\gamma_\pm\colon\{1,\ldots,k\}\to S\cup\{\emptyset\}$ for which
$\gamma_\pm(j)\neq\emptyset$ if $n_j=1$ and $\gamma_\pm(j)=\emptyset$
if $n_j=0$\footnote{This is done simply because the values
$\gamma_\pm(j)$ are never used if $n_j=0$, but we nevertheless want to
keep the indices $j$ in line with those of the sequence
$\vec{n}$.}. We denote the \emph{length} of such a sequence by
$|\vec{n}|:=k$.

Each sequence $\vec{n}=(n_1,\ldots,n_k)$ represents the diffeomorphism
type of a compact oriented $1$-dimensional submanifold of $\R^2$,
\begin{equation}
\label{eq_bdrymfd}
  C_{\vec{n}} := \bigcup_{j=1}^k I(j,n_j),
\end{equation}
where $I(j,0)$ is the circle of radius $1/4$ centred at $(j,0)\in\R^2$
and $I(j,1)=[j-1/4,j+1/4]\times\{0\}$, both equipped with the induced
orientation. Taking the disjoint union of two such manifolds
$C_{\vec{n}}$ and $C_{\vec{m}}$ is done as follows,
\begin{equation}
  C_{\vec{n}}\coprod C_{\vec{m}} := C_{\vec{n}}\cup T_{(|\vec{n}|,0)}(C_{\vec{m}}),
\end{equation}
where $T(x,y)\colon\R^2\to\R^2$ denotes the translation by $(x,y)$ in
$\R^2$.

A morphism $f\colon(\vec{n},\gamma_+,\gamma_-)\to(\vec{n}^\prime,
\gamma_+^\prime,\gamma_-^\prime)$ is a pair $f=([f],\Phi)$ consisting
of an equivalence class $[f]$ of $S$-coloured open-closed cobordisms
and a specified orientation preserving diffeomorphism
\begin{equation}
\label{eq_bdrydiff}
  \Phi\colon C_{\vec{n}}^\ast\coprod C_{\vec{n}^\prime}\to\del_0 f,
\end{equation}
such that the following conditions hold:
\begin{enumerate}
\item
  For each $j\in\{1,\ldots,|\vec{n}|\}$ for which $n_j=1$, the
  $S$-coloured open-closed cobordism $(f,\gamma)$ representing $[f]$
  equips the corresponding unit interval $\Phi(I(j,n_j))$ with the
  colours $(\gamma_+(j),\gamma_-(j))$.
\item
  For each $j\in\{1,\ldots,|\vec{n}^\prime|\}$ for which
  $n^\prime_j=1$, $(f,\gamma)$ equips the corresponding unit interval
  $\Phi(I(j,n^\prime_j))$ with the colours
  $(\gamma_+^\prime(j),\gamma_-^\prime(j))$.
\end{enumerate}

The composition $g\circ f$ of two morphisms
$f=([f],\Phi_f)\colon(\vec{n},\gamma_+,\gamma_-)\to(\vec{n}^\prime,\gamma_+^\prime,\gamma_-^\prime)$
and
$g=([g],\Phi_g)\colon(\vec{n}^\prime,\gamma_+^\prime,\gamma_-^\prime)\to
(\vec{n}^{\prime\prime},\gamma_+^{\prime\prime},\gamma_-^{\prime\prime})$
is defined as $g\circ f:=([g\coprod_{C_{\vec{n}'}}f],\Phi_{g\circ
f})$. Here $[g\coprod_{C_{\vec{n}'}}f]$ is the equivalence class of
the $S$-coloured open-closed cobordism $g\coprod_{C_{\vec{n}'}} f$
which is obtained by successively gluing $f$ and $g$ along all the
components of $C_{\vec{n}^\prime}$. $\Phi_{g\circ f}\colon
C_{\vec{n}}^\ast\coprod
C_{\vec{n}^{\prime\prime}}\to\del_0(g\coprod_{C_{\vec{n}'}}f)$ is the
obvious orientation preserving diffeomorphism obtained from
restricting $\Phi_f\colon C_{\vec{n}}^\ast\coprod
C_{\vec{n}^\prime}\to\del_0 f$ and $\Phi_g\colon
C_{\vec{n}^\prime}^\ast\coprod C_{\vec{n}^{\prime\prime}}\to\del_0 g$.

For any object $(\vec{n},\gamma_+,\gamma_-)$, the cylinder
$\id_{(\vec{n},\gamma_+,\gamma_-)}:=[0,1]\times C_{\vec{n}}$ forms an
$S$-coloured open-closed cobordism such that $\del_0
\id_{(\vec{n},\gamma_+,\gamma_-)}=C_{\vec{n}}^\ast\coprod
C_{\vec{n}}$. It plays the role of the identity morphism.

The category $\twocob$ is defined as the category $\twocob(S)$ for the
singleton set $S=\{\ast\}$. When we describe the objects of $\twocob$,
we can suppress the $\gamma_+$ and $\gamma_-$ and simply write the
sequences $\vec{n}=(n_1,\ldots,n_{|\vec{n}|})$.
\end{defn}

Examples of morphisms of $\twocob$ are depicted here,
\begin{equation}
\psset{xunit=0.3cm,yunit=0.3cm}
\begin{pspicture}[.2](2,2.5)
  \rput(0.75,0){\multl}
\end{pspicture}
\colon(1,1)\to (1),\qquad
\begin{pspicture}[.2](2,2.5)
  \rput(0.75,0){\comultc}
\end{pspicture}
\colon(0)\to(0,0),\qquad
\begin{pspicture}[.2](2,2.5)
  \rput(0.75,0){\ctl}
\end{pspicture}
\colon (0)\to (1).
\end{equation}
In these pictures, the source of the cobordism is drawn at the top and
the target at the bottom. The morphism depicted in~\eqref{eq_figure}
goes from $(1,0,1,1,1)$ to $(0,1,1,0,0)$.

The concatenation
$\vec{n}\coprod\vec{m}:=(n_1,\ldots,n_{|\vec{n}|},m_1,\ldots,m_{|\vec{m}|})$
of sequences together with the free union of $S$-coloured
open-closed cobordisms, also denoted by $\coprod$, provides the
category $\twocob(S)$ with the structure of a strict symmetric
monoidal category.

Let $k\in\N$. The symmetric group $\cal{S}_k$ acts on the subset of
objects $(\vec{n},\gamma_+,\gamma_-)\in|\twocob(S)|$ for which
$|\vec{n}|=k$. This action is defined by,
\begin{equation}
  \sigma\rhd(\vec{n},\gamma_+,\gamma_-)
  :=((n_{\sigma^{-1}(1)},\ldots,n_{\sigma^{-1}(|\vec{n}|)}),\gamma_+\circ\sigma^{-1},\gamma_-\circ\sigma^{-1}).
\end{equation}
For each object $(\vec{n},\gamma_+,\gamma_-)\in|\twocob(S)|$ and any
permutation $\sigma\in\cal{S}_{|\vec{n}|}$, we define a morphism,
\begin{equation}
\label{eq_permutation}
  \sigma^{(\vec{n},\gamma_+,\gamma_-)}\colon(\vec{n},\gamma_+,\gamma_-)\to
    \sigma\rhd(\vec{n},\gamma_+,\gamma_-),
\end{equation}
by taking the underlying $S$-coloured open-closed cobordism of the
cylinder $\id_{(\vec{n},\gamma_+,\gamma_-)}$, and replacing the
orientation preserving diffeomorphism
\begin{equation}
  \Phi\colon C_{\vec{n}}^\ast\coprod C_{\vec{n}}\to\del_0\id_{(\vec{n},\gamma_+,\gamma_-)}
\end{equation}
by one that has the components of the $C_{\vec{n}}$ for the target
permuted accordingly. For example, for $S=\{\ast\}$,
$\vec{n}=(1,0,0,1)$ and $\sigma=(234)\in\cal{S}_4$ in cycle notation,
we obtain the morphism $\sigma^{(\vec{n})}$ that is depicted
in~\eqref{eq_permute} below. As morphisms of $\twocob(S)$, \ie\ up to
the appropriate diffeomorphism, these cobordisms satisfy,
\begin{equation}
  \tau^{\sigma\rhd(\vec{n},\gamma_+,\gamma_-)}\circ\sigma^{(\vec{n},\gamma_+,\gamma_-)}
  = {(\tau\circ\sigma)}^{(\vec{n},\gamma_+,\gamma_-)}.
\end{equation}
If the source of $\sigma^{(\vec{n},\gamma_+,\gamma_-)}$ is obvious
from the context, we just write $\sigma$.

\subsubsection{Invariants for open-closed cobordisms}

In order to characterize the $S$-coloured open-closed cobordisms
of $\twocob(S)$ topologically, we need the following invariants.
The terminology is taken from Baas, Cohen, and
Ram{\'\i}rez~\cite{BCR}.

\begin{defn}
\label{def_invariants} Let $S$ be a set and
$f=([f],\Phi)\in\twocob(S)[(\vec{n},\gamma_+,\gamma_-),(\vec{n}^\prime,\gamma_+^\prime,\gamma_-^\prime)]$
be a morphism of $\twocob(S)$ from $(\vec{n},\gamma_+,\gamma_-)$
to $(\vec{n}^\prime,\gamma_+^\prime,\gamma_-^\prime)$.
\begin{enumerate}
\item
  The \emph{genus} $g(f)$ is defined to be the genus of the
  topological $2$-manifold underlying $f$.
\item
  The \emph{window number} of $f$ is a map $\omega(f)\colon S\to\N_0$
  such that $\omega(f)(s)$ is the number of components of the face
  $\del_1 f$ that are diffeomorphic to $S^1$ and that are equipped by
  $\gamma\colon\Pi_0(\del_1 f)\to S$ with the colour $s\in S$. In the
  case $S=\{\ast\}$, we write $\omega(f)\in\N_0$ instead of
  $\omega(f)(\ast)$.
\item
  Let $k$ be the number of coefficients of
  $\vec{n}\coprod\vec{n}^\prime$ that are $1$, \ie\ the number of
  components of the face $\del_0 f$ that are diffeomorphic to the unit
  interval. Number these components by $1,\ldots,k$. The \emph{open
  boundary permutation} $(\sigma(f),\gamma_\del(f))$ of $f$ consists
  of a permutation $\sigma(f)\in\cal{S}_k$ and a map
  $\gamma_\del(f)\colon\{1,\ldots,k\}\to S$. We define $\sigma(f)$ as
  a product of disjoint cycles as follows. Consider every component
  $X$ of the boundary $\del f$ of $f$ viewed as a topological
  manifold, for which $X\cap \del_0 f\cap\del_1 f\neq\emptyset$. These
  are precisely the components that contain a corner of $f$. The
  orientation of $f$ induces an orientation of $X$ and thereby defines
  a cycle $(i_1,\ldots,i_\ell)$ where the $i_j\in\{1,\ldots,k\}$ are
  the numbers of the intervals of $\del_0 f$ that are contained in
  $X$. The permutation $\sigma(f)$ is the product of these cycles for
  all such components $X$. The map
  $\gamma_\del(f)\colon\{1,\ldots,k\}\to S$ is defined such that the
  $S$-coloured open-closed cobordism $(f,\gamma)$ representing $[f]$
  equips the interval with the number $j\in\{1,\ldots,k\}$ with the
  colours $(\gamma_\del(j),\gamma_\del(\sigma^{-1}(j)))$.
\end{enumerate}
\end{defn}

For example, the morphism depicted in~\eqref{eq_figure} has $6$
components of its black boundary diffeomorphic to the unit
interval. Its open boundary permutation is
$\sigma(f)=(256)(34)\in\cal{S}_6$ if one numbers the intervals in the
source (top of the diagram) from left to right by $1,2,3,4$ and those
in the target (bottom of the diagram) from left to right by $5,6$.

\subsection{Generators}
\label{sect_generators}

Beginning with this subsection, we restrict ourselves to the case of
$\twocob$ in which there is only one boundary colour, \ie\
$S=\{\ast\}$. We use a generalization of Morse theory to manifolds
with corners in order to decompose each open-closed cobordism into a
composition of open-closed cobordisms each of which contains precisely
one critical point. The components of these form the \emph{generators}
for the morphisms of the category $\twocob$.

For the generalization of Morse theory to manifolds with corners, we
follow Braess~\cite{Br74}. We first summarize the key definitions and
results.

We need a notion of tangent space for a point $p\in\del M$ if $M$ is a
manifold with corners. Every $p\in M$ has a neighbourhood $U\subseteq
M$ which forms a submanifold of $M$ and for which there is a
diffeomorphism $\phi\colon U\to\phi(U)$ onto an open subset
$\phi(U)\subseteq\R_+^n$. Using the fact that
$\phi(p)\in\R^n_+\subseteq\R^n$, we define the tangent space of $p$ in
$M$ as $T_pM:=d\phi^{-1}_p(T_{\phi(p)}\R^n)$, \ie\ just identifying it
with that of $\phi(p)$ in $\R^n$ via the isomorphism $d\phi^{-1}_p$.

\begin{defn}
Let $M$ be a manifold with corners.
\begin{enumerate}
\item
  For each $p\in M$, we define the \emph{inwards pointing tangential
  cone} $C_pM\subseteq T_pM$ as the set of all tangent vectors $v\in
  T_pM$ for which there exists a smooth path
  $\gamma\colon[0,\epsilon]\to M$ for some $\epsilon>0$ such that
  $\gamma(0)=p$ and the one-sided derivative is:
\begin{equation}
  \lim_{t\to 0+}(\gamma(t)-\gamma(0))/t = v.
\end{equation}
\item
  Let $f\colon M\to\R$ be smooth. A point $p\in M$ is called a
  \emph{critical point} and $f(p)\in\R$ its \emph{critical value} if
  the restriction of the derivative $df_p\colon T_pM\to\R$ to the
  inwards pointing tangential cone is not surjective, \ie\ if
\begin{equation}
  df_p(C_pM)\neq\R.
\end{equation}
  The point $p\in M$ is called $(+)$-\emph{critical} if
  $df_p(C_pM)\subseteq\R_+$ and it is called $(-)$-\emph{critical} if
  $df_p(C_pM)\subseteq\R_-:=-\R_+$.
\end{enumerate}
\end{defn}

Note that $df_p\colon T_pM\to\R$ is linear and therefore maps cones to
cones, and so $df_p(C_pM)$ is either $\{0\}$, $\R_+$, $\R_-$, or
$\R$. If $p\in M$ is a critical point, then $df_p(C_pM)$ is either
$\{0\}$, $\R_+$, or $\R_-$. If $p\in M\backslash\del M$, then
$C_pM=T_pM$, and so $p$ is critical if and only if $df_p=0$. If
$p\in\del M$, $c(p)=1$, and $p$ is critical, then the restriction of
$f$ to $\del M$ has vanishing derivative, \ie\ ${d(f|_{\del M})}_p=0$.

\begin{defn}
Let $M$ be a manifold with corners and $f\colon M\to\R$ be a smooth
function.
\begin{enumerate}
\item
  A critical point $p\in M$ of $f$ is called \emph{non-degenerate} if
  the Hessian of $f$ at $p$, restricted to the kernel of $df_p$,
  has full rank, \ie\ if
\begin{equation}
  \det\Hess_p(f)|_{\ker df_p\otimes\ker df_p}\neq 0.
\end{equation}
\item
  The function $f$ is called a \emph{Morse function} if all its
  critical point are non-degenerate.
\end{enumerate}
\end{defn}

Note that if $p\in M\backslash\del M$, then the notion of
non-degeneracy is as usual. If $p\in\del M$, $c(p)=1$, and $p$ is a
non-degenerate critical point, then $p$ is a non-degenerate critical
point of the restriction $f|_{\del M}\colon\del M\to\R$ in the usual
sense. All non-degenerate critical points are isolated~\cite{Br74}.

For our open-closed cobordisms, we need a special sort of Morse
functions that are compatible with the global structure of the
cobordisms.

\begin{defn}
Let $M\in\twocob[\vec{n},\vec{m}]$ be an open-closed cobordism with
source $C_{\vec{n}}$ and target $C_{\vec{m}}$. Here we have suppressed
the diffeomorphisms from $C_{\vec{n}}$ onto a component of $\del_0 M$,
\etc, and we write $M$ for any representative of its equivalence class.
A \emph{special Morse function} for $M$ is a Morse function $f\colon
M\to\R$ such that the following conditions hold:
\begin{enumerate}
\item
  $f(M)\subseteq [0,1]$.
\item
  $f(p)=0$ if and only if $p\in C_{\vec{n}}$, and $f(p)=1$ if and only
  if $p\in C_{\vec{m}}$.
\item
  Neither $C_{\vec{n}}$ nor $C_{\vec{m}}$ contain any critical points.
\item
  The critical points of $f$ have pairwise distinct critical values.
\end{enumerate}
\end{defn}

Using the standard techniques, one shows that every open-closed
cobordism $M\in\twocob[\vec{n},\vec{m}]$ admits a special Morse
function $f\colon M\to\R$. Since $M$ is compact and since all
non-degenerate critical points are isolated, the set of critical
points of $f$ is a finite set. If neither $a\in\R$ nor $b\in\R$ are
critical values of $f$, the pre-image $N:=f^{-1}([a,b])$ forms an
open-closed cobordism with $\del_0 N=f^{-1}(\{a,b\})$ and $\del_1
N=\del_1 M\cap N$. If $[a,b]$ does not contain any critical value of
$f$, then $f^{-1}([a,b])$ is diffeomorphic to the cylinder
$f^{-1}(\{a\})\times[0,1]$.

The following proposition classifies in terms of Morse data the
non-degenerate critical points that can occur on open-closed
cobordisms.

\begin{prop}
Let $M\in\twocob[\vec{n},\vec{m}]$ be a connected open-closed
cobordism and $f\colon M\to\R$ a special Morse function such that
$f$ has precisely one critical point. Then $M$ is equivalent to
one of the following open-closed cobordisms:
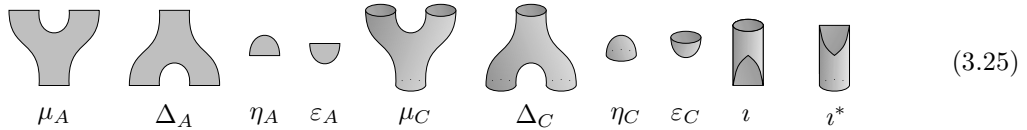
\begin{equation}
\label{eq_generators}
\begin{aligned}
\psset{xunit=.4cm,yunit=.4cm}
\begin{pspicture}[.2](27,3.5)
  \rput(0,1){\multl}
  \rput(0,0){$\mu_A$}
  \rput(4,1){\comultl}
  \rput(4,0){$\Delta_A$}
  \rput(7,2){\birthl}
  \rput(7,0){$\eta_A$}
  \rput(9,2.4){\deathl}
  \rput(9,0){$\epsilon_A$}
  \rput(12,1){\multc}
  \rput(12,0){$\mu_C$}
  \rput(16,1){\comultc}
  \rput(16,0){$\Delta_C$}
  \rput(19,2){\birthc}
  \rput(19,0){$\eta_C$}
  \rput(21,1.6){\deathc}
  \rput(21,0){$\epsilon_C$}
  \rput(23,1){\ctl}
  \rput(23,0){$\imath$}
  \rput(26,1){\ltc}
  \rput(26,0){$\imath^\ast$}
\end{pspicture}
\end{aligned}
\end{equation}
or to one of the compositions
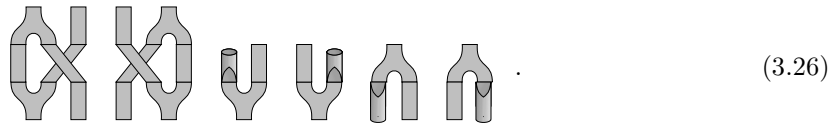
\begin{equation}
\label{eq_saddledecomposed}
\begin{aligned}
\psset{xunit=.2cm,yunit=.2cm}
\begin{pspicture}(33,7.5)
  \rput(1,0){\multl}
  \rput(4,0){\identl}
  \rput(0,2.5){\identl}
  \rput(3,2.5){\crossl}
  \rput(1,5){\comultl}
  \rput(4,5){\identl}
  \rput(7,0){\identl}
  \rput(10,0){\multl}
  \rput(8,2.5){\crossl}
  \rput(11,2.5){\identl}
  \rput(7,5){\identl}
  \rput(10,5){\comultl}
  \rput(15,0){\multl}
  \rput(16,2.5){\identl}
  \rput(13.8,2.5){\ctl}
  \rput(20,0){\multl}
  \rput(19,2.5){\identl}
  \rput(20.8,2.5){\ctl}
  \rput(25,2.5){\comultl}
  \rput(24.2,0){\ltcnew}
  \rput(26,0){\identl}
  \rput(30,2.5){\comultl}
  \rput(29,0){\identl}
  \rput(31.2,0){\ltcnew}
\end{pspicture}
\end{aligned}
.
\end{equation}
All these diagrams show open-closed cobordisms embedded in $\R^3$ and
are drawn in such a way that the vertical axis of the drawing plane is
$-f$. The source is at the top, and the target at the bottom of the
diagram.
\end{prop}

\begin{proof}
We analyze the properties of the non-degenerate critical point $p\in
M$ case by case.
\begin{enumerate}
\item
  If $p\in M\backslash\del M$, then the critical point is
  characterized by its index $i(p)$ (the number of negative
  eigenvalues of $\Hess_p(f)$) as usual; see, for
  example~\cite{Ma97}. There exists a neighbourhood $U\subseteq M$ of
  $p$ and a coordinate system $x\colon U\to\R^2$ in which the Morse
  function has the normal form,
\begin{equation}
  f(p) = -\sum_{j=1}^{i(p)} x_j^2(p) + \sum_{j=i(p)+1}^2 x_j^2(p)
\end{equation}
  for all $p\in U$.
\begin{enumerate}
\item
  If the index is $i(p)=2$, then the Morse function has a maximum at
  $p$, and so the neighbourhood (and thereby the entire open-closed
  cobordism) is diffeomorphic to $\epsilon_C$
  of~\eqref{eq_generators}. Recall that the vertical coordinate of our
  diagrams is $-f$ rather than $+f$.
\item
  If the index is $i(p)=1$, then $f$ has a saddle point. If $M$ were a
  closed cobordism, \ie\ $\del_0M=\del M$, the usual argument would
  show that $M$ is either of the form $\mu_C$ or $\Delta_C$
  of~\eqref{eq_generators}. In the open-closed case, however, the
  saddle can occur in other cases, too, depending on how the boundary
  $\del M$ is decomposed into $\del_0M$ and $\del_1M$. We proceed with
  a case by case analysis and show that in each case, this saddle is
  equivalent to one of the compositions displayed
  in~\eqref{eq_saddledecomposed}:
\begin{eqnarray}
\label{eq_saddledec1}
\begin{aligned}
\psset{xunit=.8cm,yunit=.8cm}
\begin{pspicture}(3.0,2.0)
  \rput(0,0){\bigsaddle}
\end{pspicture}
\end{aligned}
  &\cong&
\begin{aligned}
\psset{xunit=.3cm,yunit=.3cm}
\begin{pspicture}(5,7.5)
  \rput(1,0){\multl}
  \rput(4,0){\identl}
  \rput(0,2.5){\identl}
  \rput(3,2.5){\crossl}
  \rput(1,5){\comultl}
  \rput(4,5){\identl}
\end{pspicture}
\end{aligned}
  ,\\
\label{eq_saddledec2}
\begin{aligned}
\psset{xunit=.7cm,yunit=.7cm}
\begin{pspicture}(3.5,2.0)
  \rput(0,0){\bigcomposed}
\end{pspicture}
\end{aligned}
  &\cong&
\begin{aligned}
\psset{xunit=.3cm,yunit=.3cm}
\begin{pspicture}(3,5)
  \rput(1,0){\multl}
  \rput(0,2.5){\identl}
  \rput(1.8,2.5){\ctl}
\end{pspicture}
\end{aligned}
.
\end{eqnarray}
  Here we show the saddle at the left and the equivalent decomposition
  as a composition and tensor product of the cobordisms
  of~\eqref{eq_generators} with identities on the right. The saddle
  of~\eqref{eq_saddledec1} can appear in two orientations and with the
  intervals in its source and target in any ordering. In any of these
  cases, it is equivalent to one of the first two compositions
  displayed in~\eqref{eq_saddledecomposed}. The saddle
  of~\eqref{eq_saddledec2} can appear flipped upside-down or
  left-right or both, giving rise to the last four compositions
  displayed in~\eqref{eq_saddledecomposed}.

  Note that the equivalences of~\eqref{eq_saddledec1}
  and~\eqref{eq_saddledec2} relate cobordisms whose number of critical
  points differs by an odd number. This is a new feature that dos not
  occur in the case of closed cobordisms.
\item
  If $i(p)=0$, then $f$ has a minimum, and the cobordism is
  diffeomorphic to $\eta_C$ of~\eqref{eq_generators}.
\end{enumerate}
\item
  Otherwise, $p\in\del_1 M\backslash\del_0 M$, \ie\ the critical point
  is on the coloured boundary, but does not coincide with a corner of
  $M$. Consider the restriction $f|_{\del_1 M}\colon\del_1 M\to\R$
  which then has a non-degenerate critical point at $p$ with index
  $i^\prime(p)\in\{0,1\}$.
\begin{enumerate}
\item
  If $i^\prime(p)=1$, then $f|_{\del M}$ has a maximum at $p$.
\begin{enumerate}
\item
  If $p$ is a $(-)$-critical point of $f$, the cobordism is
  diffeomorphic to $\epsilon_A$ of~\eqref{eq_generators}.
\item
  If $p$ is a $(+)$-critical point of $f$, the neighbourhood of $p$
  looks as follows,
\begin{equation}
\begin{aligned}
\psset{xunit=.4cm,yunit=.4cm}
\begin{pspicture}[.2](3,3)
  \pscustom[fillcolor=lightgray,fillstyle=solid,linecolor=white]{
    \psline(1.5,2.5)(1.5,0)
    \psline(-1.5,0)
    \psline(-1.5,2.5)
    \psbezier(-.5,2.5)(-.6,1.5)(0.6,1.5)(.5,2.5)}
  \pscustom{\psbezier(-.5,2.5)(-.6,1.5)(0.6,1.5)(.5,2.5)}
  \psellipse[fillcolor=black,fillstyle=solid](0,1.73)(0.1,0.1)
  \rput(0,1.2){$p$}
  \rput(0.8,0.5){$M$}
\end{pspicture}
\end{aligned}
\end{equation}
  Consider the component of the boundary $\del M$ of $M$ as a
  topological manifold. The set of corners $\del_0 M\cap\del_1 M$
  contains at least two elements. If it contains precisely two
  elements, then the cobordism is diffeomorphic to $\imath^\ast$
  of~\eqref{eq_generators}. Otherwise, it contains six elements, and
  the cobordism is diffeomorphic to $\mu_A$ of~\eqref{eq_generators}.
\item
  Otherwise $p$ is neither $(+)$-critical nor $(-)$-critical, and so
  $df_p=0$. Non-degeneracy now means that $\Hess_p(f)$ is
  non-degenerate. Let $i^{\prime\prime}(p)\in\{0,1,2\}$ be the number
  of negative eigenvalues of $\Hess_p(f)$. The case
  $i^{\prime\prime}(p)=0$ is ruled out by the assumption that
  $i^\prime(p)=1$.
\begin{enumerate}
\item
  If $i^{\prime\prime}(p)=2$, then we are in the same situation as in
  case~2(a)i.
\item
  Otherwise $i^{\prime\prime}(p)=1$, and we are in the same situation
  as in case 2(a)ii.
\end{enumerate}
\end{enumerate}
\item
  If $i^\prime(p)=0$, then $f|_{\del M}$ has a minimum at $p$.
\begin{enumerate}
\item
  If $p$ is a $(+)$-critical point of $f$, the cobordism is
  diffeomorphic to $\eta_A$ of~\eqref{eq_generators}.
\item
  If $p$ is a $(-)$-critical point of $f$, the neighbourhood of $p$
  looks as follows,
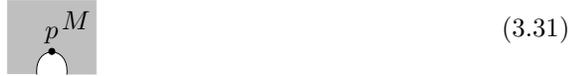
\begin{equation}
\begin{aligned}
\psset{xunit=.4cm,yunit=.4cm}
\begin{pspicture}[.2](3,3)
  \pscustom[fillcolor=lightgray,fillstyle=solid,linecolor=white]{
    \psline(1.5,0)(1.5,2.5)
    \psline(-1.5,2.5)
    \psline(-1.5,0)
    \psbezier(-.5,0)(-.6,1)(0.6,1)(.5,0)}
  \pscustom{\psbezier(-.5,0)(-.6,1)(0.6,1)(.5,0)}
  \psellipse[fillcolor=black,fillstyle=solid](0,0.77)(0.1,0.1)
  \rput(0,1.3){$p$}
  \rput(0.8,1.8){$M$}
\end{pspicture}
\end{aligned}
\end{equation}
  Similarly to case~2(a)ii above, the cobordism is either
  diffeomorphic to $\Delta_A$ or to $\imath$ of~\eqref{eq_generators}.
\item
  Otherwise $df_p=0$, and by considering $\Hess_p(f)$ similarly to
  case~2(a)iii above, we find that we are either in case~2(b)i
  or~2(b)ii.
\end{enumerate}
\end{enumerate}
\end{enumerate}
\end{proof}

The structure of arbitrary open-closed cobordisms can then be
established by using a special Morse function and decomposing the
cobordism into a composition of pieces that have precisely one
critical point each. This result generalizes the conventional handle
decomposition to the case of our sort of $2$-manifolds with corners.

\begin{prop}
\label{propGenerators} Let $f\in\twocob[\vec{n},\vec{n}^\prime]$
be any morphism. Then $[f]=[f_\ell\circ\cdots\circ f_1]$, \ie\ $f$
is equivalent to the composition of a finite number of morphisms
$f_j$ each of which is of the form $f_j=\id_{\vec{m_j}}\coprod
g_j\coprod\id_{\vec{p_j}}$ where $g_j$ is one of the morphisms
depicted in~\eqref{eq_generators} and $\id_{\vec{m_j}}$ and
$\id_{\vec{p_j}}$ are identities, \ie\ cylinders over their
source.
\end{prop}

Our pictures, for example~\eqref{eq_figure}, indicate how the
morphisms are composed from the generators. In order to keep the
height of the diagram small, we have already used relations such as
$(g\coprod\id_{\vec{n}^\prime})\circ(\id_{\vec{m}}\coprod f)=g\coprod
f$ for $f\colon\vec{n}\to\vec{n}^\prime$ and
$g\colon\vec{m}\to\vec{m}^\prime$ which obviously hold in $\twocob$.

Notice that the flat strip, twisted by $2\pi$ when we draw it as
embedded in $\R^3$, is nevertheless equivalent to the flat strip:
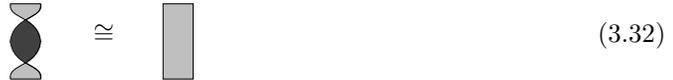
\begin{equation}
\label{doubletwist}
\psset{xunit=.4cm,yunit=.4cm}
\begin{pspicture}[.5](2,2.5)
  \rput(1,0){\twistl}
\end{pspicture}
\quad\cong\quad
\begin{pspicture}[.5](2,2.5)
  \rput(1,0){\identl}
\end{pspicture}
\end{equation}

\subsection{Relations}
\label{sect_relations}

Below we provide a list of relations that the generators of
$\twocob$ satisfy. In Section~\ref{secConsequences}, we summarize
some consequences of these relations.
In Section~\ref{secNormal}, we define a normal form for open-closed
cobordisms with a specified genus, window number and open boundary
permutation. In Section~\ref{secSufficiency}, we then provide an
inductive proof which constructs a finite sequence of
diffeomorphisms that puts an arbitrary open-closed cobordism into
the normal form using only the relations given below. Hence, we
provide a constructive proof that the relations are sufficient to
completely describe the category $\twocob$.


\begin{prop}
\label{PROPrelations}
The following relations hold in the symmetric monoidal category
$\twocob$:
\begin{enumerate}
\item
  The object $\vec{n}=(0)$, \ie\ the circle $C_{\vec{n}}= S^1$,
  forms a commutative Frobenius algebra object.
\begin{gather}
\label{algebrac}
\psset{xunit=.3cm,yunit=.3cm}
\begin{pspicture}[0.5](4,5.5)
  \rput(2,0){\multc}
  \rput(3,2.5){\multc}
  \rput(1,2.5){\curveleftc}
\end{pspicture}
\quad\cong\quad
\begin{pspicture}[0.5](4,5.5)
  \rput(2,0){\multc}
  \rput(1,2.5){\multc}
  \rput(3,2.5){\curverightc}
\end{pspicture}
\qquad\qquad
\begin{pspicture}[0.5](4,5.5)
  \rput(2,0){\multc}
  \rput(3,2.5){\smallidentc}
  \rput(1,2.5){\birthc}
\end{pspicture}
\quad\cong\quad
\begin{pspicture}[0.5](2,5.5)
  \rput(1,0){\identc}
  \rput(1,2.5){\smallidentc}
\end{pspicture}
\quad\cong\quad
\begin{pspicture}[0.5](4,5.5)
  \rput(2,0){\multc}
  \rput(1,2.5){\smallidentc}
  \rput(3,2.5){\birthc}
\end{pspicture}\\
\label{coalgebrac}
\psset{xunit=.3cm,yunit=.3cm}
\begin{pspicture}[0.5](4,5.5)
  \rput(1,0){\comultc}
  \rput(4,0){\curveleftc}
  \rput(2,2.5){\comultc}
\end{pspicture}
\quad\cong\quad
\begin{pspicture}[0.5](4,5.5)
  \rput(0,0){\curverightc}
  \rput(3,0){\comultc}
  \rput(2,2.5){\comultc}
\end{pspicture}
\qquad\qquad
\psset{xunit=.3cm,yunit=.3cm}
\begin{pspicture}[0.5](4,5.5)
  \rput(.8,0){\deathc}
  \rput(3,0){\smallidentc}
  \rput(2,1){\comultc}
\end{pspicture}
\quad\cong\quad
\begin{pspicture}[0.5](2,5.5)
  \rput(1,0){\identc}
  \rput(1,2.5){\smallidentc}
\end{pspicture}
\quad\cong\quad
\begin{pspicture}[0.5](4,5.5)
  \rput(1,0){\smallidentc}
  \rput(2.8,0){\deathc}
  \rput(2,1){\comultc}
\end{pspicture}\\
\label{frobeniusc}
\psset{xunit=.3cm,yunit=.3cm}
\begin{pspicture}[0.5](5,5.5)
  \rput(0,0){\identc}
  \rput(3,0){\multc}
  \rput(1,2.5){\comultc}
  \rput(4,2.5){\identc}
\end{pspicture}
\quad\cong\quad
\begin{pspicture}[0.5](2,5.5)
  \rput(1,0){\comultc}
  \rput(1,2.5){\multc}
\end{pspicture}
\quad\cong\quad
\begin{pspicture}[0.5](5,5.5)
  \rput(5,0){\identc}
  \rput(2,0){\multc}
  \rput(4,2.5){\comultc}
  \rput(1,2.5){\identc}
\end{pspicture}\\
\label{commutavityc}
\psset{xunit=.3cm,yunit=.3cm}
\begin{pspicture}[0.5](2,5.5)
  \rput(1,0){\multc}
  \rput(1,2.5){\crossc}
\end{pspicture}
\quad\cong\quad
\begin{pspicture}[0.5](2,5.5)
  \rput(1,0){\multc}
  \rput(0,2.5){\identc}
  \rput(2,2.5){\identc}
\end{pspicture}
\end{gather}
\item
  The object $\vec{n}=(1)$, \ie\ the interval $C_{\vec{n}}= I$,
  forms a symmetric Frobenius algebra object.
\begin{gather}
\label{algebral}
\psset{xunit=.3cm,yunit=.3cm}
\begin{pspicture}[0.5](4,5.5)
  \rput(2,0){\multl}
  \rput(3,2.5){\multl}
  \rput(1,2.5){\curveleftl}
\end{pspicture}
\quad\cong\quad
\begin{pspicture}[0.5](4,5.5)
  \rput(2,0){\multl}
  \rput(1,2.5){\multl}
  \rput(3,2.5){\curverightl}
\end{pspicture}
\qquad\qquad
\psset{xunit=.3cm,yunit=.3cm}
\begin{pspicture}[0.5](4,5.5)
  \rput(2,0){\multl}
  \rput(3,2.5){\smallidentl}
  \rput(1,2.5){\birthl}
\end{pspicture}
\quad\cong\quad
  \begin{pspicture}[0.5](2,5.5)
  \rput(1,0){\identl}
  \rput(1,2.5){\smallidentl}
\end{pspicture}
\quad\cong\quad
\begin{pspicture}[0.5](4,5.5)
  \rput(2,0){\multl}
  \rput(1,2.5){\smallidentl}
  \rput(3,2.5){\birthl}
\end{pspicture}\\
\label{coalgebral}
\psset{xunit=.3cm,yunit=.3cm}
\begin{pspicture}[0.5](4,5.5)
  \rput(1,0){\comultl}
  \rput(2,2.5){\comultl}
  \rput(4,0){\curveleftl}
\end{pspicture}
\quad\cong\quad
\begin{pspicture}[0.5](4,5.5)
  \rput(0,0){\curverightl}
  \rput(2,2.5){\comultl}
  \rput(3,0){\comultl}
\end{pspicture}
\qquad\qquad
\psset{xunit=.3cm,yunit=.3cm}
\begin{pspicture}[0.5](4,5.5)
  \rput(2,1){\comultl}
  \rput(1,1){\deathl}
  \rput(3,0){\smallidentl}
\end{pspicture}
\quad\cong\quad
\begin{pspicture}[0.5](2,5.5)
  \rput(1,0){\identl}
  \rput(1,2.5){\smallidentl}
\end{pspicture}
\quad\cong\quad
\begin{pspicture}[0.5](4,5.5)
  \rput(2,1){\comultl}
  \rput(1,0){\smallidentl}
  \rput(3,1){\deathl}
\end{pspicture}\\
\label{frobeniusl}
\psset{xunit=.3cm,yunit=.3cm}
\begin{pspicture}[0.5](5,5.5)
  \rput(0,0){\identl}
  \rput(3,0){\multl}
  \rput(1,2.5){\comultl}
  \rput(4,2.5){\identl}
\end{pspicture}
\quad\cong\quad
\begin{pspicture}[0.5](2,5.5)
  \rput(1,0){\comultl}
  \rput(1,2.5){\multl}
\end{pspicture}
\quad\cong\quad
\begin{pspicture}[0.5](5,5.5)
  \rput(5,0){\identl}
  \rput(2,0){\multl}
  \rput(4,2.5){\comultl}
  \rput(1,2.5){\identl}
\end{pspicture}\\
\label{symmetricl}
\psset{xunit=.3cm,yunit=.3cm}
\begin{pspicture}[0.5](4,5.5)
  \rput(2,.5){\deathl}
  \rput(2,.5){\multl}
  \rput(2,3){\crossl}
\end{pspicture}
\quad\cong\quad
\begin{pspicture}[0.5](4,5.5)
  \rput(1,3){\identl}
  \rput(3,3){\identl}
  \rput(2,.5){\deathl}
  \rput(2,.5){\multl}
\end{pspicture}
\end{gather}
\item
  The `zipper'
  $\psset{xunit=.25cm,yunit=.25cm}\begin{pspicture}[.4](1.8,2)\rput(1,0){\ctl}\end{pspicture}$
  forms an algebra homomorphism:
\begin{equation}
\label{ziphomo}
\psset{xunit=.3cm,yunit=.3cm}
\begin{aligned}
\begin{pspicture}[.5](4,5.5)
  \rput(2,0){\multl}
  \rput(2.9,2.5){\ctl}
  \rput(.9,2.5){\ctl}
\end{pspicture}
\quad\cong\quad
\begin{pspicture}[.5](5,5.5)
  \rput(1.9,0){\ctl}
  \rput(2.2,1.97){\multc}
\end{pspicture}
\qquad\qquad
\begin{pspicture}[.5](2,5.5)
  \rput(.9,0){\ctl}
  \rput(1.21,1.98){\birthc}
\end{pspicture}
\quad\cong\quad
\begin{pspicture}[.5](2,5.5)
  \rput(1,1){\birthl}
\end{pspicture}
\end{aligned}
\end{equation}
\item This relation describes the `knowledge' about the centre,
  \cf~\eqref{eq_kfrob1}:
\begin{equation}
\label{zipcenter}
\psset{xunit=.3cm,yunit=.3cm}
\begin{aligned}
\begin{pspicture}[.5](4,5.5)
  \rput(2,0){\multl}
  \rput(.9,2.5){\ctl}
  \rput(3,2.5){\medidentl}
\end{pspicture}
\quad\cong\quad
\begin{pspicture}[.5](4,8)
  \rput(2,0){\multl}
  \rput(2,2.5){\crossl}
  \rput(3,5){\medidentl}
  \rput(.9,5){\ctl}
\end{pspicture}
\end{aligned}
\end{equation}
\item
  The `cozipper'
  $\psset{xunit=.25cm,yunit=.25cm}\begin{pspicture}[.4](1.8,2)\rput(1,0){\ltc}\end{pspicture}$
  is dual to the zipper:
\begin{equation}
\label{eq_zipdual} \psset{xunit=.3cm,yunit=.3cm}
\begin{aligned}
\begin{pspicture}[.5](4,5.5)
  \rput(1.86,0.02){\deathc}
  \rput(2,1){\multc}
  \rput(1,3.5){\medidentc}
  \rput(2.92,3.5){\ltc}
\end{pspicture}
\quad\cong\quad
\begin{pspicture}[.5](4,5.5)
  \rput(2,1){\multl}
  \rput(2,1){\deathl}
  \rput(.9,3.5){\ctl}
  \rput(3,3.5){\medidentl}
\end{pspicture}
\end{aligned}
\end{equation}
\item
  The Cardy condition:
\begin{equation}
\label{cardy}
\psset{xunit=.3cm,yunit=.3cm}
\begin{aligned}
\begin{pspicture}[.4](3,5.5)
  \rput(1.94,1){\ctl}
  \rput(2.2,3){\ltc}
\end{pspicture}
\quad\cong\quad
\begin{pspicture}[.5](4,8)
  \rput(2,0){\multl}
  \rput(2,2.5){\crossl}
  \rput(2,5){\comultl}
\end{pspicture}
\end{aligned}
\end{equation}
\end{enumerate}
\end{prop}

\begin{proof}
It can be show in a direct computation that the depicted
open-closed cobordisms are equivalent. Writing out this proof
would be tremendously laborious, but of rather little insight. For
the Cardy condition~\eqref{cardy}, the right hand side is most
naturally depicted as:
\begin{equation}
\begin{aligned}
\psset{xunit=.25cm,yunit=.25cm}
\begin{pspicture}[.5](8,10)
   \rput(2,0){\multl}
   \rput(1,2.5){\identl}
   \rput(3.2,2.5){\twistl}
   \rput(2,5){\crossl}
   \rput(2,7.5){\comultl}
\end{pspicture}
\end{aligned}
\end{equation}
\end{proof}

\subsection{Consequences of relations}
\label{secConsequences}

In this section, we collect some additional diffeomorphisms that can
be constructed from the diffeomorphisms in
Proposition~\ref{PROPrelations}. To simplify the diagrams, we
define:
\begin{equation}
\label{defpairingl} \psset{xunit=.3cm,yunit=.3cm}
\begin{aligned}
\begin{pspicture}[0.5](2,4)
  \rput(1,1){\zigl}
\end{pspicture}
\quad:=\quad
\begin{pspicture}[0.5](2,4)
  \rput(1,1.5){\deathl}
  \rput(1,1.5){\multl}
\end{pspicture}
\qquad\qquad
\begin{pspicture}[0.5](2,4)
  \rput(1,1.5){\zagl}
\end{pspicture}
\quad:=\quad
\begin{pspicture}[0.5](2,4)
  \rput(1,0.5){\comultl}
  \rput(1,3){\birthl}
\end{pspicture}
\end{aligned}
\end{equation}
These open-closed cobordisms which we sometimes call the \emph{open
pairing}\index{pairing} and \emph{open copairing}, respectively,
satisfy the zig-zag identities:
\begin{equation}
\label{zigzagl} \psset{xunit=.3cm,yunit=.3cm}
\begin{aligned}
\begin{pspicture}[0.5](6,6)
  \rput(2,0){\zigl}
  \rput(5,0){\medidentl}
  \rput(1,2){\identl}
  \rput(3,2){\identl}
  \rput(5,2){\identl}
  \rput(4,4.5){\zagl}
  \rput(1,4.5){\medidentl}
\end{pspicture}
\quad\cong\quad
\begin{pspicture}[0.5](2,6)
  \rput(1,0){\medidentl}
  \rput(1,2){\identl}
  \rput(1,4.5){\medidentl}
\end{pspicture}
\quad\cong\quad
\begin{pspicture}[0.5](6,6)
  \rput(4,0){\zigl}
  \rput(1,0){\medidentl}
  \rput(1,2){\identl}
  \rput(3,2){\identl}
  \rput(5,2){\identl}
  \rput(2,4.5){\zagl}
  \rput(5,4.5){\medidentl}
\end{pspicture}
\end{aligned}
\end{equation}
This follows directly from the Frobenius relations, the left and
right unit laws, and the left and right counit laws. From
Equations~\eqref{symmetricl} and \eqref{algebral}, the pairing can
be shown to be symmetric and invariant,
\begin{equation}
\label{symmetricpairingl} \psset{xunit=.3cm,yunit=.3cm}
\begin{aligned}
\begin{pspicture}[0.5](3,4)
  \rput(2,0){\zigl}
  \rput(2,2){\crossl}
\end{pspicture}
\quad\cong\quad
\begin{pspicture}[0.5](2,4)
  \rput(2,0){\zigl}
  \rput(1,2){\identl}
  \rput(3,2){\identl}
\end{pspicture}
\qquad\qquad
\begin{pspicture}[0.5](4,4)
  \rput(2,0){\zigl}
  \rput(1,2){\multl}
  \rput(3,2){\curverightl}
\end{pspicture}
\quad\cong\quad
\begin{pspicture}[0.5](4,4)
  \rput(2,0){\zigl}
  \rput(3,2){\multl}
  \rput(1,2){\curveleftl}
\end{pspicture}
\end{aligned}
\end{equation}
and the same holds for the copairing. Similarly, we define the
\emph{closed pairing} and the \emph{closed copairing}:
\begin{equation}
\label{defpairingc} \psset{xunit=.3cm,yunit=.3cm}
\begin{aligned}
\begin{pspicture}[0.5](2,4)
  \rput(1,1){\zigc}
\end{pspicture}
\quad:=\quad
\begin{pspicture}[0.5](2,4)
  \rput(0.8,0.5){\deathc}
  \rput(1,1.49){\multc}
\end{pspicture}
\qquad\qquad
\begin{pspicture}[0.5](2,4)
  \rput(1,1.5){\zagc}
\end{pspicture}
\quad:=\quad
\begin{pspicture}[0.5](2,4)
  \rput(1,0.5){\comultc}
  \rput(1,3){\birthc}
\end{pspicture}
\end{aligned}
\end{equation}
These also satisfy the zig-zag identities,
\begin{equation}
\label{zigzagc} \psset{xunit=.3cm,yunit=.3cm}
\begin{aligned}
\begin{pspicture}[0.5](6,6)
  \rput(1.8,0){\zigc}
  \rput(5,0){\medidentc}
  \rput(1,2){\identc}
  \rput(3,2){\identc}
  \rput(5,2){\identc}
  \rput(4,4.5){\zagc}
  \rput(1,4.5){\medidentc}
\end{pspicture}
\quad\cong\quad
\begin{pspicture}[0.5](2,6)
  \rput(1,0){\medidentc}
  \rput(1,2){\identc}
  \rput(1,4.5){\medidentc}
\end{pspicture}
\quad\cong\quad
\begin{pspicture}[0.5](6,6)
  \rput(3.8,0){\zigc}
  \rput(1,0){\medidentc}
  \rput(1,2){\identc}
  \rput(3,2){\identc}
  \rput(5,2){\identc}
  \rput(2,4.5){\zagc}
  \rput(5,4.5){\medidentc}
\end{pspicture}
\end{aligned}
\end{equation}
and the closed pairing is symmetric and invariant,
\begin{equation}
\label{symmetricpairingc} \psset{xunit=.3cm,yunit=.3cm}
\begin{aligned}
\begin{pspicture}[0.5](3,4)
  \rput(1.8,0){\zigc}
  \rput(2,2){\crossc}
\end{pspicture}
\quad\cong\quad
\begin{pspicture}[0.5](2,4)
  \rput(1.8,0){\zigc}
  \rput(1,2){\identc}
  \rput(3,2){\identc}
\end{pspicture}
\qquad\qquad
\begin{pspicture}[0.5](4,4)
  \rput(1.8,0){\zigc}
  \rput(1,2){\multc}
  \rput(3,2){\curverightc}
\end{pspicture}
\quad\cong\quad
\begin{pspicture}[0.5](4,4)
  \rput(1.8,0){\zigc}
  \rput(3,2){\multc}
  \rput(1,2){\curveleftc}
\end{pspicture}
\end{aligned}
\end{equation}
A similar result holds for the closed copairing.

\begin{prop}
\label{propzipdual} The following open-closed cobordisms are
equivalent:
\begin{gather}
\label{zipdual0} \psset{xunit=.3cm,yunit=.3cm}
\begin{pspicture}[.5](4,4)
  \rput(1.8,0.01){\zigc}
  \rput(1,2){\medidentc}
  \rput(2.92,2){\ltc}
\end{pspicture}
\quad\cong\quad
\begin{pspicture}[.5](4,4)
  \rput(2,0.01){\zigl}
  \rput(.9,2){\ctl}
  \rput(3,2){\medidentl}
\end{pspicture}\\
\label{zipdualI} \psset{xunit=.3cm,yunit=.3cm}
\begin{pspicture}[.5](4,4)
  \rput(1.8,0.01){\zigc}
  \rput(3,2){\medidentc}
  \rput(.92,2){\ltc}
\end{pspicture}
\quad\cong\quad
\begin{pspicture}[.5](4,4)
  \rput(2,0.01){\zigl}
  \rput(2.9,2){\ctl}
  \rput(1,2){\medidentl}
\end{pspicture}\\
\label{zipcodualI} \psset{xunit=.3cm,yunit=.3cm}
\begin{pspicture}[.5](4,3.5)
  \rput(1,0){\medidentl}
  \rput(3.1,0){\ltc}
  \rput(2,1.99){\zagl}
\end{pspicture}
\quad\cong\quad
\begin{pspicture}[.5](4,3.5)
  \rput(.7,0){\ctl}
  \rput(3,0){\medidentc}
  \rput(2,1.99){\zagc}
\end{pspicture}\\
\label{zipcodualII} \psset{xunit=.3cm,yunit=.3cm}
\begin{pspicture}[.5](4,3.5)
  \rput(3,0){\medidentl}
  \rput(1.1,0){\ltc}
  \rput(2,1.99){\zagl}
\end{pspicture}
\quad\cong\quad
\begin{pspicture}[.5](4,3.5)
  \rput(2.7,0){\ctl}
  \rput(1,0){\medidentc}
  \rput(2,1.99){\zagc}
\end{pspicture}
\end{gather}
\end{prop}

\begin{proof}
Equation~\eqref{zipdual0} is just a restatement of the second axiom
in Equation~\eqref{eq_zipdual}. The proof of
Equation~\eqref{zipdualI} is as follows:
\begin{gather}
\psset{xunit=.3cm,yunit=.3cm}
\begin{aligned}
\begin{pspicture}[.5](4,4)
  \rput(1.8,0.01){\zigc}
  \rput(3,2){\medidentc}
  \rput(.92,2){\ltc}
\end{pspicture}
\quad
  \xy (0,2)*{\cong};
      (0,-2)*{\eqref{commutavityc}};
  \endxy\quad
\begin{pspicture}[.5](4,7)
  \rput(1.8,0.01){\zigc}
  \rput(2,2){\crossc}
  \rput(3,4.5){\medidentc}
  \rput(.92,4.5){\ltc}
\end{pspicture}
\quad
  \xy (0,2)*{\cong};
      (0,-2)*{{\rm Nat}};
  \endxy\quad
\begin{pspicture}[.5](4,7)
  \rput(1.8,0.01){\zigc}
  \rput(1,2){\medidentc}
  \rput(2.92,2){\ltc}
  \rput(2,4){\crossmixlc}
\end{pspicture}
\quad
  \xy (0,2)*{\cong};
      (0,-2)*{\eqref{eq_zipdual}};
  \endxy\quad
\begin{pspicture}[.5](4,7)
  \rput(2,0.01){\zigl}
  \rput(.9,2){\ctl}
  \rput(3,2){\medidentl}
  \rput(2.2,4){\crossmixlc}
\end{pspicture}
\quad
  \xy (0,2)*{\cong};
      (0,-2)*{{\rm Nat}};
  \endxy\quad
\begin{pspicture}[.5](4,7)
  \rput(2,0.01){\zigl}
  \rput(2.9,4.5){\ctl}
  \rput(1,4.5){\medidentl}
  \rput(2,2){\crossl}
\end{pspicture}
\quad
  \xy (0,2)*{\cong};
      (0,-2)*{\eqref{symmetricl}};
  \endxy\quad
\begin{pspicture}[.5](4,4)
  \rput(2,0.01){\zigl}
  \rput(2.9,2){\ctl}
  \rput(1,2){\medidentl}
\end{pspicture}
\end{aligned}
\end{gather}
By `Nat' we have denoted the obvious diffeomorphisms which,
algebraically speaking, express the naturality of the symmetric
braiding. The proof of Equation~\eqref{zipcodualI} is as follows:
\begin{gather}
\psset{xunit=.3cm,yunit=.3cm}
\begin{aligned}
\begin{pspicture}[.5](4,3.5)
  \rput(1,0){\medidentl}
  \rput(3.1,0){\ltc}
  \rput(2,1.99){\zagl}
\end{pspicture}
\quad
  \xy (0,2)*{\cong};
      (0,-2)*{\eqref{zigzagc}};
  \endxy\quad
\begin{pspicture}[.5](8,5.5)
  \rput(4,0){\zigc}
  \rput(5.2,2){\medidentc}
  \rput(7.2,2){\medidentc}
  \rput(6.2,4){\zagc}
  \rput(1,2){\medidentl}
  \rput(3.1,2){\ltc}
  \rput(2,3.99){\zagl}
\end{pspicture}
\quad
  \xy (0,2)*{\cong};
      (0,-2)*{\eqref{zipdualI}};
  \endxy\quad
\begin{pspicture}[.5](8,7)
  \rput(4,0.01){\zigl}
  \rput(1,2){\medidentl}
  \rput(4.9,2){\ctl}
  \rput(2,3.99){\zagl}
  \rput(3,2){\medidentl}
  \rput(7.2,2){\medidentc}
  \rput(6.2,4){\zagc}
\end{pspicture}
\quad
  \xy (0,2)*{\cong};
      (0,-2)*{\eqref{zigzagl}};
  \endxy\quad
\begin{pspicture}[.5](4,3.5)
  \rput(.7,0){\ctl}
  \rput(3,0){\medidentc}
  \rput(2,1.99){\zagc}
\end{pspicture}
\end{aligned}
\end{gather}
We leave the proof of Equation~\eqref{zipcodualII} as an exercise
for the reader.
\end{proof}

\begin{prop}
\label{propcomultALT} The following open-closed cobordisms are
equivalent:
\begin{equation}
\label{comultpairingIc} \psset{xunit=.3cm,yunit=.3cm}
\begin{pspicture}[.5](6,6)
  \rput(5.08,2){\identc}
  \rput(2.08,2){\multc}
  \rput(4.08,4.49){\zagc}
\end{pspicture}
\quad\cong\quad
\begin{pspicture}[.5](3,4)
  \rput(2,1){\comultc}
\end{pspicture}
\quad\cong\quad
\begin{pspicture}[.5](6,6)
  \rput(1.08,2){\identc}
  \rput(4.08,2){\multc}
  \rput(2.08,4.49){\zagc}
\end{pspicture}
\end{equation}
\begin{equation}
\label{comultpairingIIc} \psset{xunit=.3cm,yunit=.3cm}
\begin{pspicture}[.5](3,4)
  \rput(2,1){\comultc}
\end{pspicture}
\quad\cong\quad
\begin{pspicture}[.5](10,7)
  \rput(2.88,.01){\zigc}
  \rput(7.08,2){\identc}
  \rput(9.08,2){\identc}
  \rput(2.08,2){\curveleftc}
  \rput(4.08,2){\multc}
  \rput(6.08,4.49){\zagc}
  \rput(6.08,4.49){\bigzagc}
\end{pspicture}
\end{equation}
\begin{equation}
\label{comultpairingIl} \psset{xunit=.3cm,yunit=.3cm}
\begin{pspicture}[.5](6,6)
  \rput(5.08,2){\identl}
  \rput(2.08,2){\multl}
  \rput(4.08,4.49){\zagl}
\end{pspicture}
\quad\cong\quad
\begin{pspicture}[.5](3,4)
  \rput(2,1){\comultl}
\end{pspicture}
\quad\cong\quad
\begin{pspicture}[.5](6,6)
  \rput(1.08,2){\identl}
  \rput(4.08,2){\multl}
  \rput(2.08,4.49){\zagl}
\end{pspicture}
\end{equation}
\begin{equation}
\label{comultpairingIIl} \psset{xunit=.3cm,yunit=.3cm}
\begin{pspicture}[.5](3,4)
  \rput(2,1){\comultl}
\end{pspicture}
\quad\cong\quad
\begin{pspicture}[.5](10,7)
  \rput(3.08,.01){\zigl}
  \rput(7.08,2){\identl}
  \rput(9.08,2){\identl}
  \rput(2.08,2){\curveleftl}
  \rput(4.08,2){\multl}
  \rput(6.08,4.49){\zagl}
  \rput(6.08,4.49){\bigzagl}
\end{pspicture}
\end{equation}
\end{prop}

\begin{proof}
The first diffeomorphism in Equation~\eqref{comultpairingIc} is
constructed from the following sequence of diffeomorphisms:
\begin{gather}
\begin{aligned}
\psset{xunit=.3cm,yunit=.3cm}
\begin{pspicture}[.5](3,4)
  \rput(2,1){\comultc}
\end{pspicture}
\quad
  \xy (0,2)*{\cong};
      (0,-2)*{\eqref{algebrac}};
  \endxy\quad
\begin{pspicture}[.5](4,7.5)
  \rput(2.08,2){\comultc}
  \rput(2.08,4.49){\multc}
  \rput(3.08,6.98){\birthc}
\end{pspicture}
\quad
  \xy (0,2)*{\cong};
      (0,-2)*{\eqref{frobeniusc}};
  \endxy\quad
\begin{pspicture}[.5](6,7.5)
  \rput(5.08,2){\identc}
  \rput(2.08,2){\multc}
  \rput(4.08,4.49){\comultc}
  \rput(4.08,6.98){\birthc}
\end{pspicture}
\quad
  \xy (0,2)*{\cong};
      (0,-2)*{\eqref{defpairingc}};
  \endxy\quad
\begin{pspicture}[.5](6,6)
  \rput(5.08,2){\identc}
  \rput(2.08,2){\multc}
  \rput(4.08,4.49){\zagc}
\end{pspicture}
\end{aligned}
\end{gather}
The second diffeomorphism in Equation~\eqref{comultpairingIc} is
constructed similarly. The diffeomorphism in
Equation~\eqref{comultpairingIIc} is constructed as follows:
\begin{gather}
\begin{aligned}
\psset{xunit=.3cm,yunit=.3cm}
\begin{pspicture}[.5](3,4)
  \rput(2,1){\comultc}
\end{pspicture}
\quad
  \xy (0,2)*{\cong};
      (0,-2)*{\eqref{comultpairingIc}};
  \endxy\quad
\begin{pspicture}[.5](6,6)
  \rput(5.08,2){\identc}
  \rput(2.08,2){\multc}
  \rput(4.08,4.49){\zagc}
\end{pspicture}
\quad
  \xy (0,2)*{\cong};
      (0,-2)*{\eqref{zigzagc}};
  \endxy\quad
\begin{pspicture}[.5](10,7)
  \rput(2.88,.01){\zigc}
  \rput(7.08,2){\identc}
  \rput(9.08,2){\identc}
  \rput(2.08,2){\multc}
  \rput(4.08,2){\curverightc}
  \rput(6.08,4.49){\zagc}
  \rput(6.08,4.49){\bigzagc}
\end{pspicture}
\quad
  \xy (0,2)*{\cong};
      (0,-2)*{\eqref{symmetricpairingc}};
  \endxy\quad
\begin{pspicture}[.5](10,7)
  \rput(2.88,.01){\zigc}
  \rput(7.08,2){\identc}
  \rput(9.08,2){\identc}
  \rput(2.08,2){\curveleftc}
  \rput(4.08,2){\multc}
  \rput(6.08,4.49){\zagc}
  \rput(6.08,4.49){\bigzagc}
\end{pspicture}
\end{aligned}
\end{gather}
The proofs of Equations \eqref{comultpairingIl} and
\eqref{comultpairingIIl} are identical to those above with the
closed cobordisms replaced by their open counterparts.
\end{proof}

\begin{prop}
The cozipper\index{cozipper}
$\psset{xunit=.25cm,yunit=.25cm}\begin{pspicture}[.2](2,2.5)\rput(1,0){\ltc}\end{pspicture}$
is a homomorphism of coalgebras, \ie\
\begin{equation}
\label{counithomo}
\psset{xunit=.3cm,yunit=.3cm}
\begin{aligned}
\begin{pspicture}[.5](2,3.5)
  \rput(.9,0){\deathc}
  \rput(1,1){\ltc}
\end{pspicture}
\quad\cong\quad
\begin{pspicture}[.5](2,2)
  \rput(1,1){\deathl}
\end{pspicture}
\qquad\mbox{and}{\qquad}
\begin{pspicture}[.5](3,5)
  \rput(1,0){\comultc}
  \rput(.93,2.49){\ltcab{}}
\end{pspicture}
\quad\cong\quad
\begin{pspicture}[.5](3,5)
  \rput(1.9,2){\comultl}
  \rput(1,0){\ltc}
  \rput(3,0){\ltc}
\end{pspicture}
\end{aligned}
\end{equation}
\end{prop}

\begin{prop}
Open-closed cobordisms of the form
$\psset{xunit=.2cm,yunit=.2cm}\begin{pspicture}[.5](2,4)\rput(1,0){\ltc}\rput(.57,2){\ctl}\end{pspicture}$
which we sometimes call \emph{closed windows}, can be moved around
freely in any closed diagram. By this we mean that the following
open-closed cobordisms are equivalent,
\begin{equation}
\label{movewholeI} \psset{xunit=.3cm,yunit=.3cm}
\begin{aligned}
\begin{pspicture}[.5](3,7)
  \rput(.92,0){\ltcab{}}
  \rput(.7,2){\ctlab{}}
  \rput(3,0){\medidentc}
  \rput(3,2){\medidentc}
  \rput(2,4){\comultc}
\end{pspicture}
\quad\cong\quad
\begin{pspicture}[.5](3,7)
  \rput(2,0){\comultc}
  \rput(1.92,2.49){\ltc}
  \rput(1.7,4.49){\ctl}
\end{pspicture}
\quad\cong\quad
\begin{pspicture}[.5](3,7)
  \rput(2.92,0){\ltcab{}}
  \rput(2.7,2){\ctlab{}}
  \rput(1,0){\medidentc}
  \rput(1,2){\medidentc}
  \rput(2,4){\comultc}
\end{pspicture}
\end{aligned}
\end{equation}

\begin{equation}
\label{movewholeII} \psset{xunit=.3cm,yunit=.3cm}
\begin{aligned}
\begin{pspicture}[.5](3,7)
  \rput(2,0){\multc}
  \rput(.92,2.49){\ltcab{}}
  \rput(.7,4.49){\ctlab{}}
  \rput(3,2.49){\medidentc}
  \rput(3,4.49){\medidentc}
\end{pspicture}
\quad\cong\quad
\begin{pspicture}[.5](3,7)
  \rput(1.92,0){\ltcab{}}
  \rput(1.7,2){\ctlab{}}
  \rput(2,4){\multc}
\end{pspicture}
\quad\cong\quad
\begin{pspicture}[.5](3,7)
  \rput(2,0){\multc}
  \rput(2.92,2.49){\ltcab{}}
  \rput(2.7,4.49){\ctlab{}}
  \rput(1,2.49){\medidentc}
  \rput(1,4.49){\medidentc}
\end{pspicture}
\end{aligned}
\end{equation}
\begin{equation}
\label{movewholeIII} \psset{xunit=.3cm,yunit=.3cm}
\begin{aligned}
\begin{pspicture}[.5](3,9)
  \rput(2,0){\multc}
  \rput(2,2.49){\comultc}
  \rput(1.92,4.98){\ltc}
  \rput(1.7,6.98){\ctl}
\end{pspicture}
\quad\cong\quad
\begin{pspicture}[.5](3,9)
  \rput(1.92,0){\ltcab{}}
  \rput(1.7,2){\ctlab{}}
  \rput(2,4){\multc}
  \rput(2,6.5){\comultc}
\end{pspicture}
\end{aligned}
\end{equation}
\end{prop}

\begin{prop}
\label{propholel} Open-closed cobordisms of the form
$\psset{xunit=.1cm,yunit=.1cm}\begin{pspicture}[.5](3,4)\rput(1.5,0){\multl}\rput(1.5,2.5){\comultl}\end{pspicture}$
which we sometimes call \emph{open windows}, can be moved around
freely in any open diagram. More precisely, the following
open-closed cobordisms are equivalent,
\begin{equation}
\label{movelholeI} \psset{xunit=.3cm,yunit=.3cm}
\begin{aligned}
\begin{pspicture}[.5](5,7)
  \rput(1,0){\multl}
  \rput(1,2.5){\comultl}
  \rput(4,0){\identl}
  \rput(4,2.5){\curveleftl}
  \rput(2,5){\comultl}
\end{pspicture}
\quad\cong\quad
\begin{pspicture}[.5](3,7)
  \rput(1.5,0){\comultl}
  \rput(1.5,2.49){\multl}
  \rput(1.5,5){\comultl}
\end{pspicture}
\quad\cong\quad
\begin{pspicture}[.5](5,7)
  \rput(4,0){\multl}
  \rput(4,2.5){\comultl}
  \rput(1,0){\identl}
  \rput(1,2.5){\curverightl}
  \rput(3,5){\comultl}
\end{pspicture}
\end{aligned}
\end{equation}
\begin{equation}
\label{movelholeII} \psset{xunit=.3cm,yunit=.3cm}
\begin{aligned}
\begin{pspicture}[.5](5,7)
  \rput(1,2.5){\multl}
  \rput(1,5){\comultl}
  \rput(3,2.5){\curverightl}
  \rput(4,5){\identl}
   \rput(2,0){\multl}
\end{pspicture}
\quad\cong\quad
\begin{pspicture}[.5](3,7)
  \rput(1.5,2.5){\comultl}
  \rput(1.5,0){\multl}
  \rput(1.5,5){\multl}
\end{pspicture}
\quad\cong\quad
\begin{pspicture}[.5](5,7)
  \rput(4,2.5){\multl}
  \rput(4,5){\comultl}
  \rput(2,2.5){\curveleftl}
  \rput(1,5){\identl}
  \rput(3,0){\multl}
\end{pspicture}
\end{aligned}
\end{equation}
\end{prop}


\begin{prop}
\label{propholec} Open-closed cobordisms of the form
$\psset{xunit=.1cm,yunit=.1cm}\begin{pspicture}[.5](3,6)\rput(2,1){\multc}\rput(2,3.5){\comultc}\end{pspicture}$,
also called \emph{genus-one operators}\index{genus-one operator},
can be moved around freely in any closed diagram. More precisely,
\begin{equation}
\label{movecholeI} \psset{xunit=.3cm,yunit=.3cm}
\begin{aligned}
\begin{pspicture}[.5](5,7)
  \rput(1,0){\multc}
  \rput(1,2.5){\comultc}
  \rput(4,0){\identc}
  \rput(4,2.5){\curveleftc}
  \rput(2,5){\comultc}
\end{pspicture}
\quad\cong\quad
\begin{pspicture}[.5](3,7)
  \rput(1.5,0){\comultc}
  \rput(1.5,2.49){\multc}
  \rput(1.5,5){\comultc}
\end{pspicture}
\quad\cong\quad
\begin{pspicture}[.5](5,7)
  \rput(4,0){\multc}
  \rput(4,2.5){\comultc}
  \rput(1,0){\identc}
  \rput(1,2.5){\curverightc}
  \rput(3,5){\comultc}
\end{pspicture}
\end{aligned}
\end{equation}

\begin{equation}
\label{movecholeII} \psset{xunit=.3cm,yunit=.3cm}
\begin{aligned}
\begin{pspicture}[.5](5,7.2)
  \rput(2,0){\multc}
  \rput(1,2.5){\multc}
  \rput(1,5){\comultc}
  \rput(3,2.5){\curverightc}
  \rput(4,5){\identc}
\end{pspicture}
\quad\cong\quad
\begin{pspicture}[.5](3,7.2)
  \rput(1.5,0){\multc}
  \rput(1.5,2.5){\comultc}
  \rput(1.5,5){\multc}
\end{pspicture}
\quad\cong\quad
\begin{pspicture}[.5](5,7.2)
  \rput(3,0){\multc}
  \rput(4,2.5){\multc}
  \rput(4,5){\comultc}
  \rput(2,2.5){\curveleftc}
  \rput(1,5){\identc}
\end{pspicture}
\end{aligned}
\end{equation}
\end{prop}


\begin{prop}
The following open-closed cobordisms are equivalent,
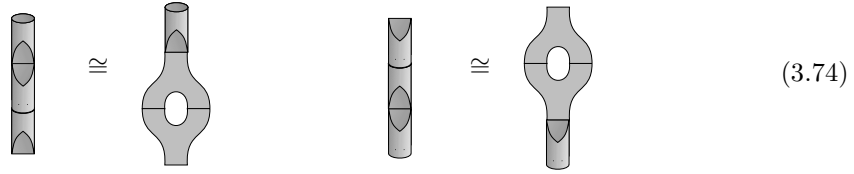
\begin{equation}
\label{symmcoparingl} \psset{xunit=.3cm,yunit=.3cm}
\begin{pspicture}[.5](3,4)
  \rput(2,2.5){\zagl}
  \rput(1,0){\identl}
  \rput(3,0){\identl}
\end{pspicture}
\quad\cong\quad
\begin{pspicture}[.5](3,4)
  \rput(2,2.5){\zagl}
  \rput(2,0){\crossl}
\end{pspicture}
\end{equation}
\begin{equation}
\begin{aligned}
\label{zipcross} \psset{xunit=.3cm,yunit=.3cm}
\begin{pspicture}[.5](3,7)
  \rput(1.4,5){\ctl}
  \rput(1.5,2.5){\comultl}
  \rput(2.5,0){\identl}
  \rput(.5,0){\identl}
\end{pspicture}
\quad\cong\quad
\begin{pspicture}[.5](3,7)
  \rput(1.4,5){\ctl}
  \rput(1.5,2.5){\comultl}
  \rput(1.5,0){\crossl}
\end{pspicture}
\qquad\qquad\qquad \psset{xunit=.3cm,yunit=.3cm}
\begin{pspicture}[.5](3,7)
  \rput(1.6,0){\ltc}
  \rput(.5,4.5){\identl}
  \rput(2.5,4.5){\identl}
  \rput(1.5,2){\multl}
\end{pspicture}
\quad\cong\quad
\begin{pspicture}[.5](3,7)
  \rput(1.6,0){\ltc}
  \rput(1.5,4.5){\crossl}
  \rput(1.5,2){\multl}
\end{pspicture}
\end{aligned}
\end{equation}
\begin{equation}
\label{whitecirc}
\begin{aligned}
\psset{xunit=.3cm,yunit=.3cm}
\begin{pspicture}[.5](3,7)
  \rput(1.5,0){\ctl}
  \rput(1.74,2){\ltc}
  \rput(1.5,4){\ctl}
\end{pspicture}
\quad\cong\quad
\begin{pspicture}[.5](3,8)
  \rput(1.4,5){\ctl}
  \rput(1.5,2.5){\comultl}
  \rput(1.5,0){\multl}
\end{pspicture}
\qquad\qquad\qquad \psset{xunit=.3cm,yunit=.3cm}
\begin{pspicture}[.5](3,7)
  \rput(1.5,0){\ltc}
  \rput(1.26,2){\ctl}
  \rput(1.5,4){\ltc}
\end{pspicture}
\quad\cong\quad
\begin{pspicture}[.5](3,8)
  \rput(1.6,0){\ltc}
  \rput(1.5,4.5){\comultl}
  \rput(1.5,2){\multl}
\end{pspicture}
\end{aligned}
\end{equation}
\end{prop}

\subsection{The normal form of an open-closed cobordism}
\label{secNormal}

In this section, we describe the normal form of an arbitrary connected
open-closed cobordism. This normal form is characterized by its genus,
window number, and open boundary permutation
(\cf~Definition~\ref{def_invariants}). For non-connected open-closed
cobordisms, the normal form has to be taken for each component
independently.

\subsubsection{The case of open source and closed target}

We begin by describing the normal form of a connected open-closed
cobordism whose source consists only of intervals and whose target
consists only of circles. More precisely, we consider those
open-closed cobordisms $f\in\twocob[\vec{n},\vec{m}]$ for which
$\vec{n}=(1,1,\ldots,1)$ and $\vec{m}=(0,0,\ldots,0)$ and denote the
set of all such cobordisms by $\twocob_{\mathrm{O\rightarrow
C}}[\vec{n},\vec{m}]$. Some examples are shown below:
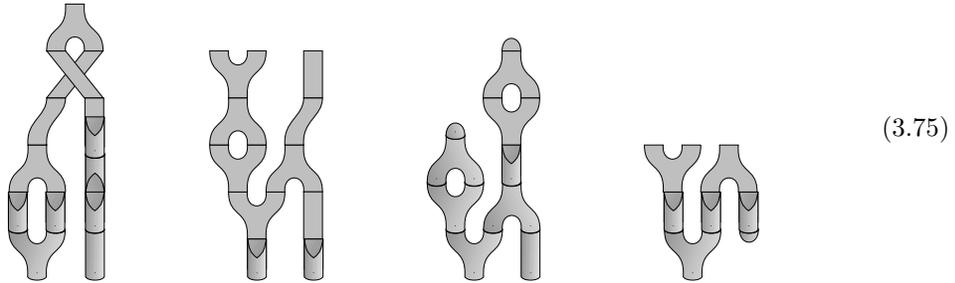
\begin{equation} \label{SampleOpenClosed}
\psset{xunit=.25cm,yunit=.25cm}
\begin{pspicture}[.5](6,15)
    \rput(4.82,7){\identl}
    \rput(5.08,0){\identc}
    \rput(2,0){\multc}
    \rput(5,2.5){\ltc}
    \rput(4.7,4.5){\ctl}
    \rput(0.9,2.5){\ltc}
    \rput(2.9,2.5){\ltc}
    \rput(5,6.5){\ltc}
    \rput(1.78,4.5){\comultl}
    \rput(1.78,7){\curverightl}\
    \rput(3.78,9.5){\crossl}
    \rput(3.78,12){\comultl}
\end{pspicture}
\qquad \qquad
 \psset{xunit=.25cm,yunit=.25cm}
\begin{pspicture}[.5](6,15)
    \rput(1.82,2){\multl}
    \rput(3.82,4.5){\comultl}
    \rput(3.82,7){\curverightl}
    \rput(4.82,9.5){\identl}
    \rput(.82,4.5){\multl}
     \rput(.82,7){\comultl}
      \rput(.82,9.5){\multl}
    \rput(4.82,2){\identl}
    \rput(2,0){\ltc}
    \rput(5,0){\ltc}
\end{pspicture}
\qquad \qquad
  \psset{xunit=.25cm,yunit=.25cm}
\begin{pspicture}[.5](6,15)
    \rput(2,0){\multc}
    \rput(1,2.5){\multc}
    \rput(1,5){\comultc}
    \rput(1,7.5){\birthc}
    \rput(5,0){\identc}
    \rput(4,2.5){\comultc}
    \rput(3.9,5){\ltc}
    \rput(3.75,7){\multl}
    \rput(3.75,9.5){\comultl}
    \rput(3.75,12){\birthl}
\end{pspicture}
\qquad \qquad \psset{xunit=.25cm,yunit=.25cm}
\begin{pspicture}[.5](6,15)
    \rput(2,0){\multc}
    \rput(0.9,2.5){\ltc}
    \rput(2.9,2.5){\ltc}
    \rput(4.76,1.5){\deathc}
    \rput(4.9,2.5){\ltc}
    \rput(.7,4.5){\multl}
    \rput(3.7,4.5){\comultl}
\end{pspicture}
\end{equation}
Once we have defined the normal form for this class of cobordisms, we
describe in Section~\ref{sect_normalall} the normal form of an
arbitrary connected open-closed cobordism by exploiting the
duality on the interval and circle, \cf~\eqref{zigzagl}
and~\eqref{zigzagc}. To provide the reader with some intuition about
the normal form, the cobordisms of~\eqref{SampleOpenClosed} are shown
in normal form below:
\begin{equation}
\qquad
 \psset{xunit=.22cm,yunit=.22cm}
\begin{pspicture}[.5](6,18)
  \rput(1,0){\comultc}
  \rput(1,2.5){\multc}
  \rput(1,5){\comultc}
  \rput(.9,7.5){\ltc}
  \rput(.53,9.5){\ctl}
  \rput(.9,11.5){\ltc}
\end{pspicture}
\qquad \qquad
\begin{pspicture}[.5](6,18)
  \rput(1,0){\comultc}
  \rput(.9,2.5){\ltc}
  \rput(.53,4.5){\ctl}
  \rput(.9,6.5){\ltc}
  \rput(.8,8.5){\multl}
  \rput(-.2,11){\multl}
\end{pspicture}
\qquad \qquad
\begin{pspicture}[.5](6,18)
  \rput(1,0){\comultc}
  \rput(1,2.5){\multc}
  \rput(1,5){\comultc}
  \rput(.9,7.5){\ltc}
  \rput(.53,9.5){\ctl}
  \rput(.9,11.5){\ltc}
  \rput(.53,13.5){\ctl}
  \rput(1,15.5){\birthc}
\end{pspicture}
\qquad \qquad  \psset{xunit=.25cm,yunit=.25cm}
\begin{pspicture}[.5](6,15)
    \rput(2,0){\multc}
    \rput(0.9,2.5){\ltc}
    \rput(2.9,2.5){\ltc}
    \rput(.7,4.5){\multl}
\end{pspicture}
\end{equation}

\begin{defn}
Let $f\in\twocob_{\mathrm{O\rightarrow C}}[\vec{n},\vec{m}]$ be
connected with open boundary permutation $\sigma(f)$, window
number $\omega(f)$, and genus $g(f)$. Write the open boundary
permutation as a product $\sigma(f)=\sigma_1\cdots\sigma_r$,
$r\in\N_0$, of disjoint cycles
$\sigma_j=(i_1^{(j)},\ldots,i_{q_j}^{(j)})$ of length $q_j\in\N$,
$1\leq j\leq r$. The normal form is the composite,
\begin{equation}
\label{eq_prenf}
  \mathrm{NF}_{\mathrm{O\rightarrow C}}(f) :=
    E_{|\vec{m}|}
    \circ D_{g(f)}
    \circ C_{\omega(f)}
    \circ B_r
    \circ\bigl(\coprod_{j=1}^r A(q_j)\bigr)
    \circ\overline{\sigma(f)},
\end{equation}
of the following open-closed cobordisms.
\begin{itemize}
\item
  For each cycle $\sigma_j$, the open-closed cobordism $A(q_j)$
  consists of $q_j-1$ flat multiplications and then a cozipper,
\begin{equation}
  A(q_j) :=
\psset{xunit=.25cm,yunit=.25cm}
\begin{pspicture}[.5](8,14)
  \rput(6,0){\ltc}
  \rput(5.8,2){\multl}
  \rput(4.8,4.5){\multl}
  \rput(2.5,9){$\ddots$}
  \rput(1.5,10){\multl}
\end{pspicture}
\end{equation}
  The normal form~\eqref{eq_prenf} contains the free union of such a
  cobordism for each cycle $\sigma_j$, $1\leq j\leq r$. Cycles of
  length one are represented by a single cozipper. If $|\vec{n}|=0$,
  then we have $r=0$, and the free union is to be replaced by the
  empty set.
\item
  If $r\geq 1$, then the open-closed cobordism $B_r$ consists of $r-1$
  closed multiplications,
\begin{equation}
  B_r :=
\psset{xunit=.25cm,yunit=.25cm}
\begin{pspicture}[.5](10.5,11)
  \rput(5.8,0){\multc}
  \rput(4.8,2.5){\multc}
  \rput(2.5,7){$\ddots$}
  \rput(1.5,8){\multc}
\end{pspicture}
\end{equation}
  and otherwise
$B_0:=\psset{xunit=.25cm,yunit=.25cm}\begin{pspicture}[.5](2,2)\rput(1,1){\birthc}\end{pspicture}$.
\item
  The open-closed cobordism $C_{\omega(f)}$ is defined as,
\begin{equation}
  C_{\omega(f)} := \underbrace{C^\prime\circ C^\prime\circ\cdots\circ C^\prime}_{\omega(f)},\qquad
  C^\prime :=
\psset{xunit=.25cm,yunit=.25cm}
\begin{pspicture}[.5](3,4)
  \rput(1.9,0){\ltc}
  \rput(1.6,2){\ctl}
\end{pspicture}
\end{equation}
  if $\omega(f)\geq 1$ and empty otherwise.
\item
  Similarly,
\begin{equation}
  D_{g(f)} := \underbrace{D^\prime\circ D^\prime\circ\cdots\circ D^\prime}_{g(f)},\qquad
  D^\prime :=
\psset{xunit=.25cm,yunit=.25cm}
\begin{pspicture}[.5](3,5)
  \rput(1.9,0){\multc}
  \rput(1.9,2.5){\comultc}
\end{pspicture}
\end{equation}
  if $g(f)\geq 1$ and empty otherwise.
\item
  $E_{|\vec{m}|}$ consists of $|\vec{m}|-1$ closed comultiplications,
\begin{equation}
  E_{|\vec{m}|} :=
\psset{xunit=.25cm,yunit=.25cm}
\begin{pspicture}[.5](10.5,11)
  \rput(3,8){\comultc}
  \rput(1.9,7){$\vdots$}
  \rput(1,0){\comultc}
  \rput(2,2.5){\comultc}
\end{pspicture}
\end{equation}
  if $|\vec{m}|\geq 1$ and a closed cup
$E_0:=\psset{xunit=.25cm,yunit=.25cm}\begin{pspicture}[.5](2,2)\rput(1,.7){\deathc}\end{pspicture}$
  otherwise.
\item
  Finally, $\overline{\sigma(f)}$ denotes the open-closed cobordism
  that represents the permutation $\overline{\sigma(f)}$ (as defined
  in~\eqref{eq_permutation}) given in the following. Let $\tau(f)$ be
  the open boundary permutation of the open-closed cobordism
\begin{equation}
\label{eq_prenormal}
    E_{|\vec{m}|}
    \circ D_{g(f)}
    \circ C_{\omega(f)}
    \circ B_r
    \circ\bigl(\coprod_{j=1}^r A(q_j)\bigr).
\end{equation}
  Since by construction both $\tau(f)$ and $\sigma(f)$ have the same
  cycle structure, characterized by the partition
  $|\vec{n}|=\sum_{j=1}^r q_j$, there exists a permutation
  $\overline{\sigma(f)}$ such that,
\begin{equation}
  \sigma(f) = {(\overline{\sigma(f)})}^{-1}\cdot\tau(f)\cdot\overline{\sigma(f)}.
\end{equation}
  Note that multiplying $\overline{\sigma(f)}$ by an element in the
  centralizer of $\sigma(f)$ yields the same open-closed cobordism
  $\mathrm{NF}_{\mathrm{O\rightarrow C}}(f)$ up to equivalence
  because of the relations~\eqref{algebral} and~\eqref{zipcross}, and
  so $\mathrm{NF}_{\mathrm{O\rightarrow C}}(f)$ is well defined.
\end{itemize}
\end{defn}

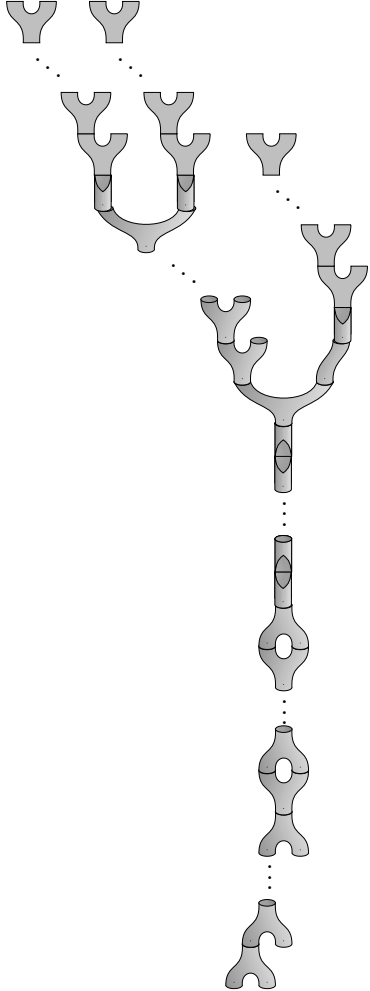
\begin{figure}
\[
\psset{xunit=.22cm,yunit=.22cm}
\begin{pspicture}[.5](5,60)
  \rput(6,8){\comultc}
  \rput(4.9,7){$\vdots$}
  \rput(4,0){\comultc}
  \rput(5,2.5){\comultc}
  \rput(6,10.5){\multc}
  \rput(6,13){\comultc}
  \rput(5.8,17){$\vdots$}
  \rput(6,18){\multc}
  \rput(6,20.5){\comultc}
  \rput(5.9,23){\ltc}
  \rput(5.53,25){\ctl}
  \rput(5.8,29){$\vdots$}
  \rput(5.9,30){\ltc}
  \rput(5.53,32){\ctl}
  \rput(6,34){\hugemultc}
  \rput(8.5,36.5){\curverightc}
  \rput(3.5,39){\topbits}
   \rput(-2.3,36.5){\medbits}
   \rput(-11,47){\topbits}
   \rput(-6,47){\topbits}
\end{pspicture}
\]
\caption{ \label{NormalForm} This figure depicts the normal form
of an open-closed cobordism in $\twocob_{\mathrm{O\rightarrow
C}}[\vec{n},\vec{m}]$ without precomposition with a permutation,
\ie\ it shows the open-closed cobordism~\eqref{eq_prenormal}.}
\end{figure}

When we prove the sufficiency of the relations in
Section~\ref{secSufficiency} below, we provide an algorithm which
automatically produces the required
$\overline{\sigma(f)}$. Figure~\ref{NormalForm} depicts the structure
of the normal form up to the $\overline{\sigma(f)}$, \ie\ it shows a
cobordism of the form~\eqref{eq_prenormal}.

Any cobordism in normal form is invariant (up to equivalence) under
composition with certain permutation morphisms as follows.

\begin{prop}
Let $[f]\in\twocob_{\mathrm{O\rightarrow C}}[\vec{n},\vec{m}]$. Then
\begin{equation}
\label{eq_perminv1}
  [\sigma^{(\vec{m})}\circ \mathrm{NF}_{\mathrm{O\rightarrow C}}(f)]
  = [\mathrm{NF}_{\mathrm{O\rightarrow C}}(f)]
\end{equation}
for any $\sigma\in\cal{S}_{|\vec{m}|}$, and
\begin{equation}
\label{eq_perminv2}
  [\mathrm{NF}_{\mathrm{O\rightarrow C}}(f)\circ\sigma_j^{(\vec{n})}]
  = [\mathrm{NF}_{\mathrm{O\rightarrow C}}(f)]
\end{equation}
for all cycles $\sigma_j\in\cal{S}_{|\vec{n}|}$ from the decomposition
of the open boundary permutation $\sigma(f)=\sigma_1\cdots\sigma_r$
into disjoint cycles.
\end{prop}

\begin{proof}
Equation~\eqref{eq_perminv1} follows from~\eqref{commutavityc},
\eqref{zigzagc}, and~\eqref{symmetricpairingc} while~\eqref{eq_perminv2}
follows from~\eqref{algebral} and~\eqref{symmetricpairingl}.
\end{proof}

\subsubsection{The case of generic source and target}
\label{sect_normalall}

We now extend the definition of the normal form of connected
cobordisms from $\twocob_{{\rm O\rightarrow C}}[\vec{n},\vec{m}]$ to
$\twocob[\vec{n},\vec{m}]$. Let $f$ be a representative of the
equivalence class $[f]$ of an open closed cobordism in
$\twocob[\vec{n},\vec{m}]$. Let $\vec{n}_0=(0,0,\ldots,0)$ and
$\vec{n}_1=(1,1,\ldots,1)$ such that $\vec{n}_0\coprod\vec{n}_1$ is a
permutation of $\vec{n}$, and similarly for $\vec{m}$.

We define a map $\Lambda\maps\twocob[\vec{n},\vec{m}]\to\twocob_{{\rm
O\rightarrow C}}[\vec{m}_1\coprod\vec{n}_1,\vec{m}_0\coprod\vec{n}_0]$
as follows: Let $\sigma_1\in\cal{S}_{|\vec{n}|}$ denote the
permutation that sends $\vec{n}$ to $\vec{n_1}\coprod\vec{n}_0$. Let
$\sigma_2\in\cal{S}_{|\vec{m}|}$ denote the permutation that sends
$\vec{m}$ to $\vec{m_1}\coprod\vec{m}_0$. For example, if
$\vec{n}=(1,0,0,1)$, then $\sigma_1$ is represented as an open-closed
cobordism by:
\begin{equation}
\label{eq_permute}
\psset{xunit=.25cm,yunit=.25cm}
  \sigma_1 \quad = \quad
\begin{pspicture}[.5](8,8)
  \rput(4,2.5){\crossmixcl}
  \rput(1,2.5){\identl}
  \rput(7,2.5){\identc}
  \rput(6,5){\crossmixcl}
  \rput(1,5){\identl}
  \rput(3,5){\identc}
\end{pspicture}
\end{equation}
For $[f]\in\twocob[\vec{n},\vec{m}]$ we define $\Lambda([f])$ to be
the open closed cobordism obtained from $[f]$ by precomposing with
$\sigma_1^{-1}$, postcomposing with $\sigma_2$, gluing closed
copairings on each circle in $\vec{n}_0$, and gluing open pairings on
each interval in $\vec{m}_1$. For example, let $[f]$ be an arbitrary
open closed cobordism from $(1,0,0,1)$ to $(0,1,1,0)$, then
$\Lambda([f])$ is illustrated below:
\begin{equation}
\Lambda\colon
\psset{xunit=.25cm,yunit=.25cm}
\begin{pspicture}[.5](8,4)
  \rput(4,1){\Generalf}
\end{pspicture}
\qquad \mapsto \qquad
\begin{pspicture}[.5](16,16)
 \rput(4,4){\Lambdaf}
\end{pspicture}
\end{equation}
Up to equivalence, this assignment does not depend on the choice
of representative in the class $[f]$; if $f'$ is a different
representative then there exists a black boundary preserving
diffeomorphism from $f$ to $f'$. Applying this diffeomorphism in
the interior of $\Lambda([f])$ shows that $\Lambda([f'])$ is
equivalent to $\Lambda([f])$,
\begin{equation}
  [\Lambda([f])] = [\Lambda([f^\prime])].
\end{equation}
$\Lambda([f])$ is connected if and only if $f$ is.

Given some extra structure, an inverse mapping can be defined. Let $g$
be a representative of a diffeomorphism class in $\twocob_{{\rm
O\rightarrow C}}[\vec{n}',\vec{m}']$ equipped with:
\begin{itemize}
\item
  a decomposition of its source into a free union
  $\vec{n}'=\vec{n}'_{t}\coprod\vec{n}'_{s}$,
\item
  a decomposition of its target into a free union
  $\vec{m}'=\vec{m}'_{t}\coprod\vec{m}'_{s}$,
\item
  an element of the symmetric group $\sigma_1'\in
  \cal{S}_{|\vec{n}_{s}'|+|\vec{m}_s'|}$, and
\item
  an element of the symmetric group $\sigma_2'\in
  \cal{S}_{|\vec{n}_{t}'|+|\vec{m}_t'|}$.
\end{itemize}
Note that the image of an $[f]\in\twocob[\vec{n},\vec{m}]$ under
the mapping $\Lambda$ is equipped with such structure. The
decompositions are given by distinguishing which intervals and
circles came from the source and the target. The permutations can
be taken to be $\sigma_1'=\sigma_1$ and $\sigma_2'=\sigma_2^{-1}$.

We define $\Lambda^{-1}([g])$ to be the open-closed cobordism in
$\twocob[\sigma_2'(\vec{n}_t'\coprod\vec{m}_t'),\sigma_1'(\vec{n}_s'\coprod\vec{m}_s')]$
given by gluing open copairings to the intervals in $\vec{n}_t'$ and
closed pairings to the circles in $\vec{m}_s'$. The result of this
gluing is then precomposed with a cobordism representing $\sigma_1'$
and postcomposed with a cobordism representing $\sigma_2'$.
\begin{equation}
\Lambda^{-1}\colon
\psset{xunit=.25cm,yunit=.25cm}
\begin{pspicture}[.5](8,3)
    \rput(4,1){\GeneralG}
\end{pspicture}
\quad \mapsto \quad
\begin{pspicture}[.5](16,16)
   \rput(2,6){\LambdaINVi}
\end{pspicture}
\end{equation}
Again, this assignment does not depend on the choice of
representative of the class $[g]$. One can readily verify that
this defines a bijection between the equivalence classes of
open-closed cobordisms in $\twocob[\vec{n},\vec{m}]$ and those of
open-closed cobordisms in $\twocob_{{\rm O\rightarrow
C}}[\vec{m}_1\coprod\vec{n}_1,\vec{m}_0\coprod\vec{n}_0]$ equipped
with the extra structure described above. One direction of this
verification,
\begin{equation}
\label{eq_lambdainv}
  [\Lambda^{-1}([\Lambda([f])])] = [f],
\end{equation}
is depicted schematically below:
\begin{equation}
\psset{xunit=.22cm,yunit=.22cm}
\begin{pspicture}[.5](5,3)
 \rput(3,1){\Generalf}
\end{pspicture}
\qquad \cong \qquad
\begin{pspicture}[.5](24,36)
 \rput(0,14){\LambdaLambda}
\end{pspicture}
\end{equation}
In Theorem~\ref{thm_normal} below, we show that for any connected
$[g]\in\twocob_{\mathrm{O\rightarrow C}}[\vec{n},\vec{m}]$, $g$ is
equivalent to its normal form, \ie\
\begin{equation}
  [g] = [\mathrm{NF}_{\mathrm{O\rightarrow C}}(g)].
\end{equation}
Applying this result to $[g]:=\Lambda([f])$ is the motivation for the
definition of the normal form for generic connected
$[f]\in\twocob[\vec{n},\vec{m}]$.

\begin{defn}
\label{def_normal} Let $[f]\in\twocob[\vec{n},\vec{m}]$ be
connected. Then we define its normal form by,
\begin{equation}
  [\mathrm{NF}(f)] :=
  \Lambda^{-1}([\mathrm{NF}_{\mathrm{O\rightarrow C}}\big( \Lambda([f]) \big)]),
\end{equation}
which can be depicted as follows
\begin{equation}
\psset{xunit=.22cm,yunit=.22cm}
\begin{pspicture}[.5](6.5,3)
  \rput(3,1){\GeneralNF}
\end{pspicture}
\quad\cong\quad \psset{xunit=.25cm,yunit=.25cm}
\begin{pspicture}[.5](16,18)
  \rput(3.95,.8){\sigmaTwoINV}
    \rput(5,8){\identl}
    \rput(7,8){\identl}
  \rput(3,7.8){\identl}
  \rput(1,7.8){\identl}
  \rput(3,5.3){\identl}
  \rput(1,5.3){\identl}
  \rput(3,2.8){\identl}
  \rput(1,2.8){\identl}
  \rput(4,10.3){\zagl}
  \rput(4,10.3){\bigzagl}
  \rput(12.0,1.2){\zigc}
  \rput(12.0,1.2){\bigzigc}
  \rput(8,4){\normalf}
  \rput(13.2,3){\identc}
  \rput(15.2,3){\identc}
  \rput(13.2,5.5){\identc}
  \rput(15.2,5.5){\identc}
  \rput(13.2,8){\identc}
  \rput(15.2,8){\identc}
  \rput(11.9,10.5){\sigmaOneINV}
\end{pspicture}
\end{equation}
\end{defn}

\subsection{Proof of sufficiency of relations}
\label{secSufficiency}

In this section, we show that any connected open-closed cobordism
$[f]\in\twocob_{\mathrm{O\rightarrow C}}[\vec{n},\vec{m}]$ can be
related to its normal form $\mathrm{NF}_{\mathrm{O\rightarrow C}}(f)$
by applying the relations of Proposition~\ref{PROPrelations} a finite
number of times. We know that $f$ is equivalent to an open-closed
cobordism of the form stated in Proposition~\ref{propGenerators}.

For convenience, we designate the following composites,
\begin{equation}
\label{extragenerators}
\psset{xunit=.25cm,yunit=.25cm}
\begin{aligned}
\begin{pspicture}[.5](1,4.5)
  \rput(0.5,0){\holel}
\end{pspicture}
\end{aligned}
\quad:=\quad
\begin{aligned}
\begin{pspicture}[.5](1,5.5)
  \rput(0.5,0){\multl}
  \rput(0.5,2.5){\comultl}
\end{pspicture}
\end{aligned}
\qquad\qquad
\begin{aligned}
\begin{pspicture}[.5](1,4.5)
  \rput(0.5,0){\holec}
\end{pspicture}
\end{aligned}
\quad:=\quad
\begin{aligned}
\begin{pspicture}[.5](1,5.5)
  \rput(0.5,0){\multc}
  \rput(0.5,2.5){\comultc}
\end{pspicture}
\end{aligned}
\qquad\qquad
\begin{aligned}
\begin{pspicture}[.5](1,4.5)
  \rput(0.5,0){\holew}
 \end{pspicture}
\end{aligned}
\quad:=\quad
\begin{aligned}
\begin{pspicture}[.5](1,4.5)
  \rput(0.9,0){\ltc}
  \rput(0.5,2){\ctl}
\end{pspicture}
\end{aligned}
\end{equation}
as being distinct generators. This simplifies the proof of the
normal form below. We continue to use the
shorthand~\eqref{defpairingl} and~\eqref{defpairingc}. However, we
do \emph{not} consider these as distinct generators.

In our diagrams, we denote an arbitrary open-closed cobordism $X$,
whose source is a general object $\vec{n}_X$ that contains at least
one $1$, as follows:
\begin{equation}
\psset{xunit=.5cm,yunit=.5cm}
\begin{pspicture}[.5](2,1.5)
  \rput(1,0){\Xl{X}}
\end{pspicture}
\end{equation}
Similarly, to denote an arbitrary open-closed cobordism $Y$, whose
source is a general object $\vec{n}_Y$ containing at least one $0$, we
use the notation:
\begin{equation}
\psset{xunit=.5cm,yunit=.5cm}
\begin{pspicture}[.5](2,1.5)
  \rput(1,0){\Xc{Y}}
\end{pspicture}
\end{equation}
Finally, for an arbitrary open-closed cobordism $Z$, whose target
$\vec{m}$ is not glued to any other cobordism in the decomposition of
$\Sigma$, we use the following notation:
\begin{equation}
\psset{xunit=.5cm,yunit=.5cm}
\begin{pspicture}[.5](2,1.5)
  \rput(1,0){\BottomX{Z}}
\end{pspicture}
\end{equation}
and similarly if the source is not glued to anything.

\begin{defn}
Let $[f]\in\twocob[\vec{n},\vec{m}]$ be written in the form of
Proposition~\ref{propGenerators}. The \emph{height} of a generator
in the decomposition of $f$ is the following number defined
inductively, ignoring all identity morphisms in the decomposition:
\begin{itemize}
\item
  $h(\psset{xunit=.25cm,yunit=.25cm}
     \begin{pspicture}[.5](4,1.7)
       \rput(2,.8){\BottomX{\scs X}}
     \end{pspicture}) := 0$
\item
  $h(\psset{xunit=.25cm,yunit=.25cm}
     \begin{pspicture}[.5](2,1.3)
       \rput(1,1.4){\deathl}
     \end{pspicture}) =
   h(\psset{xunit=.25cm,yunit=.25cm}
     \begin{pspicture}[.5](2,1.3)
       \rput(1,.4){\deathc}
     \end{pspicture}):= 0$
\item
  $h(\psset{xunit=.2cm,yunit=.2cm}
     \begin{pspicture}[.7](4,4)
       \rput(2,0){\Xl{\scs Y}}
       \rput(2,1.5){\multl}
     \end{pspicture})=
   h(\psset{xunit=.2cm,yunit=.2cm}
     \begin{pspicture}[.7](3,3)
       \rput(1.5,0){\Xl{\scs Y}}
       \rput(1.3,1.5){\ctl}
     \end{pspicture}) =
   h(\psset{xunit=.2cm,yunit=.2cm}
     \begin{pspicture}[.5](3,3)
       \rput(1.5,0){\Xl{\scs Y}}
       \rput(1.5,1.5){\birthl}
     \end{pspicture}) =
   h(\psset{xunit=.2cm,yunit=.2cm}
     \begin{pspicture}[.7](4,4)
       \rput(2,0){\Xl{\scs Y}}
       \rput(2,1.5){\holel}
     \end{pspicture}) := 1 + h(Y)$
\item
  $h(\psset{xunit=.2cm,yunit=.2cm}
     \begin{pspicture}[.7](4,4)
       \rput(2,0){\Xc{\scs Y}}
       \rput(2.3,1.5){\multc}
     \end{pspicture}) =
   h(\psset{xunit=.2cm,yunit=.2cm}
     \begin{pspicture}[.7](3,3)
       \rput(1.5,0){\Xc{\scs Y}}
       \rput(1.7,1.5){\ltc}
     \end{pspicture}) =
   h(\psset{xunit=.2cm,yunit=.2cm}
     \begin{pspicture}[.5](3,3)
       \rput(1.5,0){\Xc{\scs Y}}
       \rput(1.8,1.5){\birthc}
     \end{pspicture})=
   h(\begin{pspicture}[.6](4,5)
       \rput(2,0){\Xc{\scs Y}}
       \rput(2,1.5){\holec}
     \end{pspicture})=
   h(\begin{pspicture}[.6](4,5)
       \rput(2,0){\Xc{\scs Y}}
       \rput(2.22,1.5){\holew}
     \end{pspicture}) := 1 + h(Y)$
\item
  $h(\psset{xunit=.2cm,yunit=.2cm}
     \begin{pspicture}[.6](4,4)
       \rput(1,0){\Xc{\scs Y}}
       \rput(3,0){\Xc{\scs Z}}
       \rput(2.3,1.5){\comultc}
     \end{pspicture}) =
   h(\psset{xunit=.2cm,yunit=.2cm}
     \begin{pspicture}[.6](4,4)
       \rput(1,0){\Xl{\scs Y}}
       \rput(3,0){\Xl{\scs Z}}
       \rput(2,1.5){\comultl}
     \end{pspicture}) := h(Y)+h(Z)+1.$
\end{itemize}
\end{defn}

\begin{thm}
\label{thm_normal} Let $[f]\in\twocob_{\mathrm{O\rightarrow
C}}[\vec{n},\vec{m}]$ be a connected open-closed cobordism. Then $f$
is equivalent to its normal form, \ie\
\begin{equation}
  [f] = [\mathrm{NF}_{\mathrm{O\rightarrow C}}(f)].
\end{equation}
\end{thm}

\begin{proof}
We say \emph{decomposition} for a presentation of $f$ as a
composition of the generators as in
Proposition~\ref{propGenerators}, \ie\ for a generalized handle
decomposition. We use the term \emph{move} for the application of a
diffeomorphism from Proposition~\ref{PROPrelations}, we just say
diffeomorphism here meaning diffeomorphism relative to the black
boundary, and we use the term \emph{configuration} of a generator in
a decomposition to refer to the generators immediately pre- and
postcomposed to it.

Employing Proposition~\ref{propGenerators}, let $f$ be given by any
arbitrary decomposition. We construct a diffeomorphism from this
decomposition to $\mathrm{NF}_{\mathrm{O\rightarrow C}}(f)$ by
applying a finite sequence of the moves from
Proposition~\ref{PROPrelations}. This proceeds step by step as
follows.
\begin{enumerate}[I)]
%
%
\item
\label{stepI}
  The decomposition of $f$ is equivalent to one without any
  \emph{open cups}
\begin{math}
\psset{xunit=.25cm,yunit=.25cm}
\begin{pspicture}[.2](2,2)
  \rput(1,1.2){\deathl}
\end{pspicture}
\end{math}
  or \emph{open caps}
\begin{math}
\psset{xunit=.25cm,yunit=.25cm}
\begin{pspicture}[.2](2,2)
  \rput(1,.4){\birthl}
\end{pspicture}
\end{math}. This is achieved by applying the following moves.
\begin{enumerate}[a)]
%
%
\item
\label{stepIa}
\begin{math}
\psset{xunit=.25cm,yunit=.25cm}
\begin{pspicture}[.5](2,2)
  \rput(1,1.4){\deathl}
\end{pspicture}
\quad
  \xy
    {\ar@{|->} (-2,0);(2,0)};
    (0,-4)*{\eqref{counithomo}};
  \endxy\quad
\begin{pspicture}[.5](2,3)
  \rput(1,0){\deathc}
  \rput(1.16,1){\ltc}
\end{pspicture}
\end{math}
%
%
\item
\label{stepIb}
\begin{math}
\psset{xunit=.25cm,yunit=.25cm}
\begin{pspicture}[.5](2,2)
  \rput(1,.8){\birthl}
\end{pspicture}
\quad
  \xy
    {\ar@{|->} (-2,0);(2,0)};
    (0,-4)*{\eqref{ziphomo}};
  \endxy\quad
\begin{pspicture}[.5](2,2.5)
  \rput(1,0){\ctl}
  \rput(1.4,1.99){\birthc}
\end{pspicture}
\end{math}
\end{enumerate}
  to every instance of the open cup and cap.
%
%
\item
\label{stepII}
  The resulting decomposition of $f$ is equivalent to one in
  which every \emph{open comultiplication}
  $\psset{xunit=.2cm,yunit=.2cm}\begin{pspicture}[.2](3,2.5)\rput(1.5,0){\comultl}\end{pspicture}$
  appears in one of the following three configurations:
\begin{equation}
\label{threeforms} \psset{xunit=.25cm,yunit=.25cm}
\begin{pspicture}[.5](6,10)
  \rput(3,0){\multl}
  \rput(2,2.5){\identl}
  \rput(2,5){\identl}
  \rput(5,5){\multl}
  \rput(4,7.5){\widecomultl}
  \rput(4.8,3.7){$?$}
\end{pspicture}
\qquad\qquad
\begin{pspicture}[.5](6,10)
  \rput(5,0){\multl}
  \rput(6,2.5){\identl}
  \rput(6,5){\identl}
  \rput(3,5){\multl}
  \rput(4,7.5){\widecomultl}
  \rput(3.2,3.7){$?$}
\end{pspicture}
\qquad\qquad
\begin{pspicture}[.5](4,5)
  \rput(2,0){\comultl}
  \rput(3,2.5){\comultl}
\end{pspicture}
\end{equation}
  where the `?' may be any open-closed cobordism which may or may not
  be attached to the multiplication at the bottom. We prove this
  step-by-step by considering every possible configuration and
  providing the moves to reduce the decomposition into one of
  the above mentioned configurations.
\begin{enumerate}[a)]
%
%
\item
\label{stepIIa}
  The cases
\begin{math}
\psset{xunit=.25cm,yunit=.25cm}
\begin{pspicture}[.5](3,3)
  \rput(.5,.5){\deathl}
  \rput(1.5,.5){\comultl}
\end{pspicture}
\end{math}
  and
\begin{math}
\psset{xunit=.25cm,yunit=.25cm}
\begin{pspicture}[.5](3,3)
  \rput(2.5,.5){\deathl}
  \rput(1.5,.5){\comultl}
\end{pspicture}
\end{math}
  are excluded by step (\ref{stepIa}).
%
%
\item
\label{stepIIb}
  Wherever possible apply the move:
\begin{math}
\psset{xunit=.25cm,yunit=.25cm}
\begin{pspicture}[.5](5,5)
  \rput(3,0){\comultl}
  \rput(2,2.5){\comultl}
\end{pspicture}
  \xy
    {\ar@{|->} (-2,2);(2,2)};
    (0,-2)*{\eqref{coalgebral}};
  \endxy
\psset{xunit=.25cm,yunit=.25cm}
\begin{pspicture}[.5](5,5)
  \rput(2,0){\comultl}
  \rput(3,2.5){\comultl}
\end{pspicture}
\end{math}.
%
%
\item
\label{stepIIc}
  We consider all of the remaining possible configurations and provide
  a list of moves which either remove the open comultiplication or
  reduce its height. Since there are no longer any open cups after
  (\ref{stepIa}) and since the target of $f$ is of the form
  $\vec{m}=(0,\ldots,0)$, \ie\ a free union of circles, the open
  comultiplication is either removed from the diagram or takes
  the form claimed in \eqref{stepII} before its height is reduced to
  zero.

  Apply the following moves wherever possible:
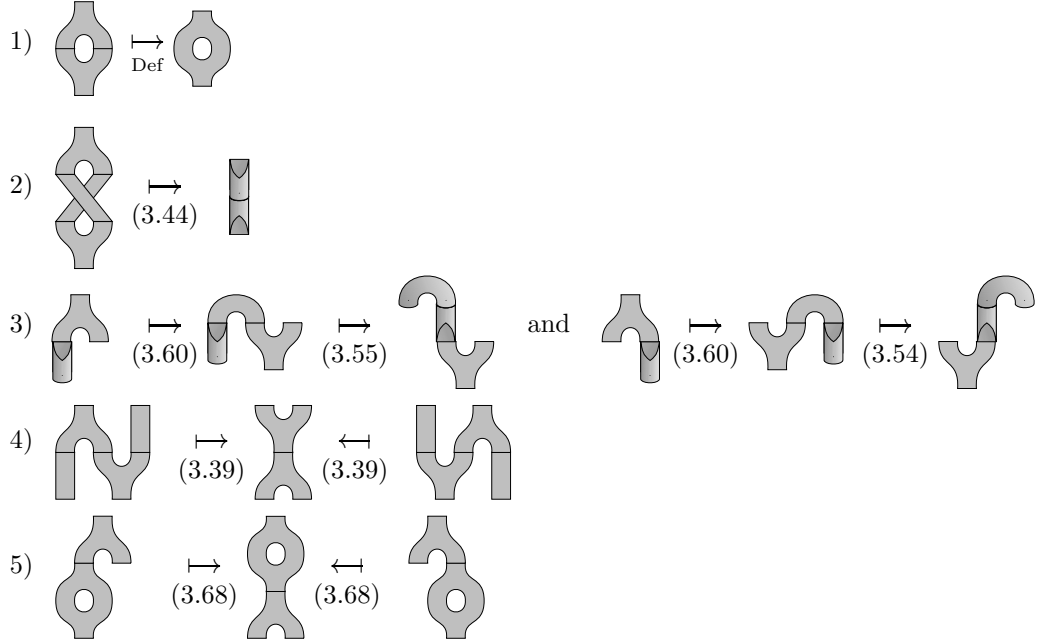
\begin{enumerate}[1)]
%
%
\item
\label{stepII1}
\begin{math}
\psset{xunit=.25cm,yunit=.25cm}
\begin{pspicture}[.5](4.5,5)
  \rput(2,0){\multl}
  \rput(2,2.5){\comultl}
\end{pspicture}
  \xy {\ar@{|->} (-2,0);(2,0)};
      (0,-3)*{{\rm \scs Def}};
  \endxy
\begin{pspicture}[.5](4,4)
  \rput(2,0){\holel}
\end{pspicture}
\end{math}
%
%
\item
\label{stepII2}
\begin{math}
\psset{xunit=.25cm,yunit=.25cm}
\begin{pspicture}[.5](4.5,8)
  \rput(2,0){\multl}
  \rput(2,2.5){\crossl}
  \rput(2,5){\comultl}
\end{pspicture}
  \xy {\ar@{|->} (-2,0);(2,0)};
                 (0,-4)*{\eqref{cardy}};
  \endxy
\begin{pspicture}[.4](3,8)
  \rput(2,1){\ctl}
  \rput(2.32,3){\ltc}
\end{pspicture}
\end{math}
%
%
\item
\label{stepII3}
\begin{math}
\psset{xunit=.25cm,yunit=.25cm}
\begin{pspicture}[.5](4.5,5)
  \rput(1,0){\ltc}
  \rput(1.8,2){\comultl}
\end{pspicture}
  \xy {\ar@{|->} (-2,0);(2,0)};
      (0,-4)*{\eqref{comultpairingIl}};
  \endxy
\begin{pspicture}[.5](6.5,4)
  \rput(4,0){\multl}
  \rput(1.15,.5){\ltc}
  \rput(2,2.5){\zagl}
\end{pspicture}
  \xy {\ar@{|->} (-2,0);(2,0)};
      (0,-4)*{\eqref{zipcodualII}};
  \endxy
\begin{pspicture}[.5](6,6)
  \rput(4.1,0){\multl}
  \rput(2.9,2.5){\ctl}
  \rput(2.3,4.5){\zagc}
\end{pspicture}
  \quad\mbox{and}\quad
\begin{pspicture}[.5](4,5)
  \rput(3,0){\ltc}
  \rput(1.8,2){\comultl}
\end{pspicture}
  \xy {\ar@{|->} (-2,0);(2,0)};
      (0,-4)*{\eqref{comultpairingIl}};
  \endxy
\begin{pspicture}[.5](6.5,4)
  \rput(2,0){\multl}
  \rput(5.15,.5){\ltc}
  \rput(4,2.5){\zagl}
\end{pspicture}
  \xy {\ar@{|->} (-2,0);(2,0)};
      (0,-4)*{\eqref{zipcodualI}};
  \endxy
\begin{pspicture}[.5](6,6)
  \rput(2,0){\multl}
  \rput(2.9,2.5){\ctl}
  \rput(4.3,4.5){\zagc}
\end{pspicture}
\end{math}
%
%
\item
\label{stepII4}
\begin{math}
\psset{xunit=.25cm,yunit=.25cm}
\begin{pspicture}[0.5](7,5.5)
  \rput(1,0){\identl}
  \rput(4,0){\multl}
  \rput(2,2.5){\comultl}
  \rput(5,2.5){\identl}
\end{pspicture}
  \xy {\ar@{|->} (-2,0);(2,0)};
      (0,-4)*{\eqref{frobeniusl}};
  \endxy
\begin{pspicture}[0.5](4,5.5)
  \rput(2,0){\comultl}
  \rput(2,2.5){\multl}
\end{pspicture}
  \xy {\ar@{|->} (2,0);(-2,0)};
      (0,-4)*{\eqref{frobeniusl}};
  \endxy
\psset{xunit=.25cm,yunit=.25cm}
\begin{pspicture}[0.5](6.5,5.5)
  \rput(6,0){\identl}
  \rput(3,0){\multl}
  \rput(5,2.5){\comultl}
  \rput(2,2.5){\identl}
\end{pspicture}
\end{math}
%
%
\item
\label{stepII5}
\begin{math}
\psset{xunit=.25cm,yunit=.25cm}
\begin{pspicture}[0.5](6.6,7)
  \rput(2,0){\holel}
  \rput(3,4){\comultl}
\end{pspicture}
  \xy {\ar@{|->} (-2,0);(2,0)};
      (0,-4)*{\eqref{movelholeI}};
  \endxy
\begin{pspicture}[0.5](4,7)
  \rput(2,0){\comultl}
  \rput(2,2.5){\holel}
\end{pspicture}
  \xy {\ar@{|->} (2,0);(-2,0)};
      (0,-4)*{\eqref{movelholeI}};
  \endxy
\psset{xunit=.25cm,yunit=.25cm}
\begin{pspicture}[0.5](6.5,7)
  \rput(4,0){\holel}
  \rput(3,4){\comultl}
\end{pspicture}
\end{math}
\end{enumerate}
%
%
\item
\label{stepIId}
  Iterate steps~\eqref{stepIIb} and~\eqref{stepIIc}. Since each
  iteration either removes the open comultiplication or reduces its
  height, this process is guaranteed to terminate with every
  comultiplication in one of the three configurations
  of~\eqref{threeforms}.
\end{enumerate}
%
%
\item
\label{stepIII}
  Now we apply a sequence of moves to the decomposition of $f$
 which reduces the number of possible configurations that need
  to be considered.
%
%
\begin{enumerate}[a)]
\item
\label{stepIIIa}
  To begin, we provide a sequence of moves to put every
  \emph{open multiplication}
  $\psset{xunit=.2cm,yunit=.2cm}\begin{pspicture}[.2](3,2.5)\rput(1.5,0){\multl}\end{pspicture}$
  in the decomposition of $f$ into one of the following
  configurations:
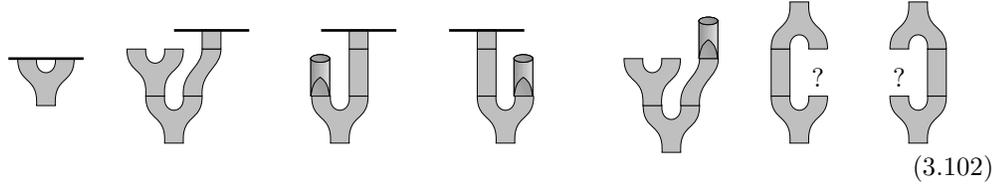
\begin{equation}
\label{formII} \psset{xunit=.25cm,yunit=.25cm}
\begin{pspicture}[0.5](3,4)
  \rput(2,1){\multl}
  \psline[linewidth=1pt](-0,3.5)(4,3.5)
\end{pspicture}
\qquad
\begin{pspicture}[0.5](6,6)
  \rput(3,0){\multl}
  \rput(2,2.5){\multl}
  \rput(4,2.5){\curverightl}
  \rput(5,5){\smallidentl}
  \psline[linewidth=1pt](3,6)(7,6)
\end{pspicture}
\qquad
\begin{pspicture}[0.5](6,6)
  \rput(3,0){\multl}
  \rput(1.8,2.5){\ctl}
  \rput(4,2.5){\identl}
  \rput(4,5){\smallidentl}
  \psline[linewidth=1pt](2,6)(6,6)
\end{pspicture}
\qquad
\begin{pspicture}[0.5](6,6)
  \rput(3,0){\multl}
  \rput(3.8,2.5){\ctl}
  \rput(2,2.5){\identl}
  \rput(2,5){\smallidentl}
  \psline[linewidth=1pt](0,6)(4,6)
\end{pspicture}
\qquad
\begin{pspicture}[0.5](5,7)
  \rput(3,0){\multl}
  \rput(2,2.5){\multl}
  \rput(4,2.5){\curverightl}
  \rput(4.83,5){\ctl}
\end{pspicture}
\qquad
\begin{pspicture}[0.5](3.5,6)
  \rput(2,0){\multl}
  \rput(1,2.5){\identl}
  \rput(2,5){\comultl}
  \rput(3,3.5){$?$}
\end{pspicture}
\qquad
\begin{pspicture}[0.5](4,6)
  \rput(2,0){\multl}
  \rput(3,2.5){\identl}
  \rput(2,5){\comultl}
  \rput(1,3.5){$?$}
\end{pspicture}
\end{equation}
  Again, we prove this claim by considering all possible
  configurations of the open multiplication.

  Apply the following moves which either removes the open
  multiplication or increases its height or attains the desired
  configuration.
%
%
\begin{enumerate}[1)]
\item
\label{stepIIIa1}
\begin{math}
\psset{xunit=.25cm,yunit=.25cm}
\begin{pspicture}[0.5](5,5)
  \rput(2,0){\multl}
  \rput(3,2.5){\multl}
\end{pspicture}
  \xy {\ar@{|->} (-2,0);(2,0)};
      (0,-4)*{\eqref{algebral}};
  \endxy
\begin{pspicture}[0.5](5,5)
  \rput(3,0){\multl}
  \rput(2,2.5){\multl}
\end{pspicture}
\end{math}
%
%
\item
\label{stepIIIa2}
\begin{math}
\psset{xunit=.25cm,yunit=.25cm}
\begin{pspicture}[0.5](5,7)
  \rput(3,0){\multl}
  \rput(2,2.5){\holel}
\end{pspicture}
  \xy {\ar@{|->} (-2,0);(2,0)};
      (0,-4)*{\eqref{movelholeII}};
  \endxy
\begin{pspicture}[0.5](4,7)
  \rput(2,4){\multl}
  \rput(2,0){\holel}
\end{pspicture}
  \xy {\ar@{|->} (2,0);(-2,0)};
      (0,-4)*{\eqref{movelholeII}};
  \endxy
\begin{pspicture}[0.5](5,7)
  \rput(3,0){\multl}
  \rput(4,2.5){\holel}
\end{pspicture}
\end{math}
%
%
\item
\label{stepIIIa3}
\begin{math}
\psset{xunit=.25cm,yunit=.25cm}
\begin{pspicture}[0.5](5,5)
  \rput(2.16,0){\multl}
  \rput(1,2.5){\ctl}
  \rput(3,2.5){\ctl}
\end{pspicture}
  \xy {\ar@{|->} (-2,0);(2,0)};
      (0,-4)*{\eqref{ziphomo}};
  \endxy
\begin{pspicture}[0.5](5,5)
  \rput(3,0){\ctl}
  \rput(3.38,1.99){\multc}
\end{pspicture}
\end{math}
%
%
\item
\label{stepIIIa4}
  The moves of steps \ref{stepII1}, \ref{stepII2}, \ref{stepII4}.
\end{enumerate}
  All other configurations are excluded by step \ref{stepI}. Since
  none of these steps increases the number of generators in the
  decomposition of $f$, iterating this process either removes all open
  multiplications or puts them into the configurations in
  \eqref{formII} as claimed above.
%
%
\item
\label{stepIIIb}
  Now we show that the source of every cozipper
  $\psset{xunit=.25cm,yunit=.25cm}\begin{pspicture}[.4](1.8,2)\rput(1,0){\ltc}\end{pspicture}$
  can be put into either of the following configurations:
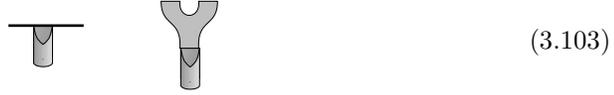
\begin{equation}
\psset{xunit=.25cm,yunit=.25cm}
\begin{pspicture}[.4](2,2)
  \rput(2,0){\ltc}
  \psline[linewidth=1pt](-0,2)(4,2)
\end{pspicture}
\qquad\qquad
\begin{pspicture}[.4](2,5)
  \rput(2.2,0){\ltc}
  \rput(2,2){\multl}
\end{pspicture}
\end{equation}
  We establish the above claim by applying the following
  sequence of moves wherever they are possible.
%
%
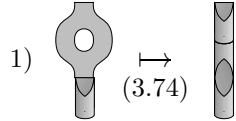
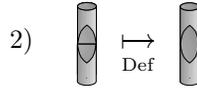
\begin{enumerate}[1)]
\item
\label{stepIIIb2}
\begin{math}
\psset{xunit=.25cm,yunit=.25cm}
\begin{pspicture}[.4](4,6.5)
  \rput(2.2,0){\ltc}
  \rput(2,2){\holel}
\end{pspicture}
  \xy {\ar@{|->} (-2,0);(2,0)};
      (0,-4)*{\eqref{whitecirc}};
  \endxy
\begin{pspicture}[0.4](4,6.5)
  \rput(2,0){\holew}
  \rput(2,4){\ltc}
\end{pspicture}
\end{math}
%
%
\item
\label{stepIIIb1}
\begin{math}
\psset{xunit=.25cm,yunit=.25cm}
\begin{pspicture}[.4](4,4.5)
  \rput(2,2){\ctl}
  \rput(2.3,0){\ltc}
\end{pspicture}
  \xy {\ar@{|->} (-2,0);(2,0)};
      (0,-3)*{{\rm \scs Def}};
  \endxy
\begin{pspicture}[0.4](4,4.5)
  \rput(2,0){\holew}
\end{pspicture}
\end{math}
%
%
\item
\label{stepIIIb3}
  The moves of step~\ref{stepII3}.
\end{enumerate}
  All other configurations are excluded by step \ref{stepI}.
%
%
\item
\label{stepIIIc}
  In this step we show that every instance of the open window
  $\psset{xunit=.25cm,yunit=.25cm}\begin{pspicture}[.4](4,4)\rput(2,0){\holel}\end{pspicture}$
  can be removed. Iterate the
  following sequence of moves wherever possible.
\begin{enumerate}[1)]
%
%
\item
\begin{math}
\psset{xunit=.25cm,yunit=.25cm}
\begin{pspicture}[.4](4,7)
  \rput(2,2.5){\holel}
  \rput(2,0){\comultl}
\end{pspicture}
  \xy {\ar@{|->} (-2,0);(2,0)};
      (0,-4)*{\eqref{movelholeI}};
  \endxy
\begin{pspicture}[.4](4,7)
  \rput(2,0){\holel}
  \rput(3,4){\comultl}
\end{pspicture}
\end{math}
%
%
\item
  The moves of steps \ref{stepIIIa2} and \ref{stepIIIb2}.
\end{enumerate}
  All other configurations are excluded by step \ref{stepI}.
  Iterating these moves is guaranteed to remove all instances of the
  open window since each iteration either removes the window or
  reduces its height. The height of the open window cannot be
  zero.
\item
\label{stepIIId}
  From the sequence of moves applied thus far, it follows that the
  target of every open multiplication is in one of the following
  configurations:
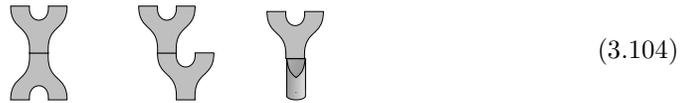
\begin{equation}
\psset{xunit=.25cm,yunit=.25cm}
\begin{pspicture}[.4](4,5.5)
  \rput(2,0){\comultl}
  \rput(2,2.5){\multl}
\end{pspicture}
\qquad
\begin{pspicture}[.4](4,5.5)
  \rput(3,0){\multl}
  \rput(2,2.5){\multl}
\end{pspicture}
\qquad
\begin{pspicture}[.4](4,5)
  \rput(2.2,0){\ltc}
  \rput(2,2){\multl}
\end{pspicture}
\end{equation}
  All other possibilities are excluded by steps~\ref{stepIIIa1},
  \ref{stepI} and \ref{stepIIIc}.
\end{enumerate}
%
%
\item
\label{stepIV}
  In this step, we apply a sequence of moves that removes all open
  comultiplications. After step~\ref{stepII}, we need to consider only
  three cases. Step~\ref{stepIII} has not changed this
  situation.  From the set of open comultiplications in the
  decomposition of $f$, choose one of minimal height.

\begin{enumerate}[a)]
\item
\label{stepIVa}
  The case
  $\psset{xunit=.2cm,yunit=.2cm}
  \begin{pspicture}[.5](5,5.1)\rput(2,0){\comultl}\rput(3,2.5){\comultl}\end{pspicture}$
  has been excluded by the assumption that the open comultiplication
  is of minimal height. Hence the only remaining configurations to
  consider are
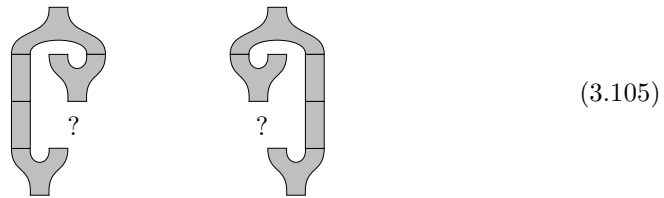
\begin{equation}
\psset{xunit=.25cm,yunit=.25cm}
\begin{pspicture}[.5](6,10)
  \rput(3,0){\multl}
  \rput(2,2.5){\identl}
  \rput(2,5){\identl}
  \rput(5,5){\multl}
  \rput(4,7.5){\widecomultl}
  \rput(4.8,3.7){$?$}
\end{pspicture}
\qquad\qquad
\begin{pspicture}[.5](6,10)
  \rput(5,0){\multl}
  \rput(6,2.5){\identl}
  \rput(6,5){\identl}
  \rput(3,5){\multl}
  \rput(4,7.5){\widecomultl}
  \rput(3.2,3.7){$?$}
\end{pspicture}
\end{equation}
  where no other open comultiplication occurs in `?' above.
%
%
\item
\label{stepIVb}
  By symmetry it suffices to consider one of the remaining
  configurations, say
\begin{math}
\psset{xunit=.2cm,yunit=.2cm}
\begin{pspicture}[.5](6,10)
  \rput(4,0){\multl}
  \rput(5,2.5){\identl}
  \rput(5,5){\identl}
  \rput(2,5){\multl}
  \rput(3,7.5){\widecomultl}
  \rput(2.2,3.7){$?$}
\end{pspicture}
\end{math}.
  We proceed by considering all possible configurations of `?'
  above. The first generator in the decomposition of `?' is determined
  by step \ref{stepIIId} and the assumption that the open
  comultiplication under consideration is of minimal height. Hence,
  only two situations are possible:
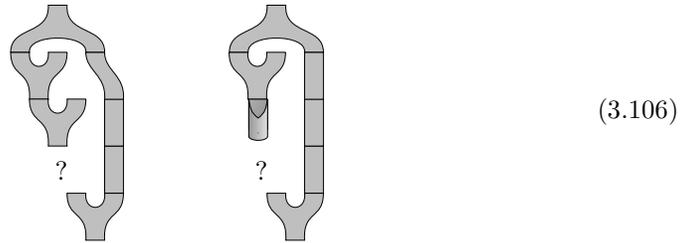
\begin{equation}
\psset{xunit=.25cm,yunit=.25cm}
\begin{pspicture}[.5](6,13)
  \rput(5,0){\multl}
  \rput(6,2.5){\identl}
  \rput(6,5){\identl}
  \rput(6,7.5){\curveleftl}
  \rput(2,7.5){\multl}
  \rput(3,5){\multl}
  \rput(3,10){\widecomultl}
  \rput(3.2,3.7){$?$}
\end{pspicture}
\qquad\qquad
\begin{pspicture}[.5](6,13)
  \rput(4,0){\multl}
  \rput(5,2.5){\identl}
  \rput(5,5){\identl}
  \rput(5,7.5){\identl}
  \rput(2,7.5){\multl}
  \rput(2.2,5.5){\ltc}
  \rput(3,10){\widecomultl}
  \rput(2.2,3.7){$?$}
\end{pspicture}
\end{equation}
%
%
\begin{enumerate}[1)]
\item
\label{stepIVb1}
   In the first case,
   iteratively apply the move
\begin{math}
\psset{xunit=.25cm,yunit=.25cm}
\begin{pspicture}[0.5](5,5)
  \rput(3,0){\multl}
  \rput(2,2.5){\multl}
\end{pspicture}
  \xy {\ar@{|->} (-2,0);(2,0)};
      (0,-4)*{\eqref{algebral}};
  \endxy
\begin{pspicture}[0.5](5,5)
  \rput(2,0){\multl}
  \rput(3,2.5){\multl}
\end{pspicture}
\end{math}
  so that the only possible configurations are
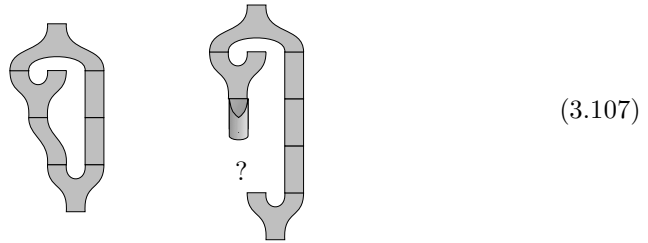
\begin{equation}
\label{eq_stepIVb1} \psset{xunit=.25cm,yunit=.25cm}
\begin{pspicture}[.5](6,10)
  \rput(5,0){\multl}
  \rput(6,2.5){\identl}
  \rput(6,5){\identl}
  \rput(3,5){\multl}
  \rput(4,7.5){\widecomultl}
  \rput(4,2.5){\curveleftl}
\end{pspicture}
\qquad\qquad
\begin{pspicture}[.5](6,13)
  \rput(4,0){\multl}
  \rput(5,2.5){\identl}
  \rput(5,5){\identl}
  \rput(5,7.5){\identl}
  \rput(2,7.5){\multl}
  \rput(2.2,5.5){\ltc}
  \rput(3,10){\widecomultl}
  \rput(2.2,3.7){$?$}
\end{pspicture}
\end{equation}
  In the following two steps we remove the open
  comultiplication from the above two situations.
%
%
\item
\label{stepIVb2}
  Consider the first case in~\eqref{eq_stepIVb1} above. The
  comultiplication is removed by the following sequence of
  moves:
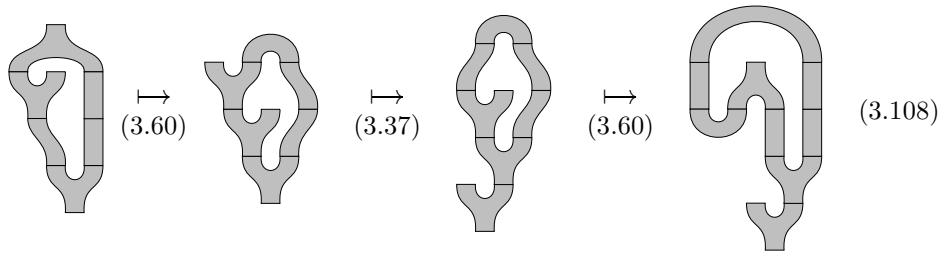
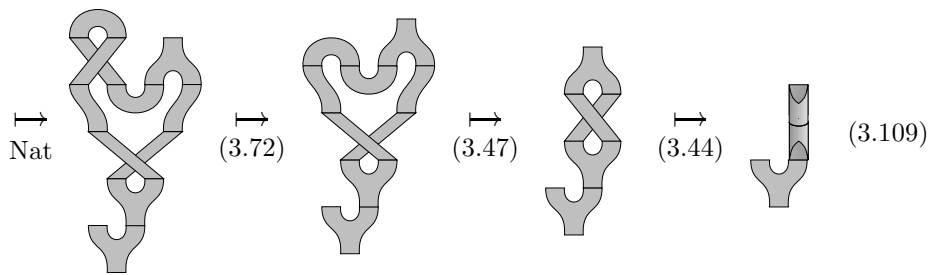
\begin{equation}
\psset{xunit=.25cm,yunit=.25cm}
\begin{pspicture}[.5](6,10)
  \rput(5,0){\multl}
  \rput(6,2.5){\identl}
  \rput(6,5){\identl}
  \rput(3,5){\multl}
  \rput(4,7.5){\widecomultl}
  \rput(4,2.5){\curveleftl}
\end{pspicture}
\quad
  \xy {\ar@{|->} (-2,2);(2,2)};
      (0,-2)*{\eqref{comultpairingIl}};
  \endxy\quad
\begin{pspicture}[.5](6,9)
  \rput(3,0){\multl}
  \rput(2,2.5){\multl}
  \rput(4,2.5){\curverightl}
  \rput(5,5){\curveleftl}
  \rput(3,7.5){\zagl}
  \rput(1,5){\multl}
\end{pspicture}
\quad
  \xy {\ar@{|->} (-2,2);(2,2)};
      (0,-2)*{\eqref{algebral}};
  \endxy\quad
\begin{pspicture}[.5](6,12)
  \rput(2,0){\multl}
  \rput(3,2.5){\multl}
  \rput(4,5){\curverightl}
  \rput(5,7.5){\curveleftl}
  \rput(2,5){\multl}
  \rput(1,7.5){\curverightl}
  \rput(3,10){\zagl}
\end{pspicture}
\quad
  \xy {\ar@{|->} (-2,2);(2,2)};
      (0,-2)*{\eqref{comultpairingIl}};
  \endxy\quad
\begin{pspicture}[.5](8,14)
  \rput(5,0){\multl}
  \rput(6,2.5){\multl}
  \rput(7,5){\identl}
  \rput(7,7.5){\identl}
  \rput(1,7.5){\identl}
  \rput(4,10){\bigzagl}
  \rput(4,7.5){\comultl}
  \rput(5,5){\identl}
  \rput(2,5.5){\zigl}
\end{pspicture}
\end{equation}
\begin{equation}
\psset{xunit=.25cm,yunit=.25cm} \quad
  \xy {\ar@{|->} (-2,2);(2,2)};
      (0,-2)*{{\rm Nat}};
  \endxy\quad
\begin{pspicture}[.5](6,14)
  \rput(2,0){\multl}
  \rput(3,2.5){\multl}
  \rput(3,5){
  \pspolygon[fillcolor=lightgray,fillstyle=solid](-1.5,0)(-.5,0)(2.5,2.5)(1.5,2.5)(-1.5,0)
  \pspolygon[fillcolor=lightgray,fillstyle=solid](1.5,0)(.5,0)(-2.5,2.5)(-1.5,2.5)(1.5,0)}
  \rput(1,10){\ucrossl}
  \rput(1,12.5){\zagl}
  \rput(1,7.5){\curveleftl}
  \rput(5,10){\comultl}
  \rput(5,7.5){\curverightl}
  \rput(3,8){\zigl}
\end{pspicture}
\quad
  \xy {\ar@{|->} (-2,2);(2,2)};
      (0,-2)*{\eqref{symmcoparingl}};
  \endxy\quad
\begin{pspicture}[.5](6,12)
  \rput(2,0){\multl}
  \rput(3,2.5){\multl}
  \rput(3,5){
  \pspolygon[fillcolor=lightgray,fillstyle=solid](-1.5,0)(-.5,0)(2.5,2.5)(1.5,2.5)(-1.5,0)
  \pspolygon[fillcolor=lightgray,fillstyle=solid](1.5,0)(.5,0)(-2.5,2.5)(-1.5,2.5)(1.5,0)}
  \rput(1,10){\zagl}
  \rput(1,7.5){\curveleftl}
  \rput(5,10){\comultl}
  \rput(5,7.5){\curverightl}
  \rput(3,8){\zigl}
\end{pspicture}
\quad
  \xy {\ar@{|->} (-2,2);(2,2)};
      (0,-2)*{\eqref{zigzagl}};
  \endxy\quad
\begin{pspicture}[.5](4.5,10)
  \rput(1.5,0){\multl}
  \rput(2.5,2.5){\multl}
  \rput(2.5,5){\crossl}
  \rput(2.5,7.5){\comultl}
\end{pspicture}
\quad
  \xy {\ar@{|->} (-2,2);(2,2)};
      (0,-2)*{\eqref{cardy}};
  \endxy\quad
\begin{pspicture}[.5](3,7)
  \rput(1.5,0){\multl}
  \rput(2.4,2.5){\ctl}
  \rput(2.7,4.5){\ltc}
  \end{pspicture}
\end{equation}
%
%
\item
\label{stepIVb3}
  Consider now the second case in~\eqref{eq_stepIVb1}. In this case
  the comultiplication is removed by the following sequence of
  moves:
  \begin{center}\makebox[0pt]{
$ \psset{xunit=.25cm,yunit=.25cm}
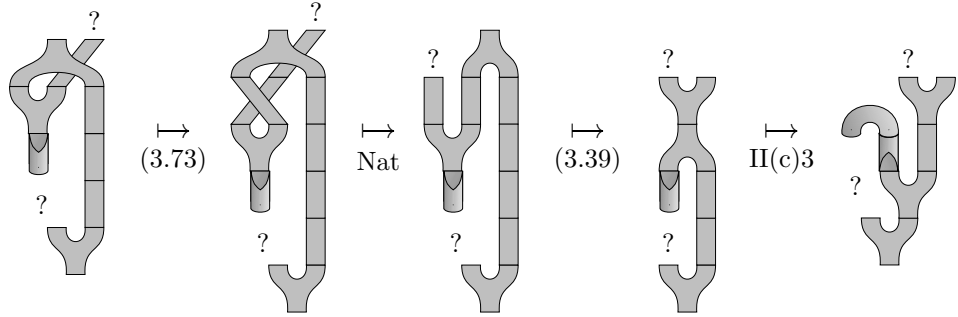
\begin{pspicture}[.5](6,13)
  \rput(4,0){\multl}
  \rput(5,2.5){\identl}
  \rput(5,5){\identl}
  \rput(5,7.5){\identl}
  \rput(2,7.5){\multl}
  \rput(4,10){
  \pspolygon[fillcolor=lightgray,fillstyle=solid](-1.5,0)(-.5,0)(1.5,2.5)(.5,2.5)(-1.5,0)}
  \rput(5,13.5){$?$}
  \rput(2.2,5.5){\ltc}
  \rput(3,10){\widecomultl}
  \rput(2.2,3.7){$?$}
\end{pspicture}
\quad
  \xy {\ar@{|->} (-2,2);(2,2)};
      (0,-2)*{\eqref{zipcross}};
  \endxy\;\;
\begin{pspicture}[.5](6,17)
  \rput(4,0){\multl}
  \rput(5,2.5){\identl}
  \rput(5,5){\identl}
  \rput(5,10){\identl}
  \rput(5,7.5){\identl}
  \rput(2,7.5){\multl}
  \rput(2,10){\crossl}
  \rput(4,12.5){
\pspolygon[fillcolor=lightgray,fillstyle=solid](-1.5,0)(-.5,0)(1.5,2.5)(.5,2.5)(-1.5,0)}
  \rput(2.2,5.5){\ltc}
  \rput(3,12.5){\widecomultl}
  \rput(2.2,3.7){$?$}
  \rput(5,16){$?$}
\end{pspicture}
\;\;\;
  \xy {\ar@{|->} (-2,2);(2,2)};
      (0,-2)*{{\rm Nat}};
  \endxy\;\;
\begin{pspicture}[.5](6,17)
  \rput(4,0){\multl}
  \rput(5,2.5){\identl}
  \rput(5,5){\identl}
  \rput(5,10){\identl}
  \rput(5,7.5){\identl}
  \rput(2,7.5){\multl}
  \rput(1,10){\identl}
  \rput(3,10){\identl}
  \rput(4,12.5){\comultl}
  \rput(2.2,5.5){\ltc}
  \rput(2.2,3.7){$?$}
  \rput(1,13.5){$?$}
\end{pspicture}
\quad
  \xy {\ar@{|->} (-2,2);(2,2)};
      (0,-2)*{\eqref{frobeniusl}};
  \endxy
\begin{pspicture}[.5](6,17)
  \rput(3.5,0){\multl}
  \rput(4.5,2.5){\identl}
  \rput(4.5,5){\identl}
  \rput(3.5,7.5){\comultl}
  \rput(3.5,10){\multl}
  \rput(2.7,5.5){\ltc}
  \rput(2.5,3.7){$?$}
  \rput(2.5,13.5){$?$}
\end{pspicture}
\;\;
  \xy {\ar@{|->} (-2,2);(2,2)};
      (0,-2)*{\ref{stepII3}};
  \endxy
\begin{pspicture}[.5](6,12)
  \rput(4,0){\multl}
  \rput(5,2.5){\multl}
  \rput(6,5){\identl}
  \rput(6,7.5){\multl}
  \rput(3.8,5){\ctl}
  \rput(3.2,7){\zagc}
  \rput(2.2,4.2){$?$}
  \rput(5,11){$?$}
\end{pspicture}
$}\end{center}
\end{enumerate}
%
%
\item
\label{stepIVc}
  Step~\ref{stepIVb} has changed the cobordism so much that the claims
  made in steps~\ref{stepII} and~\ref{stepIII} need not hold any
  longer. We therefore reapply the steps~\ref{stepII}
  and~\ref{stepIII}.
\item
\label{stepIVd}
  Then we iterate the sequence of steps~\ref{stepIVb}
  and~\ref{stepIVc} until all open comultiplications have
  disappeared. This iteration terminates because neither
  step~\ref{stepII} nor step~\ref{stepIII} (which are invoked
  in~\ref{stepIVc}) increase the number of open comultiplications, but
  step~\ref{stepIVb} always decreases this number by one.
\item
  When the last open comultiplication has disappeared in
  step~\ref{stepIVd}, step~\ref{stepIVc} ensures that the claims made
  in steps~\ref{stepII} and~\ref{stepIII} are satisfied again.
\end{enumerate}
%
%
\item
\label{stepV}
  At this stage of the proof, all open caps, open cups and open
  comultiplications have been removed from the decomposition of
  $f$. The decomposition has the following further properties.
%
%
\begin{enumerate}[a)]
\item
\label{stepVa}
  After the step~\ref{stepIIIa}, it is clear that any open
  multiplication has its source in one of the following
  configurations:
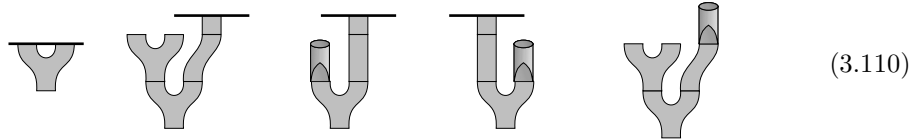
\begin{equation}
\psset{xunit=.25cm,yunit=.25cm}
\begin{pspicture}[0.5](3,4)
  \rput(2,1){\multl}
  \psline[linewidth=1pt](-0,3.5)(4,3.5)
\end{pspicture}
\qquad
\begin{pspicture}[0.5](6,6)
  \rput(3,0){\multl}
  \rput(2,2.5){\multl}
  \rput(4,2.5){\curverightl}
  \rput(5,5){\smallidentl}
  \psline[linewidth=1pt](3,6)(7,6)
\end{pspicture}
\qquad
\begin{pspicture}[0.5](6,6)
  \rput(3,0){\multl}
  \rput(1.8,2.5){\ctl}
  \rput(4,2.5){\identl}
  \rput(4,5){\smallidentl}
  \psline[linewidth=1pt](2,6)(6,6)
\end{pspicture}
\qquad
\begin{pspicture}[0.5](6,6)
  \rput(3,0){\multl}
  \rput(3.8,2.5){\ctl}
  \rput(2,2.5){\identl}
  \rput(2,5){\smallidentl}
  \psline[linewidth=1pt](0,6)(4,6)
\end{pspicture}
\qquad
\begin{pspicture}[0.5](5,7)
  \rput(3,0){\multl}
  \rput(2,2.5){\multl}
  \rput(4,2.5){\curverightl}
  \rput(4.83,5){\ctl}
\end{pspicture}
\end{equation}
%
%
\item
\label{stepVb}
  Every instance of the cozipper
  $\psset{xunit=.2cm,yunit=.2cm}\begin{pspicture}[.4](1.8,2)\rput(1,0){\ltc}\end{pspicture}$
  is in the configuration claimed in step \ref{stepIIIb}.
%
%
\item
\label{stepVc}
  All instances of
  $\psset{xunit=.2cm,yunit=.2cm}\begin{pspicture}[.4](4,4)\rput(2,0){\holel}\end{pspicture}$
  have been removed by step \ref{stepIIIc}.
%
%
\item
\label{stepVd}
  From step \ref{stepIIId} and step~\ref{stepIV}, the only possible
  configurations for the target of an open multiplication are
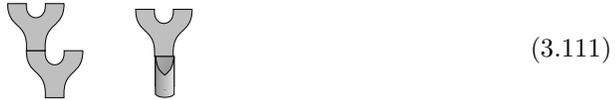
\begin{equation}
\psset{xunit=.25cm,yunit=.25cm}
\begin{pspicture}[.4](4,5.5)
  \rput(3,0){\multl}
  \rput(2,2.5){\multl}
\end{pspicture}
\qquad
\begin{pspicture}[.4](4,5)
  \rput(2.2,0){\ltc}
  \rput(2,2){\multl}
\end{pspicture}
\end{equation}
\end{enumerate}
%
%
\item
\label{stepVI}
  Now we remove every instance of the zipper
  $\psset{xunit=.2cm,yunit=.2cm}\begin{pspicture}[.4](1.8,2)\rput(1,0){\ctl}\end{pspicture}$
  in the decomposition of $f$. We consider all remaining possible
  configurations involving the zipper and provide the moves to get rid
  of it.
%
%
\begin{enumerate}[a)]
\item
\label{stepVIa}
  The following configurations:
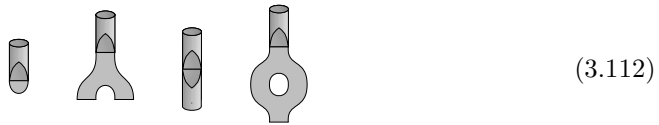
\begin{equation}
\psset{xunit=.25cm,yunit=.25cm}
\begin{pspicture}[.4](1.8,3)
  \rput(1,1){\ctl}
  \rput(1.15,1){\deathl}
\end{pspicture}
\qquad
\begin{pspicture}[.4](1.8,3)
  \rput(1,2.5){\ctl}
  \rput(1.15,0){\comultl}
\end{pspicture}
\qquad
\begin{pspicture}[.4](1.8,4)
  \rput(1,2){\ctl}
  \rput(1.3,0){\ltc}
\end{pspicture}
\qquad
\begin{pspicture}[.4](1.8,6)
  \rput(1,4){\ctl}
  \rput(1.15,0){\holel}
\end{pspicture}
\end{equation}
  are excluded by steps~\ref{stepI}, \ref{stepIV}, \ref{stepIIIb},
  and~\ref{stepIIIc}, respectively.
%
%
\item
\label{stepVIb}
  The remaining possibilities are
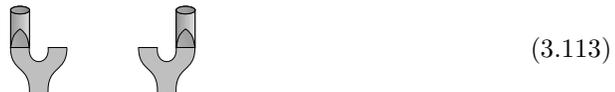
\begin{equation}
\psset{xunit=.25cm,yunit=.25cm}
\begin{pspicture}[.4](4,5)
  \rput(1,2.5){\ctl}
  \rput(2.15,0){\multl}
\end{pspicture}
\qquad
\begin{pspicture}[.4](4,5)
  \rput(3,2.5){\ctl}
  \rput(2.15,0){\multl}
\end{pspicture}
\end{equation}
  Using step~\ref{stepVd} together with possibly repeated applications
  of the following moves:
%
%
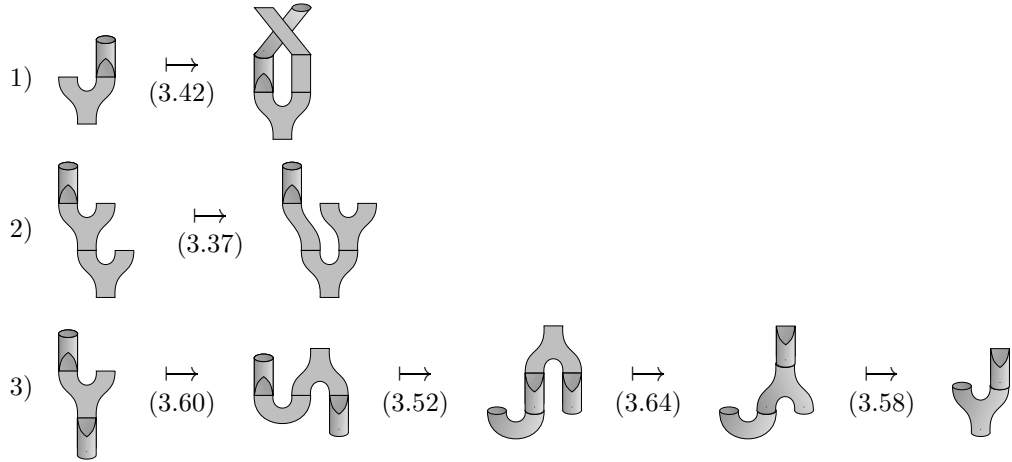
\begin{enumerate}[1)]
\item
\label{stepVIb1}
\begin{math}
\psset{xunit=.25cm,yunit=.25cm}
\begin{pspicture}[.4](4,5)
  \rput(3,2.5){\ctl}
  \rput(2.15,0){\multl}
\end{pspicture}
\quad
  \xy {\ar@{|->} (-2,2);(2,2)};
      (0,-2)*{\eqref{zipcenter}};
  \endxy\quad
\begin{pspicture}[.4](4,7)
  \rput(1,2.5){\ctl}
  \rput(2.15,0){\multl}
  \rput(3.15,2.5){\medidentl}
  \rput(2.4,4.5){\crossmixlc}
\end{pspicture}
\end{math}
%
%
\item
\label{stepVIb2}
\begin{math}
\psset{xunit=.25cm,yunit=.25cm}
\begin{pspicture}[.4](5.5,8)
  \rput(1,5){\ctl}
  \rput(2.15,2.5){\multl}
  \rput(3.15,0){\multl}
\end{pspicture}
\quad
  \xy {\ar@{|->} (-2,2);(2,2)};
      (0,-2)*{\eqref{algebral}};
  \endxy\quad
\begin{pspicture}[.4](4,8)
  \rput(1,5){\ctl}
  \rput(4.15,2.5){\multl}
  \rput(2.15,2.5){\curveleftl}
  \rput(3.15,0){\multl}
\end{pspicture}
\end{math}
%
%
\item
\label{stepVIb3}
\begin{math}
\psset{xunit=.25cm,yunit=.25cm}
\begin{pspicture}[.4](4,8)
  \rput(1,4.5){\ctl}
  \rput(2.15,2){\multl}
  \rput(2.30,0){\ltc}
\end{pspicture}
\quad
  \xy {\ar@{|->} (-2,2);(2,2)};
      (0,-2)*{\eqref{comultpairingIl}};
  \endxy\quad
\begin{pspicture}[.4](6,5)
  \rput(1,2){\ctl}
  \rput(4.15,2){\comultl}
  \rput(2.15,0){\zigl}
  \rput(5.30,0){\ltc}
\end{pspicture}
\quad
  \xy {\ar@{|->} (-2,2);(2,2)};
      (0,-2)*{\eqref{zipdual0}};
\endxy\quad
\begin{pspicture}[.4](6,7)
  \rput(2.15,0){\zigc}
  \rput(3.30,2){\ltc}
  \rput(4.15,4){\comultl}
  \rput(5.30,2){\ltc}
\end{pspicture}
\quad
  \xy {\ar@{|->} (-2,2);(2,2)};
      (0,-2)*{\eqref{counithomo}};
  \endxy\quad
\begin{pspicture}[.4](6,7)
  \rput(2.1,0){\zigc}
  \rput(4.30,2){\comultc}
  \rput(4.25,4.5){\ltc}
\end{pspicture}
\quad
  \xy {\ar@{|->} (-2,2);(2,2)};
      (0,-2)*{\eqref{comultpairingIc}};
  \endxy\quad
\begin{pspicture}[.4](4,5)
  \rput(2.30,0){\multc}
  \rput(3.25,2.5){\ltc}
\end{pspicture}
\end{math}
\end{enumerate}
  we make sure that all instances of
  $\psset{xunit=.2cm,yunit=.2cm}\begin{pspicture}[.4](1.8,2)\rput(1,0){\ctl}\end{pspicture}$
  have disappeared.
\end{enumerate}
%
%
\item
\label{stepVII}
  The resulting decomposition of $f$ is equivalent to one in
  which each \emph{closed multiplication}
  $\psset{xunit=.2cm,yunit=.2cm}\begin{pspicture}[.4](3,3)\rput(2,0){\multc}\end{pspicture}$
  has its source in one of these configurations:
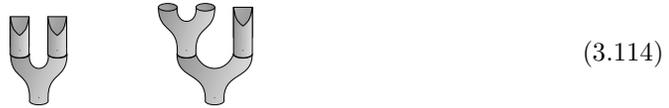
\begin{equation}
\psset{xunit=.25cm,yunit=.25cm}
\begin{aligned}
\begin{pspicture}[.4](4,5)
  \rput(2,0){\multc}
  \rput(2.9,2.5){\ltc}
  \rput(.9,2.5){\ltc}
\end{pspicture}
\end{aligned}
\qquad\qquad
\begin{aligned}
\begin{pspicture}[.4](4,5)
  \rput(1.5,0){\widemultc}
  \rput(0,2.5){\multc}
  \rput(3,2.5){\ltcnew}
\end{pspicture}
\end{aligned}
\end{equation}
%
%
\begin{enumerate}[a)]
\item
\label{stepVIIa}
  The cases
\begin{math}
\psset{xunit=.25cm,yunit=.25cm}
\begin{pspicture}[.4](4,5)
  \rput(2,0.1){\multc}
  \psline[linewidth=1pt](-0.2,2.8)(3.8,2.8)
\end{pspicture}
\end{math},
\begin{math}
\psset{xunit=.25cm,yunit=.25cm}
\begin{pspicture}[.4](4,5)
  \rput(2,.1){\multc}
  \rput(3,2.6){\medidentc}
  \rput(.9,2.6){\ltc}
  \psline[linewidth=1pt](1.8,4.8)(3.8,4.8)
\end{pspicture}
\end{math}
  and
\begin{math}
\psset{xunit=.25cm,yunit=.25cm}
\begin{pspicture}[.4](4,5)
  \rput(2,.1){\multc}
  \rput(2.9,2.6){\ltc}
  \rput(1,2.6){\medidentc}
  \psline[linewidth=1pt](-0.2,4.8)(1.8,4.8)
\end{pspicture}
\end{math}
  are excluded by the assumption that the source of $f$ is of the
  form $\vec{n}=(1,\ldots,1)$, \ie\ is a free union of intervals $I$.
%
%
\item
\label{stepVIIb}
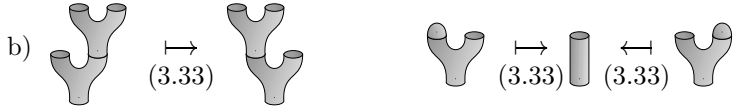
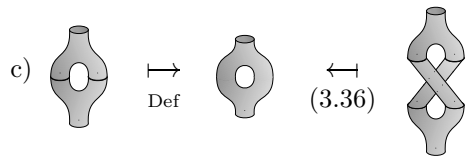
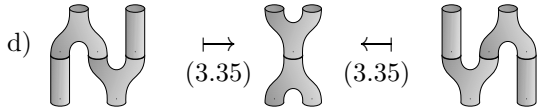
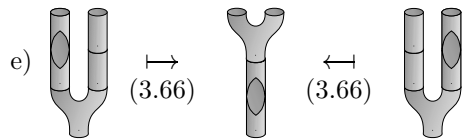
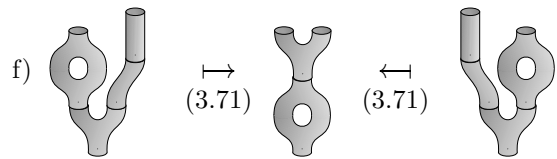
\begin{math}
\psset{xunit=.25cm,yunit=.25cm}
\begin{pspicture}[0.5](4,5)
  \rput(2,0){\multc}
  \rput(3,2.5){\multc}
\end{pspicture}
\quad
  \xy {\ar@{|->} (-2,0);(2,0)};
      (0,-4)*{\eqref{algebrac}};
  \endxy\quad
\begin{pspicture}[0.5](4,5)
  \rput(2,0){\multc}
  \rput(1,2.5){\multc}
\end{pspicture}
\qquad\qquad \psset{xunit=.25cm,yunit=.25cm}
\begin{pspicture}[0.5](4,3)
  \rput(2,0){\multc}
  \rput(1,2.5){\birthc}
\end{pspicture}
  \xy {\ar@{|->} (-2,0);(2,0)};
      (0,-4)*{\eqref{algebrac}};
  \endxy
\begin{pspicture}[0.5](2,3)
  \rput(1,0){\identc}
\end{pspicture}
  \xy {\ar@{|->} (2,0);(-2,0)};
      (0,-4)*{\eqref{algebrac}};
  \endxy
\begin{pspicture}[0.5](4,3)
  \rput(2,0){\multc}
  \rput(3,2.5){\birthc}
\end{pspicture}
\end{math}
%
%
\item
\label{stepVIIc}
\begin{math}
\psset{xunit=.25cm,yunit=.25cm}
\begin{pspicture}[0.5](4,5)
  \rput(2,0){\multc}
  \rput(2,2.5){\comultc}
\end{pspicture}
\quad
  \xy {\ar@{|->} (-2,0);(2,0)};
      (0,-4)*{{\rm \scs Def}};
  \endxy\quad
\begin{pspicture}[0.5](4,4)
  \rput(2,0){\holec}
\end{pspicture}
\quad
  \xy {\ar@{|->} (2,0);(-2,0)};
      (0,-4)*{\eqref{commutavityc}};
  \endxy\quad
\begin{pspicture}[0.5](4,8)
  \rput(2,0){\multc}
  \rput(2,2.5){\crossc}
  \rput(2,5){\comultc}
\end{pspicture}
\end{math}
%
%
\item
\label{stepVIId}
\begin{math}
\psset{xunit=.25cm,yunit=.25cm}
\begin{pspicture}[0.5](6,5.5)
  \rput(1,0){\identc}
  \rput(4,0){\multc}
  \rput(2,2.5){\comultc}
  \rput(5,2.5){\identc}
\end{pspicture}
\quad
  \xy {\ar@{|->} (-2,0);(2,0)};
      (0,-4)*{\eqref{frobeniusc}};
  \endxy\quad
\begin{pspicture}[0.5](2,5.5)
  \rput(1,0){\comultc}
  \rput(1,2.5){\multc}
\end{pspicture}
\quad
  \xy {\ar@{|->} (2,0);(-2,0)};
      (0,-4)*{\eqref{frobeniusc}};
  \endxy\quad
\begin{pspicture}[0.5](5,5.5)
  \rput(5,0){\identc}
  \rput(2,0){\multc}
  \rput(4,2.5){\comultc}
  \rput(1,2.5){\identc}
\end{pspicture}
\end{math}
%
%
\item
\label{stepVIIe}
\begin{math}
\psset{xunit=.25cm,yunit=.25cm}
\begin{pspicture}[.5](3,7)
  \rput(2,0){\multc}
  \rput(.92,2.49){\holew}
  \rput(3,2.49){\medidentc}
  \rput(3,4.49){\medidentc}
\end{pspicture}
\quad
  \xy {\ar@{|->} (-2,0);(2,0)};
      (0,-4)*{\eqref{movewholeII}};
  \endxy\quad
\begin{pspicture}[.5](3,7)
  \rput(1.92,0){\holew}
  \rput(2,4){\multc}
\end{pspicture}
\quad
  \xy {\ar@{|->} (2,0);(-2,0)};
      (0,-4)*{\eqref{movewholeII}};
  \endxy\quad
\begin{pspicture}[.5](3,7)
  \rput(2,0){\multc}
  \rput(2.92,2.49){\holew}
   \rput(1,2.49){\medidentc}
  \rput(1,4.49){\medidentc}
\end{pspicture}
\end{math}
%
%
\item
\label{stepVIIf}
\begin{math}
\psset{xunit=.25cm,yunit=.25cm}
\begin{pspicture}[.5](6,8)
  \rput(3,0){\multc}
  \rput(1.73,2.5){\holec}
  \rput(4,2.5){\curverightc}
  \rput(5,5){\identc}
\end{pspicture}
\quad
  \xy {\ar@{|->} (-2,0);(2,0)};
      (0,-4)*{\eqref{movecholeII}};
  \endxy\quad
\begin{pspicture}[.5](3,8)
  \rput(1.23,0){\holec}
  \rput(1.5,4){\multc}
\end{pspicture}
\quad
  \xy {\ar@{|->} (2,0);(-2,0)};
      (0,-4)*{\eqref{movecholeII}};
  \endxy\quad
\begin{pspicture}[.5](5,8)
  \rput(3,0){\multc}
  \rput(3.73,2.5){\holec}
  \rput(2,2.5){\curveleftc}
  \rput(1,5){\identc}
\end{pspicture}
\end{math}
\end{enumerate}
  Since each of the above moves either removes the closed
  multiplication or increases its height while not increasing the
  number of generators, iterating the above moves is guaranteed
  to terminate with the closed multiplication in one of the specified
  configurations.
%
%
\item
\label{stepVIII}
  The decomposition of $f$ is equivalent to one in which each
  closed comultiplication is in one of the following two
  configurations:
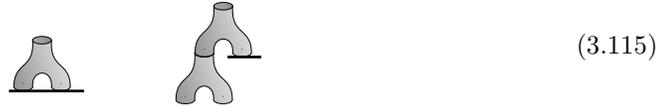
\begin{equation}
\psset{xunit=.25cm,yunit=.25cm}
\begin{pspicture}[0.5](3,4)
  \rput(2,.2){\comultc}
  \psline[linewidth=1pt](-0,0)(4,0)
\end{pspicture}
\qquad\qquad
\begin{pspicture}[0.5](4,5)
  \rput(2,0){\comultc}
  \rput(3,2.5){\comultc}
  \psline[linewidth=1pt](3,2.3)(4.8,2.3)
\end{pspicture}
\end{equation}
  We consider all possible configurations of closed
  comultiplications.
%
%
\begin{enumerate}[a)]
\item
\label{stepVIIIa}
  The cases:
\begin{equation}
\psset{xunit=.25cm,yunit=.25cm}
\begin{pspicture}[0.5](4,5)
  \rput(2,0){\multc}
  \rput(2,2.5){\comultc}
\end{pspicture}
\qquad
\begin{pspicture}[0.5](4,8)
  \rput(2,0){\multc}
  \rput(2,2.5){\crossc}
  \rput(2,5){\comultc}
\end{pspicture}
\qquad
\begin{pspicture}[0.5](6,5.5)
  \rput(1,0){\identc}
  \rput(4,0){\multc}
  \rput(2,2.5){\comultc}
  \rput(5,2.5){\identc}
\end{pspicture}
\qquad
\begin{pspicture}[0.5](5,5.5)
  \rput(5,0){\identc}
  \rput(2,0){\multc}
  \rput(4,2.5){\comultc}
  \rput(1,2.5){\identc}
\end{pspicture}
\end{equation}
  are excluded by step~\ref{stepVII}.
%
%
\item
\label{stepVIIIb} The cases:
\begin{equation}
\psset{xunit=.25cm,yunit=.25cm}
\begin{pspicture}[0.5](4,5)
  \rput(3,0){\ctl}
  \rput(2.38,1.99){\comultc}
\end{pspicture}
\qquad
\begin{pspicture}[0.5](4,5)
  \rput(1,0){\ctl}
\rput(2.38,1.99){\comultc}
\end{pspicture}
\qquad
\begin{pspicture}[0.5](4,5)
  \rput(1,0){\ctl}
  \rput(3,0){\ctl}
  \rput(2.38,1.99){\comultc}
\end{pspicture}
\end{equation}
  are excluded by step~\ref{stepVI}.
%
%
\item
\label{stepVIIIc}
  To prove the claim, we iterate the following sequences of
  moves wherever possible:
%
%
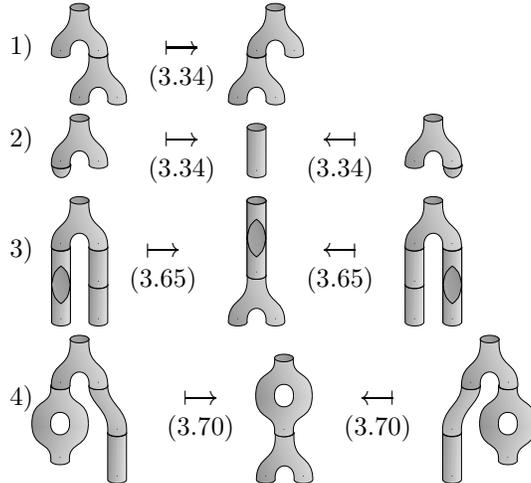
\begin{enumerate}[1)]
\item
\begin{math}
\psset{xunit=.25cm,yunit=.25cm}
\begin{pspicture}[0.5](4,5)
  \rput(3,0){\comultc}
  \rput(2,2.5){\comultc}
\end{pspicture}
\quad
  \xy {\ar@{|->} (-2,0);(2,0)};
      (0,-4)*{\eqref{coalgebrac}};
  \endxy\quad
\begin{pspicture}[0.5](4,5)
  \rput(1,0){\comultc}
  \rput(2,2.5){\comultc}
\end{pspicture}
\end{math}
%
%
\item
\begin{math}
\psset{xunit=.25cm,yunit=.25cm}
\begin{pspicture}[0.5](4,4)
  \rput(.8,0){\deathc}
  \rput(2,1){\comultc}
\end{pspicture}
\quad
  \xy {\ar@{|->} (-2,0);(2,0)};
      (0,-4)*{\eqref{coalgebrac}};
  \endxy\quad
\begin{pspicture}[0.5](2,3)
  \rput(1,0){\identc}
\end{pspicture}
\quad
  \xy {\ar@{|->} (2,0);(-2,0)};
      (0,-4)*{\eqref{coalgebrac}};
  \endxy\quad
\begin{pspicture}[0.5](4,4)
  \rput(2.8,0){\deathc}
  \rput(2,1){\comultc}
\end{pspicture}
\end{math}
%
%
\item
\label{stepVIIId}
\begin{math}
\psset{xunit=.25cm,yunit=.25cm}
\begin{pspicture}[.5](3,7)
  \rput(.93,0){\holew}
  \rput(3,0){\medidentc}
  \rput(3,2){\medidentc}
  \rput(2,3.99){\comultc}
\end{pspicture}
\quad
  \xy {\ar@{|->} (-2,0);(2,0)};
      (0,-4)*{\eqref{movewholeI}};
  \endxy\quad
\begin{pspicture}[.5](3,7)
  \rput(2,0){\comultc}
  \rput(1.92,2.49){\holew}
\end{pspicture}
\quad
  \xy {\ar@{|->} (2,0);(-2,0)};
      (0,-4)*{\eqref{movewholeI}};
  \endxy\quad
\begin{pspicture}[.5](3,7)
  \rput(2.93,0){\holew}
  \rput(1,0){\medidentc}
  \rput(1,2){\medidentc}
  \rput(2,3.99){\comultc}
\end{pspicture}
\end{math}
%
%
\item
\label{stepVIIIe}
\begin{math}
\psset{xunit=.25cm,yunit=.25cm}
\begin{pspicture}[.5](5,8)
  \rput(.73,1){\holec}
  \rput(4,0){\identc}
  \rput(4,2.5){\curveleftc}
  \rput(2,5){\comultc}
\end{pspicture}
\quad
  \xy {\ar@{|->} (-2,0);(2,0)};
      (0,-4)*{\eqref{movecholeI}};
  \endxy\quad
\begin{pspicture}[.5](3,8)
  \rput(1.5,0){\comultc}
  \rput(1.23,2.49){\holec}
\end{pspicture}
\quad
  \xy {\ar@{|->} (2,0);(-2,0)};
      (0,-4)*{\eqref{movecholeI}};
  \endxy\quad
\begin{pspicture}[.5](5,8)
  \rput(3.73,1){\holec}
  \rput(1,0){\identc}
  \rput(1,2.5){\curverightc}
  \rput(3,5){\comultc}
\end{pspicture}
\end{math}
\end{enumerate}
  This iteration is guaranteed to terminate since each move
  either decreases the height of the closed comultiplication or
  removes it.
\end{enumerate}
%
%
\item
\label{stepIX}
  In the resulting decomposition, each instance of the closed window
  $\psset{xunit=.15cm,yunit=.15cm}\begin{pspicture}[.4](1.8,4)\rput(1,0){\holew}\end{pspicture}$
  has above it one of the following:
  $\psset{xunit=.25cm,yunit=.25cm}\begin{pspicture}[.4](1.8,2)\rput(1,1.4){\birthc}\end{pspicture}$,
  $\psset{xunit=.25cm,yunit=.25cm}\begin{pspicture}[.4](1.8,4)\rput(1,0){\holew}\end{pspicture}$,
  $\psset{xunit=.25cm,yunit=.25cm}\begin{pspicture}[.4](1.8,2)\rput(1,0){\ltc}\end{pspicture}$,
  or
  $\psset{xunit=.25cm,yunit=.25cm}\begin{pspicture}[.4](3,3)\rput(1.5,0){\multc}\end{pspicture}$.
  There are only two remaining cases to consider.
%
%
\begin{enumerate}[a)]
\item
\label{stepIXa}
  The cases
\begin{math}
\psset{xunit=.25cm,yunit=.25cm}
\begin{pspicture}[.5](3,7)
  \rput(.93,0){\holew}
  \rput(3,0){\medidentc}
  \rput(3,2){\medidentc}
  \rput(2,3.99){\comultc}
\end{pspicture}
\end{math}
  and
\begin{math}
\psset{xunit=.25cm,yunit=.25cm}\begin{pspicture}[.5](4,7)
  \rput(2.93,0){\holew}
  \rput(1,0){\medidentc}
  \rput(1,2){\medidentc}
  \rput(2,3.99){\comultc}
\end{pspicture}
\end{math}
  are excluded by step~\ref{stepVIIId}.
%
%
\item
\label{stepIXb}
  The claim follows by iterating the moves
\begin{math}
\psset{xunit=.25cm,yunit=.25cm}
\begin{pspicture}[.5](3,8.5)
  \rput(2,0){\holew}
  \rput(1.83,4){\holec}
\end{pspicture}
\quad
  \xy {\ar@{|->} (-2,0);(2,0)};
      (0,-4)*{\eqref{movewholeIII}};
  \endxy\quad
\begin{pspicture}[.5](3,8.5)
  \rput(1.83,0){\holec}
  \rput(2,4){\holew}
\end{pspicture}
\end{math}.
\end{enumerate}
  At this point, the decomposition of $f$ is in the normal form
  desired. In order to see this, we need the claims made in the
  steps~\ref{stepVIII}, \ref{stepVa}, \ref{stepVI}, \ref{stepVd},
  \ref{stepVII} and the following two results.
%
%
\item
\label{stepX}
  If a \emph{closed cap}
  $\psset{xunit=.25cm,yunit=.25cm}\begin{pspicture}[.5](2,2)\rput(1,1){\birthc}\end{pspicture}$
  occurs anywhere in the resulting decomposition of $f$, then the
  source of $f$ is the object $\vec{n}=\emptyset$, and the
  $\psset{xunit=.25cm,yunit=.25cm}\begin{pspicture}[.5](2,2)\rput(1,1){\birthc}\end{pspicture}$
  has its target in one of the following configurations:
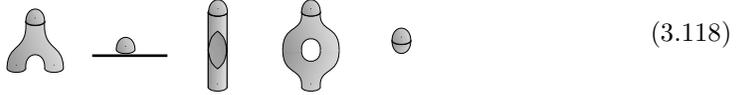
\begin{equation}
\psset{xunit=.25cm,yunit=.25cm}
\begin{pspicture}[.5](2,3)
  \rput(2,0){\comultc}
  \rput(2,2.5){\birthc}
\end{pspicture}
\qquad
\begin{pspicture}[.5](2,3)
  \psline[linewidth=1pt](-0,.7)(4,.7)
  \rput(2,1){\birthc}
\end{pspicture}
\qquad
\begin{pspicture}[.5](2,5)
  \rput(2,0){\holew}
  \rput(2.1,4){\birthc}
\end{pspicture}
\qquad
\begin{pspicture}[.5](2,5)
  \rput(2,0){\holec}
  \rput(2.24,4){\birthc}
\end{pspicture}
\qquad
\begin{pspicture}[.5](2,3)
  \rput(2,.5){\deathc}
  \rput(2.24,1.5){\birthc}
\end{pspicture}
\end{equation}
  This follows since all other possible configurations are excluded
  by steps \ref{stepVIIb} and \ref{stepVI}.
%
%
\item
\label{stepXI}
  If a \emph{closed cup}
  $\psset{xunit=.25cm,yunit=.25cm}\begin{pspicture}[.5](2,2)\rput(1,1){\deathc}\end{pspicture}$
  occurs anywhere in the resulting decomposition of $f$ then the
  target of $f$ is the object $\vec{m}=\emptyset$, and the source
  of the
  $\psset{xunit=.25cm,yunit=.25cm}\begin{pspicture}[.5](2,2)\rput(1,1){\deathc}\end{pspicture}$
  is in one of the following configurations:
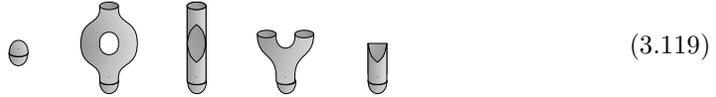
\begin{equation}
\psset{xunit=.25cm,yunit=.25cm}
\begin{pspicture}[.5](2,3)
  \rput(2,.5){\deathc}
  \rput(2.24,1.5){\birthc}
\end{pspicture}
\qquad
\begin{pspicture}[.5](2,5)
  \rput(2,0){\deathc}
  \rput(2,1){\holec}
\end{pspicture}
\qquad
\begin{pspicture}[.5](2,5)
  \rput(1.87,0){\deathc}
  \rput(2,1){\holew}
\end{pspicture}
\qquad
\begin{pspicture}[.5](2,5)
  \rput(1.82,0){\deathc}
  \rput(2,1){\multc}
\end{pspicture}
\qquad
\begin{pspicture}[.5](2,5)
  \rput(1.84,0){\deathc}
  \rput(2,1){\ltc}
\end{pspicture}
\end{equation}
  The remaining cases are excluded by step \ref{stepVIIIc}.
\end{enumerate}
  This concludes the proof.
\end{proof}

The main result for arbitrary connected open-closed cobordisms then
follows.

\begin{cor}
Let $[f]\in\twocob[\vec{n},\vec{m}]$ be connected. Then
$[f]=[\mathrm{NF}(f)]$.
\end{cor}

\begin{proof}
Using Definition~\ref{def_normal}, then applying
Theorem~\ref{thm_normal} to $\Lambda(f)$, and then
applying~\eqref{eq_lambdainv}, we find,
\begin{equation}
  [\mathrm{NF}(f)]
  = [\Lambda^{-1}([\mathrm{NF}_{\mathrm{O\rightarrow C}}(\Lambda ([f]))])]
  = [\Lambda^{-1}([\Lambda([f])])]
  = [f].
\end{equation}
\end{proof}

Since the normal form is already characterized by the invariants of
Definition~\ref{def_invariants}, we also obtain the following result.

\begin{cor}
Let $[f],[f^\prime]\in\twocob[\vec{n},\vec{m}]$ be connected such that
their genus, window number, and open boundary permutation agree, then
$[f]=[f^\prime]$.
\end{cor}

%
\section{Open-closed TQFTs}
%
\label{sect_tqft}

In this section, we define the notion of open-closed TQFTs. We show
that the categories $\twocob$ and $\Thfrob$ are equivalent as
symmetric monoidal categories which implies that the category of
open-closed TQFTs is equivalent to the category of knowledgeable
Frobenius algebras.

\begin{defn}
Let $\cal{C}$ be a symmetric monoidal category. An
\emph{open-closed Topological Quantum Field Theory (TQFT)} in
$\cal{C}$ is a symmetric monoidal functor $\twocob\to\cal{C}$. A
\emph{homomorphism} of open-closed TQFTs is a monoidal natural
transformation of such functors. By
$\cat{OC-TQFT}(\cal{C}):=\cat{Symm-Mon}(\twocob,\cal{C})$, we
denote the category of open-closed TQFTs.
\end{defn}

\begin{thm}
\label{thm_cobthfrob}
The category $\twocob$ is equivalent as a symmetric monoidal category
to the category $\Thfrob$.
\end{thm}

This theorem states the precise correspondence between topology
(Section~\ref{sect_cob}) and algebra
(Section~\ref{sect_knowfrob}). The second main result of the present
article follows from this theorem and from
Proposition~\ref{PROPthfrob}.

\begin{cor}
\label{MainThm}
Let $(\cal{C},\otimes,\1,\alpha,\lambda,\rho,\tau)$ be a symmetric
monoidal category. The category $\cat{K-Frob}(\cal{C})$ of
knowledgeable Frobenius algebras in $\cal{C}$ is equivalent as a
symmetric monoidal category to the category $\cat{OC-TQFT}(\cal{C})$.
\end{cor}

These results also guarantee that one can use the generators of
Section~\ref{sect_generators} and the relations of
Section~\ref{sect_relations} in order to perform computations in
knowledgeable Frobenius algebras. Recall that $\twocob$ is a strict
monoidal category whereas $\Thfrob$ is weak. When one translates from
diagrams to algebra, one chooses parentheses for all tensor products
and then inserts the structure isomorphisms $\alpha$, $\lambda$,
$\rho$ as appropriate. The coherence theorem of MacLane guarantees
that all ways of inserting these isomorphisms yield the same
morphisms, and so the morphisms on the algebraic side are well defined
by their diagrams.

In particular, we could have presented the second half of
Section~\ref{sect_cob}, starting with
Subsection~\ref{secConsequences}, entirely in the algebraic rather
than in the topological language.

\begin{proof}[Proof of Theorem~\ref{thm_cobthfrob}]
Define a mapping $\Xi$ from the objects of $\twocob$ to the objects of
the category $\Thfrob$ by mapping the generators as follows:
\begin{eqnarray}
\label{assign1}
  \Xi \maps \emptyset &\mapsto& \1\\
  \Xi \maps \begin{pspicture}[.3](.8,.4)
              \pscircle[linewidth=.5pt](.4,.2){.2}
            \end{pspicture} &\mapsto& C\\
  \Xi \maps \begin{pspicture}[.3](.8,.4)
              \psline[linewidth=.7pt](.1,.2)(.7,.2)
            \end{pspicture} &\mapsto& A
\end{eqnarray}
and extending to the general object $\vec{n}\in\twocob$ by mapping
$\vec{n}$ to the tensor product in $\Thfrob$ of copies of $A$ and $C$
with all parenthesis to the left. More precisely, if $\vec{n} =
(n_1,n_2,n_3,\cdots,n_k)$ with each $n_i\in\{0,1\}$, then
$\Xi(\vec{n})=\left(\left(\left(\Xi(n_1)\ten\Xi(n_2)\right)\ten\Xi(n_3)\right)\cdots\Xi(n_k)\right)$
with each $\Xi(0):=C$ and $\Xi(1):=A$. On the generating morphisms in
$\twocob$, $\Xi$ is defined as follows:
\begin{eqnarray}
\psset{xunit=.2cm,yunit=.2cm}
\begin{pspicture}[.5](2,2)
  \rput(1.2,0){\medidentc}
\end{pspicture}
  &\mapsto&
  1_C\maps C\to C\\
\psset{xunit=.2cm,yunit=.2cm}
\begin{pspicture}[.5](2,2)
  \rput(1,0){\medidentl}
\end{pspicture}
  &\mapsto&
  1_A\maps A\to A\\
\psset{xunit=.25cm,yunit=.25cm}
\begin{pspicture}[.5](2,2.5)
  \rput(1,0){\crossc}
\end{pspicture}
  &\mapsto&
  \tau_{C,C}\maps C\ten C\to C\ten C\\
\psset{xunit=.25cm,yunit=.25cm}
\begin{pspicture}[.5](2,2.5)
  \rput(1,0){\crossl}
\end{pspicture}
  &\mapsto&
  \tau_{A,A}\maps A\ten A\to A\ten A\\
\psset{xunit=.25cm,yunit=.25cm}
\begin{pspicture}[.5](2,2.5)
  \rput(1,0){\crossmixlc}
\end{pspicture}
  &\mapsto&
  \tau_{A,C}\maps A\ten C\to C\ten A\\
\psset{xunit=.25cm,yunit=.25cm}
\begin{pspicture}[.5](2,2.5)
  \rput(1,0){\crossmixcl}
\end{pspicture}
  &\mapsto&
  \tau_{C,A}\maps C\ten A\to A\ten C\\
\psset{xunit=.2cm,yunit=.2cm}
\begin{pspicture}[.5](2.2,2.5)
  \rput(1,0){\multl}
\end{pspicture}
  &\mapsto&
  \mu_A\maps A\ten A\to A\\
\psset{xunit=.2cm,yunit=.2cm}
\begin{pspicture}[.5](2,1.2)
  \rput(1,.9){\birthl}
\end{pspicture}
  &\mapsto&
  \eta_A\maps\1\to A\\
\psset{xunit=.2cm,yunit=.2cm}
\begin{pspicture}[.5](2.2,2.5)
  \rput(1,0){\comultl}
\end{pspicture}
  &\mapsto&
  \Delta_A\maps A\to A\ten A\\
\psset{xunit=.2cm,yunit=.2cm}
\begin{pspicture}[.5](2,1.2)
  \rput(1,.9){\deathl}
\end{pspicture}
  &\mapsto&
  \epsilon_A\maps A\to\1\\
\psset{xunit=.2cm,yunit=.2cm}
\begin{pspicture}[.5](2,2.5)
  \rput(1,0){\multc}
\end{pspicture}
  &\mapsto&
  \mu_C\maps C\ten C\to C\\
\psset{xunit=.2cm,yunit=.2cm}
\begin{pspicture}[.5](2,1.2)
  \rput(1.3,.5){\birthc}
\end{pspicture}
  &\mapsto&
  \eta_C\maps\1\to C\\
\psset{xunit=.2cm,yunit=.2cm}
\begin{pspicture}[.5](2,2.5)
  \rput(1,0){\comultc}
\end{pspicture}
  &\mapsto&
  \Delta_C\maps C\to C\ten C\\
\psset{xunit=.2cm,yunit=.2cm}
\begin{pspicture}[.5](2,1.2)
  \rput(1,.3){\deathc}
\end{pspicture}
  &\mapsto&
  \epsilon_C\maps C\to\1\\
\psset{xunit=.2cm,yunit=.2cm}
\begin{pspicture}[.5](2,2)
  \rput(.8,0){\ctl}
\end{pspicture}
  &\mapsto&
  \imath\maps C\to A\\
\psset{xunit=.2cm,yunit=.2cm}
\label{assign2}
\begin{pspicture}[.5](2,2)
  \rput(1.2,0){\ltc}
\end{pspicture}
  &\mapsto&
  \imath^{\ast}\maps A\to C
\end{eqnarray}
Without loss of generality we can assume that every general morphism
$f$ in $\twocob$ is decomposed into elementary generators in such
a way that each critical point in the decomposition of $f$ has a
unique critical value. We can then extend $\Xi$ to a map on all the
morphisms of $\twocob$ inductively using the following assignments:
\begin{eqnarray}
\psset{xunit=.2cm,yunit=.2cm}
\begin{pspicture}[.5](4.2,2.5)
  \rput(1,0){\multl}
  \rput(3,0){\curverightl}
\end{pspicture}
  &\mapsto&
  \mu_A\ten 1_A\maps(A\ten A)\ten A\to A\ten A\\
\psset{xunit=.2cm,yunit=.2cm}
\begin{pspicture}[.5](4.2,2.5)
  \rput(3,0){\multl}
  \rput(1,0){\curveleftl}
\end{pspicture}
  &\mapsto&
  1_A\ten\mu_A\circ\alpha_{A,A,A}\maps(A\ten A)\ten A\to A\ten A\\
\psset{xunit=.2cm,yunit=.2cm}
\begin{pspicture}[.5](4.2,2.5)
  \rput(1,0){\birthl}
  \rput(3,0){\identl}
\end{pspicture}
  &\mapsto&
  \eta_A\ten 1_A\circ\lambda_A^{-1}\maps A \to A\ten A\\
\psset{xunit=.2cm,yunit=.2cm}
\begin{pspicture}[.5](4.2,2.5)
  \rput(3,0){\birthl}
  \rput(1,0){\identl}
\end{pspicture}
  &\mapsto&
  1_A\ten \eta_A\circ\rho_A^{-1}\maps A \to A\ten A\\
\psset{xunit=.2cm,yunit=.2cm}
\begin{pspicture}[.5](2.2,2.5)
  \rput(-2,0){\comultl}
  \rput(1,0){\curveleftl}
\end{pspicture}
  &\mapsto&
  \Delta_A\ten 1_A\maps A\ten A\to (A\ten A)\ten A\\
\psset{xunit=.2cm,yunit=.2cm}
\begin{pspicture}[.5](2.2,2.5)
  \rput(1,0){\comultl}
  \rput(-2,0){\curverightl}
\end{pspicture}
  &\mapsto&
 \alpha_{A,A,A}^{-1}\circ1_A\ten\Delta_A\maps A\ten A\to(A\ten A)\ten A\\
\psset{xunit=.2cm,yunit=.2cm}
\begin{pspicture}[.5](4.2,2.5)
  \rput(1,2.5){\deathl}
  \rput(3,0){\identl}
\end{pspicture}
  &\mapsto&
  \lambda_A\circ\varepsilon_A \ten 1_A\maps A\ten A \to A\\
\psset{xunit=.2cm,yunit=.2cm}
\begin{pspicture}[.5](4.2,2.5)
  \rput(3,2.5){\deathl}
  \rput(1,0){\identl}
\end{pspicture}
  &\mapsto&
  \rho_A\circ 1_A \ten \varepsilon_A\maps A\ten A \to A \\
\psset{xunit=.2cm,yunit=.2cm}
\begin{pspicture}[.5](4,2.5)
  \rput(1,0){\multc}
  \rput(3,0){\curverightc}
\end{pspicture}
  &\mapsto&
  \mu_C\ten 1_C\maps(C\ten C)\ten C\to C\ten C\\
\psset{xunit=.2cm,yunit=.2cm}
\begin{pspicture}[.5](4,2.5)
  \rput(3,0){\multc}
  \rput(1,0){\curveleftc}
\end{pspicture}
  &\mapsto&
  1_C\ten\mu_C\circ\alpha_{C,C,C}\maps(C\ten C)\ten C\to C\ten C\\
\psset{xunit=.2cm,yunit=.2cm}
\begin{pspicture}[.5](4.2,2.5)
  \rput(1,0){\birthc}
  \rput(3,0){\identc}
\end{pspicture}
  &\mapsto&
  \eta_C\ten 1_C\circ\lambda_C^{-1}\maps C \to C\ten C\\
\psset{xunit=.2cm,yunit=.2cm}
\begin{pspicture}[.5](4.2,2.5)
  \rput(3,0){\birthc}
  \rput(1,0){\identc}
\end{pspicture}
  &\mapsto&
  1_C\ten \eta_C\circ\rho_C^{-1}\maps C\to C\ten C\\
\psset{xunit=.2cm,yunit=.2cm}
\begin{pspicture}[.5](2,2.5)
  \rput(-2,0){\comultc}
  \rput(1,0){\curveleftc}
\end{pspicture}
  &\mapsto&
  \Delta_C\ten 1_C\maps C\ten C\to(C\ten C)\ten C\\
\psset{xunit=.2cm,yunit=.2cm}
\begin{pspicture}[.5](2,2.5)
  \rput(1,0){\comultc}
  \rput(-2,0){\curverightc}
\end{pspicture}
  &\mapsto&
  \alpha_{C,C,C}^{-1}\circ1_C\ten\Delta_C\maps C\ten C\to (C\ten C)\ten
  C\\
\psset{xunit=.2cm,yunit=.2cm}
\begin{pspicture}[.5](4.2,2.5)
  \rput(1,1.5){\deathc}
  \rput(3,0){\identc}
\end{pspicture}
  &\mapsto&
  \lambda_C\circ\varepsilon_C \ten 1_C\maps C\ten C \to C\\
\psset{xunit=.2cm,yunit=.2cm}
\begin{pspicture}[.5](4.2,2.5)
  \rput(3,1.5){\deathc}
  \rput(1,0){\identc}
\end{pspicture}
  &\mapsto&
  \rho_C\circ 1_C \ten \varepsilon_C\maps C\ten C \to C
\end{eqnarray}
This assignment is well defined and extends to all the general
morphisms in $\twocob$ by the coherence theorem for symmetric monoidal
categories, which ensures that there is a unique morphism from one
object to another composed of associativity constraints and unit
constraints. The relations in Proposition~\ref{PROPrelations} and the
proof that these are all the required relations in $\twocob$ imply
that the image of $\Xi$ is in fact a knowledgeable Frobenius algebra
in $\Thfrob$.  Hence, $\Xi$ defines a functor $\twocob \to
\Thfrob$.

Define a natural isomorphism
$\Xi_2\maps\Xi(\vec{n})\ten\Xi(\vec{m})\to\Xi(\vec{n}\coprod\vec{m})$
for $X,Y\in\twocob$ as follows: Let $\vec{n}=(n_1,n_2,n_3,\ldots,n_k)$
and $\vec{m}=(m_1,m_2,m_3,\cdots,m_\ell)$ so that
\begin{eqnarray}
  \Xi(\vec{n}) &=&
    \left(\left(\left(\Xi(n_1)\ten\Xi(n_2)\right)\ten\Xi(n_3)\right)\cdots\Xi(n_k)\right),\\
  \Xi(\vec{m}) &=&
    \left(\left(\left(\Xi(m_1)\ten\Xi(m_2)\right)\ten\Xi(m_3)\right)\cdots\Xi(m_\ell)\right),\\
  \Xi(\vec{n}\coprod\vec{m}) &=&
    \left(\left(\left(\left(\left(\Xi(n_1)\ten\Xi(n_2)\right)
      \ten\Xi(n_3)\right)\cdots\Xi(n_k)\right)\ten\Xi(m_1)\right)
      \cdots\ten\Xi(m_\ell)\right).
\end{eqnarray}
Hence the map
$\Xi_2\maps\Xi(\vec{n})\ten\Xi(\vec{m})\to\Xi(\vec{n}\coprod\vec{m})$
is composed entirely of composites of the natural isomorphism
$\alpha$.  By the coherence theorem for monoidal categories, any
choice of composites from the source to the target is unique.  One
can easily verify that if $\Xi_0:=1_{\1}$, then collection
$(\Xi,\Xi_2,\Xi_0)$ defines a monoidal natural transformation.
Furthermore, our choice for the assignment by $\Xi$ of the
open-closed cobordisms generating $\twocob$'s symmetry ensures that
$(\Xi,\Xi_2,\Xi_0)$ is a symmetric monoidal functor.

Using the assignments from equations \eqref{assign1}-\eqref{assign2} we
see that the generating open-closed cobordisms in $\twocob$ define a
knowledgeable Frobenius algebra structure on the interval and circle.
Hence, by the remarks preceding Proposition~\ref{PROPthfrob} we get a
strict symmetric monoidal functor $\bar{\Xi} \maps \Thfrob \to
\twocob$. In this case, if $X$ is related to $Y$ in $\Thfrob$ by a
sequence of associators and unit constraints then $X$ and $Y$ are
mapped to the same object in $\twocob$.  We now show that $\Xi$ and
$\bar{\Xi}$ define an equivalence of categories.

Let $\vec{n}$ be a general object in $\twocob$. From the discussion
above we have that $\bar{\Xi}\Xi(\vec{n})=\vec{n}$, so that
$\bar{\Xi}\Xi(\vec{n}) = 1_{\twocob}$.  If $X$ is an object of
$\Thfrob$ then $X$ is a parenthesized word consisting of the symbols
$\1,A,C,\ten$.  Let $\bar{\Xi}(X) = (n_1,n_2,\cdots,n_n)$ where the
ordered sequence $(n_1,n_2,\cdots,n_n)$ corresponds to the ordered
sequence of $A$'s and $C$'s in $X$.  Hence, $\Xi\bar{\Xi}(X)$ is the
word obtained from $X$ by removing all the symbols $\1$ and putting
all parenthesis to the left. Thus, $\Xi\bar{\Xi}(X)$ is isomorphic
to $X$ by a sequence of associators and unit constraints. We have
therefore established the desired monoidal equivalence of symmetric
monoidal categories.
\end{proof}

The following special cases are covered by Corollary~\ref{MainThm}.

\begin{defn}
Let $\cat{2Cob}^{\rm open}$, $\cat{2Cob}^{\rm closed} =
\cat{2Cob}$, and $\cat{2Cob}^{\rm planar}$ denote the
subcategories of $\twocob$ consisting only of purely open cobordisms,
purely closed cobordisms, and purely open cobordisms that can be
embedded into the plane. An open (respectively closed, planar open)
TQFT is a functor from $\cat{2Cob}^{\rm open}$ (respectively
$\cat{2Cob}^{\rm closed}$, $\cat{2Cob}^{\rm planar}$) into a symmetric
monoidal category $\cal{C}$ ($\cal{C}$ need not be symmetric in the
planar open context).
\end{defn}

\begin{cor}
Let $\cal{C}$ be a symmetric monoidal category. The category of open
TQFTs in $\cal{C}$ is equivalent as a symmetric monoidal category to
the category of symmetric Frobenius algebras in $\cal{C}$.
\end{cor}

The following well-known result on $2$-dimensional closed
TQFTs~\cite{Abrams1,Kock} follows from Corollary~\ref{MainThm}, as
does the $2$-dimensional case of~\cite{Lau2}.

\begin{cor}
Let $\cal{C}$ be a symmetric monoidal category. The category of closed
TQFTs in $\cal{C}$ is equivalent as a symmetric monoidal category to
the category of commutative Frobenius algebras in $\cal{C}$.
\end{cor}

\begin{cor}
Let $\cal{C}$ be a monoidal category. The category of planar open
topological quantum field theories in $\cal{C}$ is equivalent to the
category of Frobenius algebras in $\cal{C}$.
\end{cor}

%
\section{Boundary labels}
%
\label{sect_labels}

In this section, we generalize the results on knowledgeable Frobenius
algebras and on open-closed cobordisms to free boundaries labeled with
elements of some set $S$. The proofs of these results are very similar
to the unlabeled case, and so we state only the results.

\begin{defn}
Let $S$ be a set. An $S$-\emph{coloured knowledgeable Frobenius
algebra}\index{knowledgeable Frobenius algebra!$S$-coloured}
\begin{equation}
  ({\{A_{ab}\}}_{a,b\in S},{\{\mu_{abc}\}}_{a,b,c\in S},
  {\{\eta_a\}}_{a\in S},{\{\Delta_{abc}\}}_{a,b,c\in S},
  {\{\epsilon_a\}}_{a\in S},C,{\{\imath_a\}}_{a\in S},
  {\{\imath^\ast_a\}}_{a\in S})
\end{equation}
in some symmetric monoidal category
$(\cal{C},\otimes,\1,\alpha,\lambda,\rho,\tau)$ consists of,
\begin{itemize}
\item
  a commutative Frobenius algebra object
  $(C,\mu,\eta,\Delta,\epsilon)$ in $\cal{C}$,
\item
  a family of objects $A_{ab}$ of $\cal{C}$, $a,b\in S$,
\item
  families of morphisms $\mu_{abc}\colon A_{ab}\otimes A_{bc}\to
  A_{ac}$, $\eta_a\colon\1\to A_{aa}$, $\Delta_{abc}\colon A_{ac}\to
  A_{ab}\otimes A_{bc}$, $\epsilon_a\colon A_{aa}\to\1$,
  $\imath_a\colon C\to A_{aa}$ and $\imath_a^\ast\colon A_{aa}\to C$
  of $\cal{C}$ for all $a,b,c\in S$ such that the following conditions
  are satisfied for all $a,b,c,d\in S$:
\begin{eqnarray}
  \mu_{abd}\circ (\id_{A_{ab}}\otimes\mu_{bcd})\circ\alpha_{A_{ab},A_{bc},A_{cd}}
    &=& \mu_{acd}\circ (\mu_{abc}\otimes\id_{A_{cd}}),\label{generalizedASS}\\
  \mu_{aab}\circ (\eta_{aa}\otimes\id_{A_{ab}})
    &=& \lambda_{A_{ab}},\label{generalizedUNITi}\\
  \mu_{abb}\circ (\id_{A_{ab}}\otimes\eta_{bb})
    &=& \rho_{A_{ab}},\label{generalizedUNITii}\\
  \alpha_{A_{ab},A_{bc},A_{cd}}\circ(\Delta_{abc}\otimes\id_{A_{cd}})\circ\Delta_{acd}
    &=& (\id_{A_{ab}}\otimes\Delta_{bcd})\circ\Delta_{abd},\\
  (\epsilon_{aa}\otimes\id_{A_{ab}})\circ\Delta_{aab}
    &=& \lambda_{A_{ab}}^{-1},\\
  (\id_{A_{ab}}\otimes\epsilon_{bb})\circ\Delta_{abb}
    &=& \rho_{A_{ab}}^{-1},\\
  \Delta_{abd}\circ\mu_{acd}
    &=& (\id_{A_{ab}}\otimes\mu_{bcd})\circ\alpha_{A_{ab},A_{bc},A_{cd}}\circ(\Delta_{abc}\otimes\id_{A_{cd}})\nn\\
    &=& (\mu_{abc}\otimes\id_{A_{cd}}) \circ\alpha_{A_{ab},A_{bc},A_{cd}}^{-1}\nn\\
    &&\quad\circ(\id_{A_{ab}}\otimes\Delta_{bcd}),\\
  \epsilon_{aa}\circ\mu_{aba}
    &=& \epsilon_{bb}\circ\mu_{bab}\circ\tau_{A_{ab},A_{ba}},\\
  \mu_{aaa}\circ(\imath_a\otimes\imath_a)
    &=& \imath_a\circ\mu,\\
  \eta_{aa}
    &=& \imath_a\circ\eta,\\
  \mu_{aab}\circ(\imath_a\otimes\id_{A_{ab}})
    &=& \mu_{abb}\circ\tau_{A_{bb},A_{ab}}\circ(\imath_b\otimes\id_{A_{ab}}),\\
  \epsilon\circ\mu\circ(\id_C\otimes\imath^\ast_a)
    &=& \epsilon_{aa}\circ\mu_{aaa}\circ(\imath_a\otimes\id_{A_{aa}}),\\
  \imath_a\circ\imath_b^\ast
    &=& \mu_{aba}\circ\tau_{A_{ba},A_{ab}}\circ\Delta_{bab}.
\end{eqnarray}
\end{itemize}
\end{defn}

It is easy to see that the notion of an $S$-coloured knowledgeable
Frobenius algebra precisely models the topological relations of
Proposition~\ref{PROPrelations} for all possible ways of labeling the
free boundaries with elements of the set $S$. The following
consequences of this definition are not difficult to see from the
diagrams of Proposition~\ref{PROPrelations}.

\begin{cor}
Let $(\{A_{ab}\},\{\mu_{abc}\},\{\eta_a\},\{\Delta_{abc}\},
\{\epsilon_a\},C,\{\imath_a\},\{\imath^\ast_a\})$ be an $S$-coloured
knowledgeable Frobenius algebra in some symmetric monoidal category
$\cal{C}$.
\begin{enumerate}
\item
  Each $A_{ab}$, $a,b\in S$, is a rigid object of $\cal{C}$ whose
  left- and right-dual is given by $A_{ba}$.
\item
  Each $A_{aa}$, $a\in S$, forms a symmetric Frobenius algebra object
  in $\cal{C}$.
\item
  Each $\imath_a\colon C\to A_{aa}$, $a\in S$, forms a homomorphism of
  algebras in $\cal{C}$.
\item
  Each $\imath^\ast_a\colon A_{aa}\to C$, $a\in S$, forms a
  homomorphism of coalgebras in $\cal{C}$.
\item
  Each $A_{ab}$ forms an $A_{aa}$-left-$A_{bb}$-right-bimodule in
  $\cal{C}$.
\item
  Each $A_{ab}$ forms an $A_{aa}$-left-$A_{bb}$-right-bicomodule in
  $\cal{C}$.
\end{enumerate}
\end{cor}

\begin{defn}
A \emph{homomorphism}
\begin{eqnarray}
  f&\colon&
(\{A_{ab}\},\{\mu_{abc}\},\{\eta_a\},\{\Delta_{abc}\},\{\epsilon_a\},C,
\{\imath_a\},\{\imath^\ast_a\})\nn\\
   &\to&(\{A^\prime_{ab}\},
\{\mu^\prime_{abc}\},\{\eta^\prime_a\},\{\Delta^\prime_{abc}\},
\{\epsilon^\prime_a\},C^\prime,\{\imath^\prime_a\},\{{\imath^\prime}^\ast_a\})
\end{eqnarray}
of $S$-coloured knowledgeable Frobenius algebras is a pair
$f=({\{f_{ab}\}}_{a,b\in S},f_C)$ consisting of a homomorphism of
Frobenius algebras $f_C\colon C\to C^\prime$ and a family of morphisms
$f_{ab}\colon A_{ab}\to A^\prime_{ab}$, $a,b\in S$ that satisfy the
following conditions for all $a,b,c\in S$:
\begin{eqnarray}
  \mu^\prime_{abc}\circ(f_{ab}\otimes f_{bc})
    &=& f_{ac}\circ\mu_{abc},\\
  \eta^\prime_a
    &=& f_{aa}\circ\eta_a,\\
  \Delta^\prime_{abc}\circ f_{ac}
    &=& (f_{ab}\otimes f_{bc})\circ\Delta_{abc},\\
  \epsilon^\prime_a\circ f_{aa}
    &=& \epsilon_a,\\
  \imath^\prime_a\circ f_C
    &=& f_{aa}\circ\imath_a,\\
  {\imath^\prime_a}^\ast\circ f_{aa}
    &=& f_C\circ\imath^\ast_a.
\end{eqnarray}
\end{defn}

\begin{defn}
By $\cat{K-Frob}^{(S)}(\cal{C})$ we denote the category of
$S$-coloured knowledgeable Frobenius algebras in some symmetric
monoidal category $\cal{C}$ and their homomorphisms.
\end{defn}

\begin{defn}
The category of \emph{open-closed TQFTs} in some symmetric monoidal
category $\cal{C}$ \emph{with free boundary labels} in some set $S$ is
the category
\begin{equation}
  \cat{OC-TQFT}^{(S)}(\cal{C}):=\cat{Symm-Mon}(\twocob(S),\cal{C}).
\end{equation}
\end{defn}

In the $S$-coloured case, the correspondence between the algebraic and
the topological category of Corollary~\ref{MainThm} generalizes to the
following result.

\begin{thm}
Let $S$ be some set and $\cal{C}$ be a symmetric monoidal
category. The categories $\cat{K-Frob}^{(S)}(\cal{C})$ and
$\cat{OC-TQFT}^{(S)}(\cal{C})$ are equivalent as symmetric monoidal
categories.
\end{thm}

One can verify~\cite{LP} that the groupoid algebra of a finite
groupoid gives rise to an $S$-coloured knowledgeable Frobenius
algebra for which $S$ is the set of objects of the groupoid.

%
\section{Conclusion}
%
\label{sect_conclusion}

In this paper, we have extended the results of classical cobordism
theory to the context of $2$-dimensional open-closed cobordisms. Using
manifolds with faces with a particular global structure, rather than
the full generality of manifolds with corners, we have defined an
appropriate category of open-closed cobordisms. Using a generalization
of Morse theory to manifolds with corners, we have found a
characterization of this category in terms of generators and
relations. In order to prove the sufficiency of the relations, we have
explicitly constructed the diffeomorphism between an arbitrary
cobordism and a normal form which is characterized by topological
invariants.

All of the technology outlined above is defined for manifolds with
faces of arbitrary dimension. Thus, our work suggests a natural
framework for studying extended topological quantum field theories in
dimensions three and four.  Using 3-manifolds or 4-manifolds with
faces, one can imagine defining a category (most likely
higher-category) of extended three or four dimensional cobordisms.  In
both cases, gluing will produce well defined composition operations
using the existing technology for manifolds with faces.  One could
then extract a list of generating cobordisms, again using a suitable
generalization of Morse theory.

The main difficulty in obtaining a complete generators and relations
description of these higher-dimensional extended cobordism categories
is the lack of general theory producing the relations. Specifically,
the handlebody theory for manifolds with boundaries and corners is not
as advanced as the standard Morse theory for closed manifolds. For the
$2$-dimensional case, we were able to use relations previously
proposed in the literature and to show the sufficiency of these
relations by finding the appropriate normal form for $2$-dimensional
open-closed cobordisms. Our induction proof shows that the proposed
relations are in fact necessary and sufficient to reduce an arbitrary
open-closed cobordism to the normal form. To extend these results to
higher-dimensions, it is expected that a more sophisticated procedure
will be required, most likely involving a handlebody theory for
manifolds with faces.

We close by commenting on a different approach to TQFTs with
corners. In the literature, for example~\cite{CY2}, extended TQFTs
are often defined for manifolds with corners in which the basic
building blocks have the shape of bigons~\cite{CY2} with only one
sort of boundary along which one can always glue. This is a
special case of our definition which is obtained if every coloured
boundary between two corners is shrunk until it disappears and
there is a single corner left that now separates two black
boundaries.

\subsubsection*{Acknowledgements}

HP would like to thank Emmanuel College for a Research Fellowship.
HP has performed a part of this work at Emmanuel College and at the
Department of Applied Mathematics and Theoretical Physics (DAMTP),
University of Cambridge, UK, and at the Max-Planck Institute for
Gravitational Physics, Potsdam, Germany. AL is grateful to Nils Baas
for introducing him to the subject of open-closed TQFT. Both authors
would like to thank John Baez and Martin Hyland for helpful
discussions and advice and Gerd Laures, Ingo Runkel, and Jonathan
Woolf for correspondence. Both authors are grateful to the European
Union Superstring Theory Network for support.

\appendix

%
\section{Symmetric monoidal categories}
%
\label{app_moncat}

In this appendix, we collect some key definitions for easier reference.

\begin{defn}
Let $(\cal{C},\otimes,\1,\alpha,\lambda,\rho)$ and
$(\cal{C}',\otimes,\1',\alpha',\lambda',\rho')$ be monoidal
categories. A \emph{monoidal functor}
$\psi\colon\cal{C}\to\cal{C}^\prime$ is a triple
$\psi=(\psi,\psi_2,\psi_0)$ consisting of,
\begin{itemize}
\item
  a functor $\psi \maps C \to C'$,
\item
  a natural isomorphism $\psi_{2}\maps\psi(X)\ten\psi(Y)\to\psi(X\ten
  Y)$, where for brevity we suppress the subscripts indicating the
  dependence of this isomorphism on $X$ and $Y$, and
\item
  an isomorphism $\psi_0 \maps \1' \to \psi(\1)$,
\end{itemize}
such that the following diagrams commute for all objects
$X,Y,Z\in\cal{C}$:
\begin{gather}
\begin{aligned}
\xymatrix@!C{
 (\psi(X) \ten \psi(Y)) \ten \psi(Z)
   \ar[r]^>>>>>>>{\psi_{2} \ten \1}
   \ar[d]_{a_{\psi(X), \psi(Y), \psi(Z)}}
& \psi(X \ten Y) \ten \psi(Z)
   \ar[r]^{\psi_{2}}
& \psi((X \ten Y) \ten Z)
   \ar[d]^{\psi(\alpha_{X,Y,Z})}   \\
 \psi(X) \ten (\psi(Y) \ten \psi(Z))
   \ar[r]^>>>>>>>{\1 \ten \psi_{2}}
& \psi(X) \ten \psi(Y \ten Z)
   \ar[r]^{\psi_{2}}
& \psi(X \ten (Y \ten Z))
 }
\end{aligned}\\
\begin{aligned}
\xymatrix{
 \1' \ten \psi(X)
   \ar[r]^{\lambda'_{\psi(X)}}
   \ar[d]_{\psi_{0} \ten \1}
&  \psi(X)  \\
 \psi(\1) \ten \psi(X)
   \ar[r]^{\psi_{2}}
&  \psi(\1 \ten X)
   \ar[u]_{\psi(\lambda_{X})}
}
\end{aligned}\\
\begin{aligned}
\xymatrix{
 \psi(X) \ten \1'
   \ar[r]^{\rho'_{\psi(X)}}
   \ar[d]_{\1 \ten \psi_{0}}
&  \psi(X)    \\
 \psi(X) \ten \psi(\1)
   \ar[r]^{\psi_{2}}
&  \psi(X \ten \1)
   \ar[u]_{\psi(\rho_X)}
}
\end{aligned}
\end{gather}
The monoidal functor is called \emph{strict} if $\psi_2$ and $\psi_0$
are identities.
\end{defn}

\begin{defn}
Let $(\cal{C},\otimes,\1,\alpha,\lambda,\rho,\tau)$ and
$(\cal{C}',\otimes,\1',\alpha',\lambda',\rho',\tau')$ be symmetric
monoidal categories. A \emph{symmetric monoidal functor}
$\psi\colon\cal{C}\to\cal{C}^\prime$ is a monoidal functor for which
the following additional diagram commutes for all $X,Y\in\cal{C}$:
\begin{equation}
\begin{aligned}
\xymatrix{
 \psi(X) \ten \psi(Y)
   \ar[r]^{\tau'_{X,Y}}
   \ar[d]_{\psi_2}
&  \psi(Y) \ten \psi(X)
    \ar[d]_{\psi_2}  \\
 \psi(X \ten Y)
   \ar[r]^{\psi(\tau)}
&  \psi(Y \ten X) }
\end{aligned}
\end{equation}
\end{defn}

\begin{defn}
Let $(\cal{C},\otimes,\1,\alpha,\lambda,\rho)$ and
$(\cal{C}',\otimes,\1',\alpha',\lambda',\rho')$ be monoidal categories
and $\psi\colon\cal{C}\to\cal{C}^\prime$ and
$\psi^\prime\colon\cal{C}\to\cal{C}^\prime$ be monoidal functors. A
\emph{monoidal natural transformation}
$\phi\colon\psi\Rightarrow\psi^\prime$ is a natural transformation
such that for all objects $X,Y$ of $\cal{C}$, the following diagrams
commute,
\begin{equation}
\begin{aligned}
\xymatrix{
  \psi(X) \ten \psi(Y) \ar[rr]^{\phi_{X} \ten \phi_Y}
  \ar[d]_{\psi_2}
  && \psi'(X) \ten \psi'(Y) \ar[d]^{\psi'_2}   \\
  \psi(X \ten Y) \ar[rr]^{\phi_{X \ten Y}}
  && \psi'(X \ten Y)
  }
\end{aligned}
\qquad\mbox{and}\qquad
\begin{aligned}
\xymatrix{
 \1' \ar[d]_{\psi_0} \ar[dr]^{\psi'_0}
 \\
 \psi(\1) \ar[r]^{\phi_{\1}} & \psi'(\1)
 }
\end{aligned}
\end{equation}
\end{defn}

\begin{defn}
Let $\cal{C}$ be a small symmetric monoidal category and let
$\cal{C}^\prime$ be an arbitrary symmetric monoidal category. We
denote by $\cat{Symm-Mon}(\cal{C},\cal{C}^\prime)$ the category of
symmetric monoidal functors $\cal{C}\to\cal{C}^\prime$ and
monoidal natural transformations between them. It is clear that
the tensor product of symmetric monoidal functors and monoidal
natural transformations defines a symmetric monoidal structure on
the category $\cat{Symm-Mon}(\cal{C},\cal{C}^\prime)$.
\end{defn}

\begin{defn}
Let $\cal{C}$ and $\cal{C}^\prime$ be monoidal categories. We say that
$\cal{C}$ and $\cal{C}^\prime$ are \emph{equivalent as monoidal
categories} if there is an equivalence of categories
$\cal{C}\simeq\cal{C}^\prime$ given by functors
$F\colon\cal{C}\to\cal{C}^\prime$ and
$G\colon\cal{C}^\prime\to\cal{C}$ and natural isomorphisms $\eta\colon
1_\cal{C}\Rightarrow GF$ and $\epsilon\colon FG\Rightarrow
1_{\cal{C}^\prime}$ such that both $F$ and $G$ are monoidal functors
and $\eta$ and $\epsilon$ are monoidal natural transformations.

If $\cal{C}$ and $\cal{C}^\prime$ are symmetric monoidal categories,
we say that they are \emph{equivalent as symmetric monoidal
categories} if in addition $F$ and $G$ are symmetric monoidal
functors.
\end{defn}

%
\end{document}